\documentclass{kybernetika_spec} 
\usepackage{amsmath,amsfonts,amssymb,amsthm}
\usepackage{color,xspace,curves,graphicx}
        \makeindex
   \DeclareMathAlphabet{\mathsfsl}{OT1}{cmss}{m}{sl}
   \DeclareMathAlphabet{\mathbfsl}{OT1}{cmr}{bx}{it}
   \newfont{\mitq}{cmmi10 at 9pt}    % p265 Kopka1
      \renewcommand{\ge}{\geqslant}\renewcommand{\le}{\leqslant}
      \renewcommand{\geq}{\geqslant}\renewcommand{\leq}{\leqslant}
      \def\R{\mathbb{R}}
      \newcommand{\sm}{\setminus}\newcommand{\pdm}{\subseteq}\newcommand{\pmn}{\emptyset}
         \newcommand{\te}{\ensuremath{\tau}\xspace}        \newcommand{\teob}{\ensuremath{\theta}\xspace}
         \newcommand{\vte}{\ensuremath{\vartheta}\xspace}  \newcommand{\bnu}{\ensuremath{{\bf0}}}
   \newcommand{\vynech}[1]{}    
\newtheorem{theorem}{Theorem}[section]
\newtheorem{lemma}[theorem]{Lemma}
\newtheorem{proposition}[theorem]{Proposition}
\newtheorem{corollary}[theorem]{Corollary}
\newtheorem{remark}[theorem]{Remark}

\newtheorem{example}[theorem]{Example}
\newtheorem{definition}[theorem]{Definition}

    \newcommand{\ttt}{\ensuremath{z}\xspace}
    \newcommand{\TT}{\ensuremath{Z}\xspace}
    \newcommand{\calT}{\ensuremath{\mathcal Z}\xspace}

    \newcommand{\cl}[1]{\ensuremath{\mathsfsl{cl}(#1)}}
    \newcommand{\ri}[1]{\ensuremath{\mathsfsl{ri}(#1)}}
    \newcommand{\s}[1]{\ensuremath{\mathsfsl{s}(#1)}}
    \newcommand{\cc}[1]{\ensuremath{\mathsfsl{cc}(#1)}}  
         
       \newcommand{\cnc}[1]{\ensuremath{\mathsfsl{cnc}(#1)}}
       \newcommand{\cns}[1]{\ensuremath{\mathsfsl{cns}(#1)}}
           \newcommand{\conv}[1]{\ensuremath{\mathsfsl{conv}(#1)}}
           \newcommand{\cone}[1]{\ensuremath{\mathsfsl{cone}(#1)}}
    
    \newcommand{\cs}[1]{\ensuremath{\mathsfsl{cs}(#1)}}

    \newcommand{\sg}[1]{\ensuremath{\mathsfsl{sgn}(#1)}}
        \newcommand{\sgs}{\ensuremath{\mathsfsl{sgn}}}
    
                 \newcommand{\exn}[1]{\ensuremath{\mathsfsl{exn}(#1)}}
    \newcommand{\dom}[1]{\ensuremath{\mathsfsl{dom}(#1)}}

       \newcommand{\PCQ}{{\sc pcq}\xspace}  \newcommand{\DCQ}{{\sc dcq}\xspace}
              \newcommand{\ME}{{\sc m{\footnotesize ax}e{\footnotesize nt}}\xspace}
       \newcommand{\efa}{\ensuremath{g^*_a}\xspace}
       \newcommand{\pya}{\ensuremath{\tilde{g}_a}\xspace}
        \newcommand{\pyah}{\ensuremath{\tilde{g}_{\nmh,a}}\xspace}
       \newcommand{\gea}{\ensuremath{\hat{g}_a}\xspace}
        \newcommand{\geah}{\ensuremath{\hat{g}_{\nmh,a}}\xspace}
       
       \renewcommand{\Theta}{\ensuremath{\varTheta}\xspace}
  \newcommand{\clF}{\ensuremath{\cl{F}}\xspace}

  \newcommand{\riF}{\ensuremath{\ri{F}}\xspace}
  \newcommand{\zat}[2]{\ensuremath{\langle{#1},{#2}\rangle}}
  \newcommand{\zlo}[2]{\ensuremath{\mbox{\large$\frac{#1}{#2}$}}}
     \newcommand{\spbox}[1]{\raisebox{-0.15ex}{#1}}
 \newcommand{\intt}{\ensuremath{\spbox{\Large$\int$}\:}}
  \newcommand{\inttt}[1]{\ensuremath{\mbox{\large$\int$}_{\!\!\!#1}}}
 
 \newcommand{\norm}[1]{\ensuremath{|\!|#1|\!|}\xspace}

 \newcommand{\calG}{\ensuremath{\mathcal G}\xspace}
 \newcommand{\calH}{\ensuremath{\mathcal H}\xspace}
 \newcommand{\moma}{\ensuremath{\varphi}\xspace}
 \newcommand{\ga}{\ensuremath{\gamma}\xspace}     \newcommand{\coga}{\gamma^*}
 \newcommand{\gat}{\ensuremath{[\gamma t]}\xspace}     \newcommand{\cogat}{\gat^*}
 \newcommand{\nm}{\ensuremath{\beta}\xspace}     \newcommand{\conm}{\ensuremath{\nm^*\xspace}\xspace}
 \newcommand{\nmh}{\ensuremath{[\beta h]}\xspace}\newcommand{\conmh}{\ensuremath{\nmh^*\xspace}\xspace}
 
 \newcommand{\fuc}{\ensuremath{\mathsfsl H}}
 \newcommand{\vafu}{{\mathsfsl J}}   \newcommand{\revafu}{{\mathsfsl I}}
 \newcommand{\dufu}{{\mathsfsl K}}
 \newcommand{\vare}{\varepsilon}
 \newcommand{\asli}{Z_{\nm,\,\mathrm{al}}}
 \newcommand{\cofi}{Z_{\nm,\,\mathrm{cf}}}
 \newcommand{\trn}{{\scriptstyle\text{\mitq\symbol{1}}}_{\ga}}
 \newcommand{\trnga}[1]{{\scriptstyle\text{\mitq\symbol{1}}}_{\ga}(#1)}
 \newcommand{\trnarg}[2]{{\scriptstyle\text{\mitq\symbol{1}}}_{#1}(#2)}
 \newcommand{\trndu}[1]{{\scriptstyle\text{\mitq\symbol{1}}}_{\coga}(#1)}
 \newcommand{\lad}{{\scriptstyle\text{\mitq\symbol{7}}}_{\ga}}
 \newcommand{\ladga}[1]{{\scriptstyle\text{\mitq\symbol{7}}}_{\ga}(#1)}
 \newcommand{\ladarg}[2]{{\scriptstyle\text{\mitq\symbol{7}}}_{#1}(#2)}
 \newcommand{\trnnin}{{\scriptstyle\text{\mitq\symbol{1}}}_{\nm}}

 \newcommand{\ladnin}{{\scriptstyle\text{\mitq\symbol{7}}}_{\nm}}
 \newcommand{\egy}{\mbox{${\mathit1}\hspace*{-2.4mm}{\mathit1}\!$}}
 \newcommand{\gadiv}[2]{\ensuremath{{\mathsfsl D}_\ga(#1,#2)}}
 \newcommand{\Cot}{\mathsfsl{C}}    \newcommand{\Dore}{\mathsfsl{D}}
 \newcommand{\Bre}[1]{\ensuremath{{\mathsfsl  B}_{\nm}(#1)}}
 \newcommand{\Brega}[1]{\ensuremath{{\mathsfsl  B}_{\ga}(#1)}}
         \newcommand{\Bres}{\ensuremath{{\mathsfsl B}_{\nm}}}
         \newcommand{\Breh}[1]{\ensuremath{{\mathsfsl B}_{\nmh}(#1)}}
 
   \newcommand{\BreZ}[2]{\ensuremath{{\mathsfsl  B}_{#2,\nm}(#1)}}
 \newcommand{\Crc}[1]{\ensuremath{\omega_{#1}}}
  \newcommand{\FF}{\ensuremath{\mathcal{F}}\xspace}
  \newcommand{\Fa}[1]{\ensuremath{\boldsymbol{F}_{\!#1}}\xspace}
     
      \newcommand{\cnmoma}[1]{\ensuremath{\mathsfsl{cn}_{\moma}(#1)}}
       \newcommand{\cnmomapl}[1]{\ensuremath{\mathsfsl{cn}_{\moma}^+(#1)}}
 \newcommand{\cgd}{u}

\renewcommand {\thesection}{\arabic{section}}
\renewcommand {\thesubsection}{\arabic{section}.{\small\Alph{subsection}}}
\begin{document}
\pagestyle{myheadings}

\title{Generalized minimizers\newline of convex integral functionals,\newline Bregman distance,
Pythagorean identities}

\author{Imre Csisz\'ar and Franti\v{s}ek Mat\'{u}\v{s}}

\contact{Imre}{Csisz\'ar}{A.\ R\'{e}nyi Institute of Mathematics, Hungarian Academy of Sciences,
                          H-1364 Budapest, P.O.Box 127. Hungary.}{csiszar.imre@renyi.mta.hu}
\contact{Franti\v{s}ek}{Mat\'{u}\v{s}}{Institute of Information Theory and Automation --
                          Academy of Sciences of the Czech Republic, Pod Vod\'arenskou v\v{e}\v{z}\'{\i}~4,
                          182\,08 Praha 8. Czech Republic.}{matus@utia.cas.cz}

\markboth{I. Csisz\'ar and F. Mat\'{u}\v{s}} {Generalized minimizers of convex integral functionals}

\maketitle

\indent {\it This paper is dedicated to the memory of Igor Vajda (1942--2010)}

\bigskip

\begin{abstract}
    Integral functionals based on convex normal integrands are minimized
    subject to finitely many moment constraints. The integrands are finite
    on the positive and infinite on the negative numbers, strictly convex
    but not necessarily differentiable. The minimization is viewed as a primal
    problem and studied together with a dual one in the framework of convex
    duality. The effective domain of the value function is described by a conic
    core, a modification of the earlier concept of convex core. Minimizers and
    generalized minimizers are explicitly constructed from solutions of modified
    dual problems, not assuming the primal constraint qualification. A~generalized
    Pythagorean identity is presented using Bregman distance and a correction term
    for lack of essential smoothness in integrands. Results are applied to
    minimization of Bregman distances. Existence of a generalized dual solution
    is established whenever the dual value is finite, assuming the dual constraint
    qualification. Examples of `irregular' situations are included, pointing to
    the limitations of generality of certain key results.
\end{abstract}

\keywords{maximum entropy, moment constraint, generalized
    primal/dual solutions, normal integrand, minimizing sequence,
    convex duality, Bregman projection, conic core,
    generalized exponential family, inference principles}

\classification{94A17, 49J53, 49K30, 62B10, 65K10, 90C46}

\begin{center}
\begin{tabular}{rlr}
\multicolumn{3}{c}{{\bf\sc Contents}}\\
\multicolumn{3}{c}{}\\
1. & Introduction & \pageref{S:intro}\\
2. & Preliminaries & \pageref{S:prelim}\\
3. & Preliminaries on the dual problem & \pageref{S:dualproblem}\\
4. & The constraint qualifications & \pageref{S:roleCQ}\\
5. & Conic cores & \pageref{S:cnc}\\
6. & The effective domain of the value function &
    \pageref{S:dom}\\
7. & Dispensing with the \PCQ in the primal problem &
    \pageref{S:main}\\
8. & Bregman projections & \pageref{S:BREGproj}\\
9. & Generalized solutions of the dual problem &
    \pageref{S:dual}\\
\end{tabular}
\end{center}    \normalsize

\begin{center}
\begin{tabular}{rlr}
10. & Examples & \pageref{S:examples}\\
11. & Relation of this work to previous ones &
    \pageref{S:relation}\\
    \multicolumn{2}{l}{Appendix A \ Integral representation} & \pageref{S:intrep}\\
    \multicolumn{2}{l}{Appendix B \ Restricted value function} & \pageref{S:revafu}\\
    \multicolumn{2}{l}{Appendix C \ \ga-divergences} & \pageref{S:gadiv}\\
    \multicolumn{2}{l}{Acknowledgement} & \pageref{acknow}\\
    \multicolumn{2}{l}{References} & \pageref{references}\\
    \multicolumn{2}{l}{Index} & \pageref{index}
\end{tabular}
\end{center}    \normalsize

\section{Introduction}\label{S:intro}

 \paragraph{1.A.} Let $\mu$ be a $\sigma$-finite measure on a measurable space
 $(\TT,\calT)$\index{$(\TT,\calT)$~~underlying measurable space}\label{underlying measurable space}
 and $\moma\colon\TT\to\R^d$ a \calT-measurable vector-valued
 function referred\index{$\moma$~~moment mapping}\label{moment mapping}
 to as the \emph{moment mapping}. The linear space of the \calT-measurable
 functions $g\colon\TT\to\R$ with $\mu$-integrable $\moma g$ is denoted by
 $\calG$, and\index{$\triangleq$~~equal by definition}\label{equal by definition}
 \[
    \calG_a\triangleq\big\{g\in\calG\colon\mbox{$\int_\TT$}\,\moma g\,{\mathrm d}\mu=a\big\}\,,\qquad a\in\R^d\,.
 \]
 Here, $a$ is the \emph{moment vector} of $g\in\calG_a$\index{$\calG$~~class of functions with a moment}\label{class of functions with a moment}
 \index{$\calG_a$~~class of functions with the moment $a$}\label{class of functions with the moment a}
 while the functions $g\notin\calG$ have no moment vectors.
 The set of nonnegative functions in $\calG$/$\calG_a$ is denoted by
 $\calG^+$/$\calG_a^+$.\index{$\calG^+$~~class of nonnegative functions in $\calG$}\label{class of nonnegative functions in}

 This work studies the minimization of \emph{integral functionals}
 of the form\index{$\fuc_{\nm}$~~integral functional}\label{integral functional}
 \begin{equation}\label{E:fuc}
   \fuc_{\nm}(g)\triangleq \inttt{\TT}\:\nm(\ttt,g(\ttt))\;\mu({\mathrm d}\ttt)
 \end{equation}
 subject to $g\in\calG_a$. Here, $\nm\colon\TT\times\R\to(-\infty,+\infty]$
 is a \emph{normal convex integrand}~\cite[Chapter~14]{RW} such that for
 $\ttt\in\TT$ the function $t\mapsto\nm(\ttt,t)$ is finite and strictly convex
 when $t>0$ and equals $+\infty$ when $t<0$. The positive and negative parts
 of the integral in~\eqref{E:fuc} may be both infinite, in which case the
 integral is taken to be $+\infty$ by convention. A necessary condition for
 $\fuc_{\nm}(g)<+\infty$ is the nonnegativity of~$g$, thus the minimization
 of~$\fuc_\nm$ is actually over the family
 $\calG_a^+$.\index{$\nm$~~integrand with $\nm(\ttt,\cdot)\in\varGamma$}\label{integrand with nm}

% INFERENCE

\paragraph{1.B.} Minimization problems of this kind emerge across various
scientific disciplines, notably
 in \emph{inference}. When $g$ is an unknown nonnegative function on $\TT$ whose moment
 vector $\int_\TT\,\moma g\,{\mathrm d}\mu$ can be measured in an experiment providing a vector
 $a$, typical \emph{inference principles} call for adopting, as `best guess' of~$g$,
 a minimizer of $\fuc_{\nm}$ over $\calG_a^+$, for a specific choice of $\nm$. The unknown
 $g$ may be a probability density, or its integral may be known otherwise, in which case
 one coordinate function of the moment mapping $\moma$ is taken to be identically~$1$. Most
 often \emph{autonomous} integrands are used, which means that $\nm$ does not depend on the
 first coordinate $\ttt\in\TT$. Typical choices are $t\ln t$ or $-\ln t$ or~$t^2$ giving
 $\fuc_{\nm}(g)$ equal to the negative Shannon or Burg entropy\footnote{Here,
 `entropy' is understood in a wide sense. Shannon entropy in the strict sense
 refers to the case when $\mu$ is the counting measure on a finite or countable
 set and $g$ is a probability mass function.} or squared $L^2$-norm of
 $g\geq0$.

 If a `prior guess' $h$ for $g$ is available, that would be adopted before
 the measurement, related inference principles suggest  to take as `best guess'
 after the measurement the minimizer of some `distance' of $g$ from~$h$
 subject to $g\in\calG_a^+$. Two kinds of non-metric distance often used in this
 context are \emph{Bregman distances}, see eq.~\eqref{E:Brebe}, and\index{$\ga$-divergence}\label{ga divergence}
 \emph{$\ga$-divergences}, see eq.~\eqref{E:gadiv}. The most familiar
 is the information ($I$-) divergence, also called Kullback--Leibler distance
 or relative entropy, that belongs to both families. For $h$ fixed, both kinds
 of distance are nonnegative integral functionals in $g$ of the form~\eqref{E:fuc},
 with non-autonomous integrands. The minimization of $\ga$-divergences can be
 easily reduced to that of integral functionals with autonomous integrands,
 see Appendix~\ref{S:gadiv}, but this is not possible for Bregman distances
 except in special cases. It will become evident below that Bregman distances
 inevitably enter the minimization of $\fuc_{\nm}$ over $\calG_a^+$, even
 in the autonomous case.

 Another common approach to inference problems as above is to specify
 a priori a family of functions $f_\vte$ parameterized by some $\vte$
 and search in that family a function whose moment vector equals the
 experimentally measured vector $a$, thus solve the equation
 $\int_\TT\,\moma f_\vte\,{\mathrm d}\mu=a$ in the parameter $\vte$. There
 is a close relationship between this approach and the one based
 on the minimization of $\fuc_{\nm}$ given moment constraints.
 Indeed, the latter will suggest to use the parametric family
 defined after eq.~\eqref{E:fte}. If some function in that family has
 moment vector equal to~$a$ then this function minimizes
 $\fuc_{\nm}$ on~$\calG_a^+$. Even if no such function exists,
 it is usually possible to specify a `best' function in the family,
 which is also a `generalized solution' of the minimization problem.

% DUALITY

 \paragraph{1.C.} The minimization of $\fuc_{\nm}$ over $\calG_a$, or equivalently
 over $\calG_a^+$, is approached here by convex duality theory,
 as in~\cite{Bo.Le.dua,Bo.Le}. A strategy is to introduce
 the \emph{value function} $\vafu_{\nm}$ by\index{$\vafu_{\nm}$~~value function}\label{nm value function}
 \begin{equation}\label{E:priva}
    \vafu_{\nm}(a)\triangleq{\inf}_{g\in\calG_a^+}\: \fuc_{\nm}(g)\,,
    \qquad  a\in\R^d\,,
 \end{equation}
 and to study its conjugate and biconjugate.\index{$^*$~~conjugate of convex function}\label{conjugate of convex function}
 The value function ranges in $[-\infty,+\infty]$ and is convex. The case when it is identically
 $+\infty$ is often excluded, writing $\vafu_{\nm}\not\equiv+\infty$, but it is sometimes
 not straightforward to recognize. Usually, the value function is \emph{proper}, thus not
 identically $+\infty$ and never equal to $-\infty$. No general description of the effective
 domain $\dom{\vafu_{\nm}}$ of the value function,\index{$\mathsfsl{dom}$~~effective domain}\label{effective domain} thus the set
 of $a\in\R^d$ with $\vafu_{\nm}(a) <+\infty$, seems to be available in literature. This domain
 is contained in the set of the moment vectors\index{$\cnmoma{\mu}$~~$\moma$-cone of $\mu$}\label{cone of}
 $\int_{\TT}\:\moma g \,{\mathrm d}\mu$ of the functions $g\in\calG^+$, that is called
 here the \emph{$\moma$-cone $\cnmoma{\mu}$ of $\mu$}. Theorem~\ref{T:domain} describes
 $\dom{\vafu_{\nm}}$ in terms of faces of  $\cnmoma{\mu}$. A crucial point is to represent
 the $\moma$-cone via a new concept of \emph{conic core} for Borel measures on $\R^d$,
 introduced in Section~\ref{S:cnc} similarly to the convex cores in~\cite{Csi.Ma.cc}.

 The minimization in~\eqref{E:priva} is the \emph{primal problem} and the infimum
 $\vafu_{\nm}(a)$ is the \emph{primal value} for~$a$.\index{$\vafu_{\nm}(a)$~~primal value for $a$}\label{primal value for a}
 The value is attained if a minimizer exists. Since $\nm$ is strictly convex,
 if $\vafu_{\nm}(a)$ is finite then such a minimizer is unique\footnote{
         in the sense that any two minimizers are $\mu$-a.e.\ equal.
         As a rule, equality of functions is understood $\mu$-a.e., unless
         $\ttt\in\TT$ is included in the notation.}
 and it is referred to as the \emph{primal solution} $g_a$ for $a$.\index{$g_a$~~primal solution for $a$}\label{primal solution for a}
 A first goal is to recognize whether the primal value is finite, then whether it is attained
 in which case a construction of the primal solution is desirable. A~second goal is
 to understand the behavior of minimizing sequences $g_n$ in $\calG_a^+$ for which
 $\fuc_{\nm}(g_n)$ converges to the primal value $\vafu_{\nm}(a)$. When all minimizing
 sequences converge to a common limit locally in measure then the limit function
 will\index{$\rightsquigarrow$~~local convergence in measure}\label{local convergence in measure} be called the \emph{generalized
 primal solution} and denoted by~$\gea$.\index{$\gea$~~generalized primal solution for $a$}\label{generalized primal solution for a}
 This convergence, denoted by $g_n\rightsquigarrow \gea$, means that
 $\mu(Y\cap\{|g_n-\gea|>\vare\})\to0$ for every $Y\in\calT$ of finite
 $\mu$-measure and every $\vare>0$. The fact justifying the terminology that each primal solution
 is also a generalized primal solution is discussed in Subsection~1.E.
 after eq.~\eqref{E:Pyth}, see also Corollary~\ref{C:uvid}.

 The \emph{convex conjugate} $\vafu_{\nm}^*$ of the value function
 is defined by \index{$\zat{\cdot}{\cdot}$~~inner product}\label{inner product}\label{conjugate of convex function 2}
 \[
    \vafu_{\nm}^*(\vte)\triangleq
            {\sup}_{a\in\R^d}\:\big[\zat{\vte}{a}-\vafu_{\nm}(a)\big]\,,
                    \qquad\vte\in\R^d\,,
 \]
 where $\zat{\cdot}{\cdot}$ is the scalar product on $\R^d$.
 The conjugate $\conm$ of $\nm$,
 \[
    \conm(\ttt,r)\triangleq{\sup}_{t\in\R}
        \:\big[\, r t-\nm(\ttt,t)\,\big]\,,
          \qquad\ttt\in\TT\,,\;r\in\R\,,
 \]
 is a convex normal integrand, giving rise to the integral functional
 $\fuc_{\conm}$ and the convex function $\dufu_{\nm}$ given by
 \[
    \dufu_{\nm}(\vte)\triangleq\inttt{\TT}\:
        \conm\big(\ttt,\zat{\vte}{\moma(\ttt)}\big)\:\mu({\mathrm d}\ttt)
                    =\fuc_{\conm}(\zat{\vte}{\moma})\,, \qquad\vte\in\R^d\,.
 \]
 The following key fact is referred to as the \emph{integral representation}
 of $\vafu_{\nm}^*$. Its proof, building on~\cite{Rock.cofu,Rock.cofu.dual,RW},
 is presented in Appendix A.

\begin{theorem}\label{T:repre}
    If $\vafu_{\nm}\not\equiv+\infty$ then $\vafu_{\nm}^*=\dufu_{\nm}$.
\end{theorem}

 The convex conjugate $\vafu_{\nm}^*$ is proper if and only if $\vafu_{\nm}$
 is proper \cite[Theorem 12.2]{Rock}, which takes place if and only if
 $\dom{\vafu_{\nm}}$ and $\dom{\vafu_{\nm}^*}$ are both nonempty. When
 $\dom{\vafu_{\nm}}=\pmn$, thus Theorem~\ref{T:repre} does not apply,
 $\dufu_{\nm}$ may differ from $\vafu_{\nm}^*\equiv -\infty$ and may be
 a proper convex function, see Example~\ref{Ex:vafu.dual}. A sufficient
 condition for $\vafu_{\nm}$ to be proper is the finiteness of $\dufu_{\nm}$
 on an open set, see Corollary~\ref{C:dufuint}, a new result below.

 The \emph{biconjugate}\index{biconjugate}\label{biconjugate} of $\vafu_{\nm}$ is obtained
 by conjugating $\vafu_{\nm}^*$,
 \[
    \vafu_{\nm}^{**}(a)\triangleq
            {\sup}_{\vte\in\R^d}\:\big[\zat{\vte}{a}
                        - \vafu_{\nm}^*(\vte)\big]\,,
                                \qquad a\in\R^d\,.
 \]
 If $\vafu_{\nm}\not\equiv+\infty$ then $\vafu_{\nm}^{**}=\dufu_{\nm}^*$,
 by Theorem~\ref{T:repre}.

 The maximization in the conjugation of $\dufu_\nm$
 \begin{equation}\label{E:duva}
    \dufu_{\nm}^*(a)={\sup}_{\vte\in\R^d}\:
            \Big[\zat{\vte}{a}-\inttt{\TT}\:
                \conm\big(\ttt,\zat{\vte}{\moma(\ttt)}\big)\:\mu({\mathrm d}\ttt)\Big]\,,
                    \qquad a\in\R^d\,,
 \end{equation}
 is called the \emph{dual problem} for $a$, also when $\vafu_{\nm}\equiv+\infty$, thus
 when Theorem~\ref{T:repre} does not apply. The supremum $\dufu_{\nm}^*(a)$ in~\eqref{E:duva}
 is\index{$\dufu^*_{\nm}(a)$~~dual value for $a$}\label{dual value for} a \emph{dual value}. If it is finite and attained,
 each maximizer is a \emph{dual solution}. The latter situation is also referred to as existence of
 Lagrange multipliers, see~\cite{Bo.Le}. There is an intimate relationship between the primal and
 dual problems discussed in detail below. The primal value always dominates the dual one,
 see Lemma~\ref{L:duleqpri}. Their distance is the \emph{duality gap}. If the gap is zero,
 thus the primal and dual values coincide, the dual problem provides valuable
 information on the primal one. What makes the strategy effective is
 that the dual problem is finite dimensional and unconstrained.

% CONSTRAINTS

\paragraph{1.D.} Standard results are typically proved under the pair of
conditions\index{\PCQ~~primal constraint qualification}\label{primal constraint qualification}
 \begin{equation}\tag{{\sc pcq}}
    \text{$\vafu_{\nm}$ is proper and $a\in\ri{\dom{\vafu_{\nm}}}$}
 \end{equation}
 referred to jointly as the \emph{primal constraint qualification}.
 Here, $\mathsfsl{ri}$ stands for the relative interior. A convex
 function that takes the value $-\infty$ somewhere, does so
 everywhere in the relative interior of its effective domain, and thus
 the \PCQ can be equivalently stated replacing the first condition by
 $\vafu_{\nm}(a)>-\infty$. By Remark~\ref{R:Slater}, the second condition
 in \PCQ is equivalent to the existence of a positive function $g$ in
 $\calG_a$. Under the \PCQ for $a$, the duality gap is zero,
 $\vafu_{\nm}(a)=\dufu_{\nm}^*(a)$ \cite[Theorems 7.4 and 12.2]{Rock}.

 A special role will be played by the set $\Theta_{\nm}$ of those
 $\vte\in\dom{\dufu_{\nm}}$ for which the function $r\mapsto\conm(\ttt,r)$
 is finite in a neighborhood of $\zat{\vte}{\moma(\ttt)}$ for $\mu$-a.a.\
 $\ttt\in\TT$. This set is convex but possibly empty.
 The assumption\index{$\Theta_{\nm}$~~a subset of \dom{\dufu_{\nm}}}\label{a subset of dom}
 \begin{equation}\tag{\mbox{\sc dcq}}
    \text{$\Theta_{\nm}$ is nonempty}
 \end{equation}
 is referred to as the \emph{dual constraint qualification}\index{\DCQ~~dual constraint qualification}\label{dual constraint qualification}
 (\DCQ). For sufficient conditions of its validity see Remark~\ref{R:momas}. If the \DCQ holds,
 maximization in the dual problem \eqref{E:duva} can be restricted to $\Theta_{\nm}$ without
 changing the dual value or loosing a dual solution, see Lemmas~\ref{L:Theta}
 and~\ref{L:dualsol}.

 Computation of directional derivatives of $\dufu_{\nm}$ features
 the following functions $f_\vte$ of $\ttt\in\TT$,
 \begin{equation}\label{E:fte}
  f_\vte(\ttt)\triangleq\begin{cases}\displaystyle(\conm)'(\ttt,\zat{\vte}{\moma(\ttt)})\,,
              \quad&\displaystyle\text{if $\conm(\ttt,\cdot)$ is differentiable at
                    $\zat{\vte}{\moma(\ttt)}$,}\\
                                        0\,,&\text{otherwise.}
                           \end{cases}
 \end{equation}
 The family $\FF_{\nm}\triangleq\{f_\vte\colon\vte\in\Theta_{\nm}\}$
 will play a similar role as exponential families do in the case
 of the negative Shannon entropy functional \cite{BaNi,Chen}.
 \index{$\FF_{\nm}$~~family of functions $f_\vte$}\label{family of functions}

% UNDER PCQ

 \paragraph{1.E.} Let the \PCQ hold for~$a\in\R^d$. Then, the primal and dual values for
 $a$ are finite, coincide, $\vafu_{\nm}(a)=\dufu_{\nm}^{*}(a)$, and a dual
 solution $\vte\in\dom{\dufu_{\nm}}$ exists, by Lemma~\ref{L:pri=du}.
 If the \DCQ fails then no primal solution exists and the generalized
 primal solution does not exist either, see Lemma~\ref{L:ainri} and
 Theorem~\ref{T:geprisol}. Otherwise, if $\Theta_{\nm}\neq\pmn$, each
 dual solution $\vte$ belongs to~$\Theta_{\nm}$ and gives rise to the
 same function $f_\vte$, by Corollary~\ref{C:nondep}. This unique function
 from $\FF_{\nm}$ is called here the \emph{effective dual solution}
 for~$a$\index{$\efa$~~effective dual solution}\label{effective dual solution}
 and is denoted by~$\efa$, see Remark~\ref{R:efduso}. The primal solution
 $g_a$ exists if and only if $\int_\TT\,\moma\efa\,{\mathrm d}\mu$ exists and
 equals~$a$, in which case $g_a=\efa$, see Lemma~\ref{L:ainri}.
 Alternatively, by the same lemma, the primal solution $g_a$ exists
 if and only if $\FF_{\nm}$ intersects $\calG_a$, thus the equation
 $\int_\TT\,\moma f_\vte\,{\mathrm d}\mu=a$ has a solution $\vte\in\Theta_{\nm}$.
 In this case, $\FF_{\nm}\cap \calG_a=\{g_a\}$. Subject to the \PCQ
 and \DCQ, these conditions are always satisfied if $\dufu_{\nm}$
 is essentially  smooth, then $g_a$ exists and equals $\efa$, see
 Corollary~\ref{C:prisolclassic}. Otherwise, $g_a$ may not exist,
 for~$\efa$ need not have moment vector, or its moment vector may
 differ from~$a$, see Examples~\ref{Ex:momefa} and \ref{Ex:Burg}.
 Under the \PCQ for~$a$, however, the \DCQ is necessary and sufficient
 for the existence of the generalized primal solution $\gea$, which
 then coincides with the effective dual solution $\efa$, see
 Theorem~\ref{T:geprisol}, a new result.

 The main results of this paper include extensions of the above assertions to
 the cases when the \PCQ is relaxed to the finiteness of~$\vafu_{\nm}(a)$, see
 Section~\ref{S:main}. These are relevant when the effective domain of $\vafu_{\nm}$
 includes a nontrivial relative boundary. Depending on the position
 of~$a$ in the convex set $\dom{\vafu_{\nm}}$, the dual problem is modified,
 restricting the integration to a subset of $\TT$ that corresponds to a face
 of the convex cone $\cnmoma{\mu}$. In the modified problem a solution
 exists and the above assertions have appropriate reformulations, see
 Theorem~\ref{T:primal.dual}. In particular, a primal solution exists
 if and only if an extension of the family $\FF_{\nm}$ intersects $\calG_a$,
 see Corollary~\ref{C:primal.dual}. This resolves existence of the primal
 and generalized primal solutions for $a\in\R^d$ and their construction
 without the \PCQ, whenever $\vafu_{\nm}(a)$ is finite,
 even if $\vafu_{\nm}(b)=-\infty$ for some $b\neq a$.

 Another main result is the \emph{generalized Pythagorean
 identity},\index{generalized Pythagorean identity}\label{generalized Pythagorean identity}
 see~Theorem~\ref{T:main}, asserting that for any $a\in\R^d$ with
 $\vafu_{\nm}(a)$ finite there exists a unique function $\pya$ such that
 \begin{equation}\label{E:Pyth}
       \fuc_{\nm}(g)=\vafu_{\nm}(a)+\Bres(g,\pya)+\Cot_\nm(g)\,,
            \qquad g\in\calG_{a}^+\,,
 \end{equation}
 under a condition not stronger than the \DCQ. Under the \PCQ for $a$
 and \DCQ, the function $\pya$ equals $g_a^*$ while in general $\pya$
 is constructed as the effective dual solution of a modified dual problem.
 In \eqref{E:Pyth}, $\Bres$ denotes Bregman distance defined by
 eq.~\eqref{E:Brebe} and $\Cot_\nm$ is a nonnegative correction
 functional, defined in special cases by eq.~\eqref{E:Corr} and in
 general by eq.~\eqref{E:Corr2}. The idea to involve a correction is
 new even under the \PCQ and \DCQ. When \nm is essentially smooth then
 $\Cot_\nm$ is identically zero and \eqref{E:Pyth} without the correction
 becomes a \emph{Pythagorean identity}. In general, omitting $\Cot_\nm$
 in \eqref{E:Pyth} a \emph{Pythagorean inequality} arises. The inequality
 allows for the conclusions that the generalized primal solution $\gea$
 exists and equals $\pya$, and that if the primal solution $g_a$ exists
 then $g_a=\pya=\gea$, using Corollary~\ref{C:inmeas}.

 In absence of the \PCQ, generalized solutions are introduced also
 for the dual problem, and their existence is proved under general conditions,
 see Theorem~\ref{T:GMLE}. The generalized primal and dual solutions coincide
 if the duality gap is zero. In general, their Bregman distance is not larger
 that the duality gap, see Remark~\ref{R:ga=ha}.

% OUTLINE

 \paragraph{1.F.} This work is organized as follows. Section~\ref{S:prelim} collects
 definitions, technicalities, auxiliary lemmas, and presents general results
 on the normal integrands and Bregman distances. In Section~\ref{S:dualproblem},
 the function $\dufu_\nm$ is studied, its directional derivatives computed,
 and a new sufficient condition for $\vafu_{\nm}\not\equiv +\infty$ is
 presented in terms of this function.

 Section~\ref{S:roleCQ} summarizes results about the primal a dual
 problems, mostly familiar in the case of autonomous and essentially
 smooth integrands. These results cover the case when the
 \PCQ holds, but some describe also the more general situation when
 the primal and dual values coincide and a dual solution exists.
 The concepts of effective dual solutions and of the correction functional
 are introduced and a first restricted version of the generalized Pythagorean
 identity is elaborated, which appears new already in this form. Another
 new result relates the existence of generalized primal solutions
 to the~\DCQ.

 Conic cores are introduced and studied in Section~\ref{S:cnc}.
 In Section~\ref{S:dom} a geometric description of the effective domain
 of the value function is given via the $\moma$-cone of $\mu$.
 Section~\ref{S:main} formulates the main results on the primal problem
 without the \PCQ, including the generalized Pythagorean identity. The
 main results are specialized to the problem of Bregman projections
 in Section~\ref{S:BREGproj}, and general Pythagorean identities are
 also treated there. Section~\ref{S:dual} is devoted to the dual problem,
 its main result is a theorem on existence of generalized dual solutions.
 All examples are collected in Section~\ref{S:examples}. The relations
 of this work to previous ones are discussed in~Section~\ref{S:relation}.

 Appendix~A presents a proof of Theorem~\ref{T:repre}, extending a standard
 result about the interchange of integration and minimization. Appendix~B
 describes how the usual approach to Shannon entropy maximization is embedded
 into the framework. Appendix~C addresses $\ga$-divergences, and presents
 a lemma that would admit to restrict attention to finite measures $\mu$
 throughout this paper.

%22222222222222222222222222222222222222222222222222222222222222222222222222222
\section{Preliminaries}\label{S:prelim}

 The terminology and notation of \cite{Rock} are mostly adopted. If
 $C\pdm\R^d$ then $\cl{C}$ is the closure and $\ri{C}$ the relative
 interior of $C$, thus the interior in the topology of the affine
 hull of~$C$.\index{$\cl{F}$~~closure of $F$}\label{closure of}\index{$\ri{F}$~~relative interior of $F$}\label{relative interior of}
 Subsets of $\TT$ on which certain relations hold are
 denoted briefly by these relations in the curly brackets. For example,
 the level set $\{\ttt\in\TT\colon g(\ttt)>t\}$ of a function
 $g\colon\TT\to\R$ is denoted by $\{g>t\}$. Shorthand notations for
 $\mu$-almost everywhere are $\mu$-a.e.\ or $[\mu]$.\index{$[\mu]$ ~~$\mu$-almost everywhere}\label{almost everywhere}
 The function $\sgs\colon\R\mapsto\{+,-\}$ assigns $+$ to the \emph{nonnegative}
 and $-$ to the negative numbers.\index{$\sg{r}$~~sign of $r$, equals $+$ if $r=0$}\label{sign of r}

% FUNCTION GAMMA

 \paragraph{2.A.} Let $\varGamma$ denote the family of functions\index{$\varGamma$~~class of convex functions $\ga$}\label{class of convex functions}
 $\ga\colon\R\to(-\infty,+\infty]$\index{$\ga$~~strictly convex lsc function in $\varGamma$}\label{strictly convex lsc function in}
 that are finite and strictly convex for $t>0$, equal to $+\infty$ for $t<0$, and
 satisfy $\ga(0)=\lim_{t\downarrow0}\ga(t)$. In terms of~\cite{Rock}, $\ga$ is proper
 and closed, thus lower semi\-continuous (lsc).\index{lower semicontinuous, lsc}\label{lower semicontinuous, lsc}
 The effective domain \dom{\ga} equals $(0,+\infty)$ or $[0,+\infty)$.
 The left/right derivatives of \ga at $t>0$ are
 finite, $\ga'_-(t)\leq\ga'_+(t)$, and both $\ga'_-$ and $\ga'_+$ increase. Let
 $\ga'_-(0)\triangleq-\infty$ and $\ga'_+(0)\triangleq\lim_{t\downarrow0}\ga'_+(t)$,
 which is the standard right derivative at $0$ if $\ga(0)<+\infty$. Further, let
 $\ga(+\infty)\triangleq\lim_{t\uparrow+\infty}\ga(t)$ and
 $\ga'(+\infty)\triangleq\lim_{t\uparrow+\infty}\ga'_+(t)$.
 If $\ga'(+\infty)=+\infty$ then $\ga$ is called \emph{cofinite}.\index{cofinite function}\label{cofinite function}
 Otherwise, the function $t\mapsto t\ga'(+\infty)-\ga(t)$ is increasing.
 If it has a finite limit as $t\uparrow+\infty$ then $\ga$ is called
 \emph{asymptotically linear}.\index{asymptotically linear function}\label{asymptotically linear function}

 The convex conjugate $\coga$ of $\ga\in\varGamma$ is given by
 $\coga(r)={\sup}_{t>0}\,[r\,t-\ga(t)]$, $ r\in\R$. It is finite
 and non\-decreasing in the interval $(-\infty,\ga'(+\infty))$,
 and $\coga(r)=+\infty$ for $r>\ga'(+\infty)$. When $\ga$ is not cofinite
 then $\coga(\ga'(+\infty))=\lim_{t\uparrow+\infty}\,[t\ga'(+\infty)-\ga(t)]$,
 thus $\coga$ is finite at $\ga'(+\infty)$ if and only if $\ga$
 is asymptotically linear. If $r\downarrow-\infty$ then
 $\coga(r)\downarrow-\ga(0)$ where the limit is denoted also
 by $\coga(-\infty)$. If $\ga(0)$ is finite then
 $\coga(r)=-\ga(0)$ for $r\leq\ga'_+(0)$. The strict convexity
 of $\ga$ implies that $\coga$ is \emph{essentially smooth}
 \cite[Theorem~26.3]{Rock}, thus $\coga$ is differentiable in
 $(-\infty,\ga'(+\infty))$ and if $\ga$ is not cofinite then
 $(\coga)'(r)\uparrow+\infty$ as $r\uparrow\ga'(+\infty)$. The
 latter holds also when $\ga$ is cofinite.

 Let $\cgd$ denote the function defined for $r<\ga'(+\infty)$
 by $\cgd(r)=(\coga)'(r)$. The following lemma contains an elementary
 reformulation of the fact that the subgradient mappings of $\ga$
 and $\coga$ are mutually inverse \cite[Corollary~23.5.1]{Rock}.

\begin{lemma}\label{L:gamma0}
    Let $\ga\in\varGamma$.  For $t\geq0$ and $r\in\R$, if $\ga'_-(t)\leq r\leq \ga'_+(t)$
    then $r<\ga'(+\infty)$, $\cgd(r)=t$ and $\coga(r)=tr-\ga(t)$.
\end{lemma}

 The function $\cgd$ is nondecreasing on $(-\infty,\ga'(+\infty))$.
 It is strictly increasing if and only if $\coga$ is strictly
 convex which is equivalent to the essential smoothness of $\ga$, thus the
 differentiability of $\ga$ in $(0,+\infty)$ together with
 $\ga'_+(0)=-\infty$. Further, $\cgd(r)\downarrow0$ as $r\downarrow-\infty$,
 where the limit $0$ is interpreted as $(\coga)'(-\infty)$,
 and $\cgd(r)\uparrow+\infty$ as $r\uparrow\ga'(+\infty)$.
 The function $\cgd$ vanishes on the
 interval $(-\infty,\ga'_+(0)]$.

  Another reformulation of Lemma~\ref{L:gamma0}
 is convenient for future references.

\begin{lemma}\label{L:gamma}
    For $\ga\in\varGamma$, $\cgd$ defined as above and $r<\ga'(+\infty)$
    \par\indent (i)~~~$\ga'_-(\cgd(r))\leq r \leq\ga'_+(\cgd(r))\,$,
    \par\indent (ii)~~$\ga(\cgd(r))=r\cgd(r)-\coga(r)$.
\end{lemma}

% DELTA

 \paragraph{2.B.} Bregman distances will be defined by means of the following functions
 of two variables. Given $\ga\in\varGamma$, for $s,t\geq0$ let
 \begin{equation}\label{E:delta}
    \trnga{s,t}\triangleq\ga(s)-\ga(t)-\ga'_{\sg{s-t}}(t)[s-t]
                                        \quad   \text{if $\ga'_+(t)$ is finite,}
 \end{equation}
 and $\trnga{s,0}\triangleq s\cdot(+\infty)$ otherwise. This definition of $\trn$
 is extended to all $(s,t)\in\R^2$, letting $\trnga{s,t}\triangleq+\infty$
 if $s<0$ or $t<0$. Beyond these cases, $\trnga{s,t}$ equals $+\infty$ if and only if
 $0=s<t$ and $\ga(0)=+\infty$, or $s>t=0$ and $\ga'_+(0)=+\infty$. The strict convexity
 of $\ga$ implies that $\trnga{s,t}\geq0$ with the equality if and only if $s=t\ge 0$.
 \index{$\trn$, $\trnnin$~~Bregman integrand}\label{Bregman integrand}

\begin{lemma}\label{L:lsc}
  For $\ga\in\varGamma$ the function $\trn\colon\R^2\to[0,+\infty]$ is lower semicontinuous.
\end{lemma}

\begin{Proof}
 By definition, $\trn$ is lsc at $(s,t)$ if $s<0$ or $t<0$. By nonnegativity,
 $\trn$ is lsc at $(s,t)$ if $s=t\geq0$. Otherwise, for $s,t\geq0$ different
 let $s_n\to s$ and $t_n\to t$ such that the sequence $\trnga{s_n,t_n}$ has
 a finite limit $r$. Thus, $s_n$ and $t_n$ are eventually nonnegative. If
 $t>0$  then the sequence $\ga'_{\sg{s_n-t_n}}(t_n)(s_n-t_n)$ has at most
 two accumulation points $\ga_{\pm}'(t)(s-t)$. Hence, $\trn$ is lsc at
 $(s,t)$ because $\ga$ is lsc at $s$ and continuous at $t$. If $t=0$ then
 $\ga'_+(0)<+\infty$ since $r$ is finite. Therefore, $\trnga{s_n,t_n}$ converges
 to $\trnga{s,0}=\ga(s)-\ga(0)-\ga'_+(0)\cdot s$ whence $\trn$ is lsc at $(s,t)$.
\end{Proof}

\begin{lemma}\label{L:limiting}
  If $s,t\geq 0$ then there exist sequences $s_n$ and $t_n$ of positive rational
  numbers such that $s_n\to s$, $t_n\to t$ and $\trnga{s_n,t_n}\to\trnga{s,t}$.
\end{lemma}

\begin{Proof}
 The assertion is trivial for $s=t$, taking $s_n=t_n$ rational.
 Otherwise, if $0\leq s<t$ then limiting along sequences $s_n\downarrow s$,
 $t_n\uparrow t$ works by the continuity of $\ga'_-$ from the left. Analogously,
 if $0\leq t<s$ then $t_n\downarrow t$, $s_n\uparrow s$ works.
\end{Proof}

\begin{lemma}\label{L:mgamma}
   If $K>0$ and $\vare>0$ then
   \[
       m_\ga^{K,\vare}\triangleq\min\big\{\min_{s\leq K}\:\trnga{s,s+\vare},
                        \min_{t\leq K}\:\trnga{t+\vare,t}\big\}
   \]
   is a positive lower bound on $\trnga{s,t}$ whenever $0\leq\min\{s,t\}\leq K$
   and $|s-t|\geq\vare$.
\end{lemma}

\begin{Proof}
 The two minima are finite and attained since $\trn$ is nonnegative and lsc,
 by Lemma~\ref{L:lsc}. Then, they cannot vanish whence $m_\ga^{K,\vare}>0$.
 For any $s\geq0$ the function $t\mapsto\trnga{s,t}$ is non-decreasing
 in $[s,+\infty)$, and thus lower bounded by $\trnga{s,s+\vare}$ for
 $t\geq s+\vare$. For any $t\geq0$ the function $s\mapsto\trnga{s,t}$
 is non-decreasing in $[t,+\infty)$, and thus lower bounded by
 $\trnga{t+\vare,t}$  when $s\geq t+\vare$. It follows for $s,t\geq0$
 that if $|s-t|\geq\vare$ then $\trnga{s,t}$ is lower bounded by
 $\trnga{s,s+\vare}$ or $\trnga{t+\vare,t}$ which implies the
 assertion.\rule{3cm}{0cm}
\end{Proof}

% GAT

 \paragraph{2.C.} For fixed $s$ with $\ga(s)$ finite, the function $t\mapsto\trnga{s,t}$
 need not be continuous on $(0,+\infty)$, and it need not be convex even
 if $\ga$ is differentiable. For fixed $t\geq0$, with $t=0$ allowed only
 when $\ga'_+(0)$ is finite, the function
 \[
    \gat\colon s\mapsto\trnga{s,t}
 \]
 differs from $\ga$  on $(0,+\infty)$ by a piecewise linear function.
 This notation is introduced for the purposes of Section~\ref{S:BREGproj}.

\begin{lemma}\label{L:gat}
  Let $\ga\in\varGamma$ and $t>0$, or $t=0$ if $\ga'_+(0)$ is finite. Then,
  $\gat\in\varGamma$, $\gat'(t)=0$ for $t>0$, $\gat'_+(t)=0$ for $t=0$,
  \begin{align*}
        \gat'_{\pm}(s)&=\ga'_{\pm}(s)-\ga'_{\sg{s-t}}(t)\,,\qquad &&s\ge0\,,\,t\neq s\,,\\
        \gat'(+\infty)&=\ga'(+\infty)-\ga'_+(t),&&\\
        \cogat(r)&=\coga(r+\ga'_{\sg{r}}(t))-\coga(\ga'_{\sg{r}}(t))\,,\qquad && r\in\R\,,\\
      (\cogat)'(r)&=(\coga)'(r+\ga'_{\sg{r}}(t))\,,\qquad &&r<\gat'(+\infty)\,.
  \end{align*}
\end{lemma}

 If $t=0$ and $r<0$ then $\ga'_{\sg{r}}(t)=\ga'_{-}(0)=-\infty$,
 and thus the values $\coga(\ga'_{\sg{r}}(t))$ and $\coga(-\infty)$
 are equal to $-\ga(0)$, see Subsection~2.A. Hence, the right-hand sides
 of the last two equations in Lemma~\ref{L:gat} equal~$0$ when $t=0$ and $r<0$.

\begin{Proof}
 The assertions on derivatives of $\gat$ follow by differentiation
 in~\eqref{E:delta} and limiting. Further,  $\cogat(0)=0$ because $\gat$
 has the global minimum $0$ attained at $t$. When computing the conjugate
 \[
    \cogat(r)={\sup}_{s>0}\;\big[rs-\ga(s)+\ga(t)+\ga'_{\sg{s-t}}(t)[s-t]\big]
 \]
 of $\gat$ at $r>0$, the supremum can be restricted to $s>t$, thus
 \[
    \cogat(r)=[\ga(t)-\ga_+'(t)t]+{\sup}_{s>t}\;\big[s[r+\ga_+'(t)]-\ga(s)\big]\,.
 \]
 By Lemma~\ref{L:gamma0}, the first term is equal to $-\coga(\ga_+'(t))$,
 and the second one to $\coga(r+\ga_+'(t))$ since $r>0$. The conjugate
 is computed similarly at $r<0$ when $t>0$, in which case the supremum
 can be restricted to $s<t$. If $r<0$ and $t=0$ then $\cogat(r)$ is
 equal to $\ga(0)+\coga(r+\ga_+'(0))=0$. The last assertion follows
 from the last but one by differentiation.
\end{Proof}

\begin{lemma}\label{L:compB}
  Let $\ga\in\varGamma$, $s\ge 0$ and $t>0$, or $t=0$ if $\ga'_+(0)$ is finite.
  Then, $\trnarg{\gat}{s,t}$ equals $\trnga{s,t}$ and for $r\ge0$ different from $t$
  \[
   \trnga{s,r}=\trnarg{\gat}{s,r}+
        [\ga'_{\sg{s-t}}(t)-\ga'_{\sg{r-t}}(t)][s-t]\,.
  \]
\end{lemma}

\begin{Proof}
 Excluding the case $r=0$, $\ga'_+(0)=-\infty$, Lemma~\ref{L:gat} implies that
 $\gat'_{+}(r)$ is finite and
 \[
   \trnarg{\gat}{s,r}=\trnga{s,t}-\trnga{r,t}-[\ga'_{\sg{s-r}}(r)-\ga'_{\sg{r-t}}(t)][s-r]
 \]
 when $r\neq t$. The assertion follows by a simple calculation.
 In the excluded case both sides are $0$ or $+\infty$ according as $s=0$ or $s>0$.
 \end{Proof}

% LADGA

 \paragraph{2.D.} The correction term in \eqref{E:Pyth} will be constructed by means of the function
 \begin{equation}\label{E:upsi}
    \ladga{s,r}\triangleq[\ga'_{\sg{s-\cgd(r)}}(\cgd(r))-r][s-\cgd(r)]\,,\quad
        \quad s\ge 0\,,\,r<\ga'(+\infty)\,,
 \end{equation}
 where $\ga\in\varGamma$ and $\cgd(r)=(\coga)'(r)$ as in Subsection~2.A.
 By Lemma~\ref{L:gamma}\emph{(i)}, $\ladga{s,r}\geq0$. This quantity is identically zero
 if $\ga$ is essentially smooth, in which case $\cgd(r)>0$ and $\ga'(\cgd(r))=r$ for
 $r<\ga'(+\infty)$. If $\ga$ is differentiable on $(0, +\infty)$ then
 $\ladga{s,r}$ equals $|\ga'_+(0)-r|_+\cdot s$.%
 \index{$\lad$, $\ladnin$~~integrand behind the correction}\label{integrand behind the correction}

\begin{lemma}\label{L:keyid}
  For $\ga\in\varGamma$, $s\ge0$ and $r<\ga'(+\infty)$
  \[
      \ga(s)+\coga(r) = {r}{s}+ \trnga{s,\cgd(r)}+ \ladga{s,r}\,.
  \]
\end{lemma}

\begin{Proof}
 The assumptions and Lemma~\ref{L:gamma}\emph{(i)} imply that
 $\ga'_+(\cgd(r))$ is finite. Then, by the definitions of $\trn$ and
 $\lad$, the right-hand side equals
 \[
    {r}{s}+ \ga(s)-\ga(\cgd(r))-r[s-\cgd(r)]= \ga(s)-\ga(\cgd(r))+r\cgd(r)\,.
 \]
 Hence, the assertion follows by Lemma~\ref{L:gamma}\emph{(ii)}.
\end{Proof}

\bigskip

 In Section~\ref{S:dual}, the following analogue of $\trn$
 \[
    \trndu{r_2,r_1}\triangleq\coga(r_2)-\coga(r_1)-\cgd(r_1)[r_2-r_1]\,,
        \qquad r_1,r_2<\ga'(+\infty)\,,
 \]
 with $\coga$ replacing $\ga$ is needed.

\begin{lemma}\label{L:Bld}
   For $\ga\in\varGamma$ and $r_1,r_2<\ga'(+\infty)$
   \[
        \trndu{r_2,r_1}=\trnga{\cgd(r_1),\cgd(r_2)}+\ladga{\cgd(r_1),r_2}\,.
   \]
\end{lemma}

\begin{Proof}
 By definition and Lemma~\ref{L:gamma}\emph{(ii)},
 $\trndu{r_2,r_1}=\coga(r_2)+\ga(\cgd(r_1))-r_2\cgd(r_1)$.
 Hence the assertion follows by  Lemma~\ref{L:keyid}.
\end{Proof}

% INTERGRANDS

\paragraph{2.E.} Following \cite[Chapter~14D]{RW},
 a function $f\colon\TT\times\R^k\to[-\infty,+\infty]$ is an \emph{integrand}
 if for every $x\in\R^k$ the function $\ttt\mapsto f(\ttt,x)$ on~$\TT$,
 denoted also by $f(\cdot,x)$, is $\calT$-measurable. An integrand
 as above is convex/lsc if $f(\ttt,\cdot)$ is convex/lsc for every
 $\ttt\in\TT$, and this convention extends naturally. A function
 $f\colon\TT\times\R^k\to[-\infty,+\infty]$ is a \emph{normal
 integrand} if $f(\ttt,\cdot)$ is lsc for each $\ttt\in\TT$ and
 $\ttt\mapsto\inf_{x\in D}f(\ttt,x)$ is $\calT$-measurable for each open ball
 $D\subset\R^k$. This is not the original definition of \cite[Chapter~14D]{RW}
 but is equivalent to it by \cite[Proposition~14.40 and Proposition~28]{RW}.
 In particular, a normal integrand is an integrand.\index{$B$~~the class of integrands $\nm$}\label{the class of integrands}

 The class of integrands $\nm\colon\TT\times\R\to(-\infty,+\infty]$ such
 that $\nm(\ttt,\cdot)\in\varGamma$ for all $\ttt\in\TT$ is denoted by $B$.
 The assumption $\nm\in B$ is the only restriction on $\nm$ adopted throughout
 this paper. By \cite[Proposition 14.39]{RW}, each $\nm\in B$ is a normal
 convex integrand. The conjugate $\conm$ of $\nm\in B$ inherits this property
 \cite[Theorem 14.50]{RW}. Hence, the function $\ttt\mapsto\nm(\ttt,g(\ttt))$
 in the definition of $\fuc_\nm$, and $\ttt\mapsto\conm(\ttt,\zat{\vte}{\moma(\ttt)})$
 in the definition of $\dufu_{\nm}$, are $\calT$-measurable \cite[Proposition 14.28]{RW}.

 By \cite[Proposition 14.56]{RW}, for $t>0$ the functions $\nm'_+(\cdot,t)$
 and $\nm'_-(\cdot,t)$ are $\calT$-measurable and so are their monotone limits
 $\nm'_+(\cdot,0)$ and $\nm'(\cdot,+\infty)$. Hence, the set\index{cofinite function}\label{cofinite function 2}
 $\{\nm'(\cdot,+\infty)=+\infty\}$ of $\ttt\in\TT$ for which $\nm(\ttt,\cdot)$
 is cofinite, denoted by $\cofi$, belongs to $\calT$. Let $\asli$ denote
 the set of those $\ttt\in\TT\sm\cofi$ for which $\nm(\ttt,\cdot)$
 is asymptotically linear. Then, $\asli\in\calT$ since\index{asymptotically linear function}\label{asymptotically linear function 2}
 $\ttt\mapsto\conm(\ttt,\nm'(\ttt,+\infty))$ is
 $\calT$-measurable on $\TT\sm\cofi$. \index{$\cofi$, $\asli$~~special subsets of $\TT$}\label{special subsets of TT}

 Recall the notation $\trn$ from Subsection~2.D. Assuming the
 integrand $\nm$ is in $B$, the function $(\ttt,s,t)\mapsto\trnarg{\ga}{s,t}$,
 where $\ga=\nm(\ttt,\cdot)$ depends on $\ttt$, is denoted by $\trnnin$.

\begin{lemma}\label{L:nindelta}
   If $\nm\in B$ then $\trnnin\colon\TT\times\R^2\to[0,+\infty]$
   is a normal integrand.
\end{lemma}

\begin{Proof}
 By Lemma~\ref{L:lsc}, for each $\ttt\in\TT$ the function $(s,t)\mapsto\trnnin(\ttt,s,t)$ is lsc.
 Hence, it suffices to prove that the function $\ttt\mapsto\inf_{(s,t)\in D}\,\trnnin(\ttt,s,t)$
 is $\calT$-measurable for each open ball $D\subset\R^2$. The infimum is identically
 $+\infty$ if $D$ contains no points with positive coordinates. Otherwise,
 Lemma~\ref{L:limiting} implies that to each $\ttt\in\TT$ and $(s,t)\in D$
 with $s,t\geq0$ there exist sequences of positive rational numbers $s_n\to s$
 and $t_n\to t$ such that $\trnnin(\ttt,s_n,t_n)\to\trnnin(\ttt,s,t)$.
 It follows that the above infimum does not change when admitting only
 $(s,t)\in D$ with positive rational coordinates. The measurability of
 this infimum over countably many pairs is implied by the measurability
 of the individual functions $\ttt\mapsto\trnnin(\ttt,s,t)$, $s,t>0$,
 which in turn follows from the measurability of $\nm(\cdot,s)$,
 $\nm(\cdot,t)$, $\nm'_+(\cdot,t)$ resp.\ $\nm'_-(\cdot,t)$.
\end{Proof}

% BREGMAN

\paragraph{2.F.} For $\nm\in B$ the \emph{Bregman distance} of \calT-measurable
 functions $g,h$ is defined as\index{$\Bres$~~Bregman distance}\label{Bregman distance}
 \begin{equation}\label{E:Brebe}
    \Bre{g,h}\triangleq\inttt{\TT}\:\trnnin(\ttt,g(\ttt),h(\ttt))\:\mu({\mathrm d}\ttt)
 \end{equation}
 where the integrated function is \calT-measurable by Lemma~\ref{L:nindelta}
 and \cite[Proposition 14.28]{RW}. The quantity $\Bre{g,h}$ can be finite only
 for  $g$ and $h$ nonnegative $\mu$-a.e., by the definition of $\trnnin$.
 The Bregman distance is not a metric on nonnegative functions, nevertheless
 $\Bre{g,h}\ge 0$ with the equality if and only if  $g=h\ge 0\,[\mu]$.

 The Bregman distances corresponding to the autonomous integrands
 $t\ln t,\; -\ln t,\; t^2$ mentioned in Subsection~1.B.
 are the Kullback--Leibler distance ($I$-divergence), Itakura-Saito
 distance, and squared $L^2$-distance.

\begin{remark}\label{R:Bre}\rm
  In the literature the term Bregman distance frequently refers to
  nonmetric  distances on $\R^k$ associated with convex functions
  $\phi\colon\R^k\to(-\infty,+\infty]$. Typically, $\phi$ is assumed
  differentiable in the interior of its effective domain and the Bregman
  distance of $x$ in $\dom{\phi}$ from $y$ in the interior is
  defined as $\phi(x)-\phi(y)-\zat{\nabla\phi(y)}{x-y}$. For $k=1$ and
  $\phi=\ga\in\varGamma$ differentiable in $(0,+\infty)$ this reduces
  to $\trnga{x,y}$ from Subsection~2.D. The special case
  of \eqref{E:Brebe} for $\TT= \{1,\ldots,k\}$ gives the Bregman
  distance in this sense of the vector $x=(g(1),\ldots,g(k))$ from
  $y=(h(1),\ldots,h(k))$, associated with the convex function
  \[
    \phi(x_1,\ldots,x_k)=\mbox{$\sum$}_{\ttt\in\TT}\;\mu(\ttt)\,\nm(\ttt,x_\ttt)\,,
        \qquad (x_1,\ldots,x_k)\in\R^k\,.
  \]
  As this $\phi$ is a sum of convex functions of individual coordinates,
  the associated Bregman distance is \emph{separable}.
\end{remark}

\begin{lemma}\label{L:Brelsc}\index{$\rightsquigarrow$~~local convergence in measure}\label{local convergence in measure 2}
   For $\nm\in B$ the Bregman distance $\Bres$ is lsc for local convergence
   in measure, jointly in both coordinates.
\end{lemma}

\begin{Proof}
 Assuming $g_n\rightsquigarrow g$, $h_n\rightsquigarrow h$
 and $\liminf_{n\to +\infty}\Bre{g_n,h_n}=r$, there exists
 an increasing sequence $n_k$ such that the subsequence $g_{n_k}$
 converges to $g$ and $h_{n_k}$ to $h$, both $\mu$-a.e.
 Then,
 \[
   r\geq\inttt{\TT}\:\liminf\nolimits_{k\to\infty}\:
                \trnnin(\ttt,g_{n_k}(\ttt),h_{n_k}(\ttt))\:\mu({\mathrm d}\ttt)
                \ge\Bre{g,h}\,,
 \]
 by Fatou lemma and the lower semicontinuity of $\trn$,
 see Lemma~\ref{L:lsc}.
\end{Proof}

\begin{lemma}\label{L:inmeas}
  Given $\nm\in B$, to any set $C\in\calT$ of finite $\mu$-measure
  and positive numbers $K$, $\xi$ and $\vare$ there exists
  $\delta>0$ such that for \calT-measurable functions $g$ and $h$
  either of $\Bre{g,h}\leq\delta$ or $\Bre{h,g}\leq\delta$ implies
         \[
            \mu(C\cap\{|g-h|>\vare\})<\xi + \mu(C\cap\{g>K\})\,.
         \]
\end{lemma}

\begin{Proof}
 The function $m_\nm^{K,\vare}$ given by
 \[
       m_\nm^{K,\vare}(\ttt)\triangleq
            \min\big\{\min_{s\leq K}\:\trnnin(\ttt,s,s+\vare),
              \min_{t\leq K}\:\trnnin(\ttt,t+\vare,t)\big\}\,,\qquad \ttt\in\TT\,,
 \]
 is positive by Lemma~\ref{L:mgamma} and
 \calT-measurable by Lemma~\ref{L:nindelta} and \cite[Proposition 14.37]{RW}.
 This and $\mu(C)<+\infty$ imply that $\mu(C\cap \{m_\nm^{K,\vare}<\eta\})<\frac12\xi$
 for $\eta>0$ sufficiently close~to~$0$. Let $\delta$ be equal to $\frac12 \eta\xi$.
 By the definition of Bregman distance and Lemma~\ref{L:mgamma},
 whenever $\delta\geq\min\{\Bre{g,h},\Bre{h,g}\}$
 \[\begin{split}
    \tfrac12 \eta\xi&\geq
        \inttt{\{g\leq K,|g-h|>\vare\}}\:
        \min\big\{\trnarg{\nm}{\ttt,g(\ttt),h(\ttt)},\trnarg{\nm}{\ttt,h(\ttt),g(\ttt)}\big\} \;\mu({\mathrm d}\ttt)\\
            &\geq \inttt{\{g\leq K,|g-h|>\vare\}}\:  m_\nm^{K,\vare} \;{\mathrm d}\mu
            \geq \eta\cdot\mu\big(C\cap\{m_\nm^{K,\vare}\geq\eta\}\cap\{g\leq K,|g-h|>\vare\}\big)\,.
    \end{split}
 \]
 Therefore,
 \[
    \mu(C\cap\{|g-h|>\vare\})\leq
       \mu(C\cap\{m_\nm^{K,\vare}<\eta\})+\mu(C\cap\{g> K\})+\tfrac12\xi
 \]
 and the assertion follows by the choice of $\eta$.
\end{Proof}

\begin{corollary}\label{C:inmeas}
  If a sequence of $\calT$-measurable functions $g_n$ converges
  to a $\calT$-measur\-able function~$g$ either in the sense
  $\Bre{g_n,g}\to0$ or $\Bre{g,g_n}\to0$, then $g_n\rightsquigarrow g$.
\end{corollary}

\noindent
 Note that for certain integrands $\nm\in B$ the hypothesis of this corollary
 admits to conclude even $L_1$-convergence, see \cite[Lemma~3]{Csi.genproj}.

%333333333333333333333333333333333333333333333333333333333333333333333333333333333333333333
\section{Preliminaries on the dual problem}\label{S:dualproblem}

 This section collects auxiliary results on the function $\dufu_{\nm}$ and the
 subset $\Theta_{\nm}$ of its effective domain. Proposition~\ref{P:properH}
 provides a new sufficient condition for $\vafu_{\nm}\not\equiv +\infty$.

\begin{lemma}\label{L:domdufu}
    If $\vte\in\dom{\dufu_{\nm}}$ then
    $\zat{\vte}{\moma}\leq\nm'(\cdot,+\infty)\,[\mu]$
    with the strict inequality on $\TT\sm\asli$.
\end{lemma}

\begin{Proof}
 The effective domain of $\dufu_{\nm}$ consists of those
 $\vte\in \R^d$ for which the positive part of the integral
 $\int_{\TT}\:\conm(\ttt,\zat{\vte}{\moma(\ttt)})\:\mu({\mathrm d}\ttt)$
 is finite. For such $\vte$ the integrand is finite $\mu$-a.e.,
 and since $\nm\in B$, the assertion follows.
\end{Proof}

\bigskip

 Recall that the set $\Theta_{\nm}$ consists of those $\vte\in\dom{\dufu_{\nm}}$
 for which $\conm(\ttt,\cdot)$ is finite around $\zat{\vte}{\moma(\ttt)}$
 for $\mu$-a.a.\ $\ttt\in\TT$. Since $\nm\in B$, by properties of $\coga$
 discussed in Subsection~2.A,
 \begin{equation}\label{E:DCQ}
      \text{$\vte\in\dom{\dufu_{\nm}}$ belongs to $\Theta_{\nm}$ if and only if
        $\zat{\vte}{\moma}<\nm'(\cdot,+\infty)\,[\mu]$.}
 \end{equation}
 The inequality in \eqref{E:DCQ} is equivalent to existence of the derivative
 $(\conm)'(\ttt,\cdot)$ at $\zat{\vte}{\moma(\ttt)}$.

\begin{corollary}\label{C:Thetaindom}
    $\Theta_{\nm}=\{\vte\in\dom{\dufu_{\nm}}\colon
    \zat{\vte}{\moma}<\nm'(\cdot,+\infty)\:\text{$\mu$-a.e.\ on $\asli$}\}$.
\end{corollary}

\begin{remark}\label{R:momas}\rm
 The set $\Theta_{\nm}$ may be empty even if $\dufu_{\nm}$ is proper,
 see Examples~\ref{Ex:vafu.dual} and \ref{Ex:emqua}. However,
 $\mu(\asli)=0$ implies $\Theta_{\nm}=\dom{\dufu_{\nm}}$, by
 Corollary~\ref{C:Thetaindom}. In case $\dom{\dufu_{\nm}}\neq\pmn$,
 another sufficient condition for $\Theta_{\nm}\neq\pmn$ is
 the \emph{moment assumption}\index{moment assumption}\label{moment assumption}
 \begin{equation}\label{E:ma}
      \zat{\teob}{\moma}>0\,\;\;[\mu]\,,
                  \quad\text{for some $\teob\in\R^d$.}
  \end{equation}
 In fact, \eqref{E:ma} implies $\vte-t\teob\in\Theta_{\nm}$ for
 $\vte\in\dom{\dufu_{\nm}}$ and $t>0$, by the monotonicity of each
 $\conm(\ttt,\cdot)$ and \eqref{E:DCQ}. The moment assumption holds
 for example when one coordinate of $\moma$ is a nonzero constant.
 If $\nm'(\cdot,+\infty)<0$ $[\mu]$ then \eqref{E:ma} holds for
 each $\teob\in\R^d$ with $-\teob\in\dom{\dufu_{\nm}}$.

 If $\nm'(\cdot,+\infty)>0$ $[\mu]$ then a trivial sufficient condition
 for $\Theta_{\nm}\neq\pmn$ is the $\mu$-integrability of $\conm(\ttt,0)$,
 for it implies that $\vte=\bnu$ belongs to $\Theta_{\nm}$.
\end{remark}

\begin{lemma}\label{L:Theta}
   Under the \DCQ  $\Theta_{\nm}\neq\pmn$, the sets
   $\dom{\dufu_{\nm}}$ and $\Theta_{\nm}$ have the same relative interior,
   and dual values do not change when the maximization in \eqref{E:duva}
   is restricted to~$\Theta_{\nm}$.
\end{lemma}

\begin{Proof}
 Let $\vte\in\dom{\dufu_{\nm}}$, $\teob\in\Theta_{\nm}$ and $0<t<1$.
 Then, $\te_t=t\teob+(1-t)\vte$ belongs to $\dom{\dufu_{\nm}}$.
 By Lemma~\ref{L:domdufu}, $\zat{\vte}{\moma}\leq\nm'(\cdot,+\infty)\,[\mu]$,
 and $\zat{\teob}{\moma}<\nm'(\cdot,+\infty)\,[\mu]$ follows from \eqref{E:DCQ}.
 Then, the latter inequality holds for $\te_t$ instead of $\teob$.
 By~\eqref{E:DCQ}, $\te_t$ belongs to~$\Theta_{\nm}$. Since $\te_t\to\vte$
 as $t\downarrow0$, this proves that $\vte\in\cl{\Theta_{\nm}}$. Therefore,
 $\dom{\dufu_{\nm}}$ is contained in $\cl{\Theta_{\nm}}$. As
 $\Theta_{\nm}\pdm\dom{\dufu_{\nm}}$ is convex, this proves
 the first assertion. The second assertion follows from the
 convexity of $\dufu_{\nm}$.
\end{Proof}

\bigskip

 Each function $f_{\vte}\in\FF_{\nm}$, $\vte\in\Theta_{\nm}$, defined in \eqref{E:fte}
 is $\calT$-measurable by \cite[Proposition 14.56]{RW}. It is nonnegative, vanishes
 on the set $\{\zat{\vte}{\moma}\le\nm'_+(\cdot,0)\}$, and is $\mu$-a.e.\ positive
 on the complement of this set. In particular, if $\nm$ is essentially smooth
 then $f_{\vte}>0\;[\mu]$.

\begin{remark}\label{L:unique.dual}\rm
  The parametrization of $\FF_\nm$ in~\eqref{E:fte} is not bijective, in
  general. However, if $\nm$ is essentially  smooth then $(\conm)'(\ttt,\cdot)$
  is strictly increasing and maps $(-\infty,\nm'(\ttt,+\infty))$ onto
  $(0,+\infty)$, hence the function $f_\vte\in\FF_\nm$ does determine the
  function $\zat{\vte}{\moma}$. In this case, the parametrization is bijective
  when the component functions of \moma are linearly independent under $\mu$.
  This has not been assumed because in modified primal and dual problems
  of Section~\ref{S:main} restrictions of $\mu$ occur under which
  the independence is lost anyhow.
\end{remark}

\begin{lemma}\label{L:limint}
    For $\teob,\vte\in\R^d$ with $\dufu_{\nm}(\teob)$
    and $\dufu_{\nm}(\vte)$ finite, the directional derivative
    \[
        \dufu'_{\nm}(\vte;\,\teob-\vte)\triangleq
        \lim_{t\downarrow 0}\tfrac{1}{t}[\dufu_{\nm}(\vte+t(\teob-\vte))-\dufu_{\nm}(\vte)]
    \]
    is equal to $-\infty$ if the set $\{\zat{\vte}{\moma}=\nm'(\cdot,+\infty)\,,\,
    \zat{\teob}{\moma}\neq\zat{\vte}{\moma}\}$
    has positive $\mu$-measure. Otherwise, it equals
    \begin{equation}\label{E:integr}
        \inttt{\{\zat{\vte}{\moma}<\nm'(\cdot,+\infty)\}}\:
                \zat{\teob-\vte}{\moma(\ttt)}
                    \cdot(\conm)'(\ttt,\zat{\vte}{\moma(\ttt)})\:\mu({\mathrm d}\ttt)
    \end{equation}
    where the positive part of the integral is finite.
\end{lemma}

\begin{Proof}
 Since $\dufu_{\nm}(\teob)$ and $\dufu_{\nm}(\vte)$ are finite,
 $\conm(\ttt,\zat{\vte}{\moma(\ttt)})$ and
 $\conm(\ttt,\zat{\teob}{\moma(\ttt)})$ are finite
 for $\mu$-a.a.\ $\ttt\in\TT$. For those $\ttt$, the function
 $\phi_\ttt$ given by
 \[
   \phi_\ttt(t)=\tfrac1t\big[\conm\big(\ttt,\zat{\vte}
        {\moma(\ttt)}+t\zat{\teob-\vte}{\moma(\ttt)}\big)
             -\conm\big(\ttt,\zat{\vte}{\moma(\ttt)}\big)\big]\,,
             \qquad t>0\,,
 \]
 is non-increasing as $t\downarrow0$, by convexity.
 If $\zat{\vte}{\moma(\ttt)}<\nm'(\ttt,+\infty)$ then $\conm(z,\cdot)$
 is differentiable at $\zat{\vte}{\moma(\ttt)}$ and $\phi_\ttt(t)$ tends to $\zat{\teob-\vte}{\moma(\ttt)}\cdot(\conm)'(\ttt,\zat{\vte}{\moma(\ttt)})$.
 If $\zat{\vte}{\moma(\ttt)}=\nm'(\ttt,+\infty)$, thus $\ttt\in\asli$,
 then $\phi_\ttt(t)$  tends to $0$ or $-\infty$ according to
 $\zat{\teob}{\moma(\ttt)}=\zat{\vte}{\moma(\ttt)}$ or not.
 It follows that
 \[
   \dufu'_{\nm}(\vte;\,\teob-\vte)=
   \lim_{t\downarrow0}\:\inttt{\{\zat{\teob}{\moma}\neq\zat{\vte}{\moma}\}}\:\phi_\ttt(t)\:\mu({\mathrm d}\ttt)
 \]
 where the limiting and integration can be interchanged by monotone
 convergence. Then, the limit is $-\infty$ if the set
 $\{\zat{\vte}{\moma}=\nm'(\cdot,+\infty)\,,\,
    \zat{\teob}{\moma}\neq\zat{\vte}{\moma}\}$
 is not $\mu$-negligible. Otherwise, $\zat{\vte}{\moma}<\nm'(\cdot,+\infty)\,[\mu]$
 on $\{\zat{\teob}{\moma}\neq\zat{\vte}{\moma}\}$ on account of Lemma~\ref{L:domdufu},
 and the second assertion follows by the interchange. The integral \eqref{E:integr}
 cannot be $+\infty$ by monotonicity.
\end{Proof}

\begin{remark}\label{R:limint}\rm
    If $\vte\in\Theta_{\nm}$ then the integral \eqref{E:integr} is equal to
    $\int_{\TT}\:\zat{\teob-\vte}{\moma}f_\vte\:{\mathrm d}\mu$, see~\eqref{E:DCQ}.
\end{remark}

\begin{corollary}\label{C:moma0}
    If $\dufu_{\nm}$ is finite in a neighborhood of $\vte$
    then $\nm^*(\ttt,\cdot)$ is differentiable at $\zat{\vte}{\moma(\ttt)}$
    for $\mu$-a.a.\ $\ttt\in\TT$ with $\moma(\ttt)\neq\bnu$,
    $\dufu_{\nm}$ is differentiable at $\vte$, and
    \[
        \nabla\dufu_{\nm}(\vte)=
        \inttt{\{\moma\neq\bnu\}}\:\moma(\ttt)
                    \cdot(\conm)'(\ttt,\zat{\vte}{\moma(\ttt)})\:\mu({\mathrm d}\ttt)\,.
    \]
    If additionally the set $\{\moma=\bnu\}$ is $\mu$-negligible then $\vte\in\Theta_\nm$.
\end{corollary}

\begin{Proof}
 Since $\dufu_{\nm}$ is convex, the hypothesis implies that all
 directional derivatives at $\vte$ are finite. Therefore,
 Lemma~\ref{L:limint} implies for each $\teob$ sufficiently close to
 $\vte$ that $\zat{\teob}{\moma}\neq\zat{\vte}{\moma}$ $[\mu]$ on the
 set $\{\zat{\vte}{\moma}=\nm'(\cdot,+\infty)\}$. It follows that on
 this set $\moma=\bnu\,[\mu]$. Hence, recalling Lemma~\ref{L:domdufu},
 the integral \eqref{E:integr} can be equivalently taken over
 $\{\moma\neq\bnu\}$ and the assertion follows.
 If $\mu(\{\moma=\bnu\})=0$ then $\vte\in\Theta_\nm$ by~\eqref{E:DCQ}.
\end{Proof}

\begin{corollary}\label{C:neededinBoston}
    Under finiteness of $\dufu_{\nm}$ on an open set, the \DCQ holds if and only if
    the set where $\moma=\bnu$ and $\lim_{t\uparrow+\infty} \nm(\cdot,t)$ is finite
    is $\mu$-negligible.
\end{corollary}

\begin{Proof}
 This follows from \eqref{E:DCQ} and Corollary~\ref{C:moma0}.
\end{Proof}

\begin{proposition}\label{P:properH}
    If $\dufu_{\nm}$ is finite in a neighborhood of $\vte$
    then $\vafu_{\nm}(\nabla\dufu_{\nm}(\vte))<+\infty$.
\end{proposition}

\begin{Proof}
 Since $\dufu_{\nm}(\vte)$ is finite, both integrals in the sum
 \[
    \dufu_{\nm}(\vte)=\inttt{\{\moma\neq\bnu\}}\:
        \conm\big(\ttt,\zat{\vte}{\moma(\ttt)}\big)\:\mu({\mathrm d}\ttt)
                    +\inttt{\{\moma=\bnu\}}\:\conm(\cdot,0)\:{\mathrm d}\mu
 \]
 are finite. Then, the function $\inf_t \nm(\cdot,t)=-\conm(\cdot,0)$
 is $\mu$-integrable on the set \mbox{$\{\moma=\bnu\}$}. By Lemma~\ref{L:pomstand},
 $\int_{\{\moma=\bnu\}}\:\nm(\ttt,h(\ttt))\:\mu({\mathrm d}\ttt)$ is
 finite for some $\calT$-measurable function $h$. Let $g$
 be the function defined by $g(\ttt)=h(\ttt)$ if
 $\moma(\ttt)=\bnu$ and $g(\ttt)=(\conm)'(\ttt,\zat{\vte}{\moma(\ttt)})$
 otherwise, $\ttt\in\TT$. Corollary~\ref{C:moma0} implies that this derivative
 exists $\mu$-a.e.\ and $\int_\TT\,\moma g\,{\mathrm d}\mu$ equals $\nabla\dufu_{\nm}(\vte)$.
 By Lemma~\ref{L:gamma}\emph{(ii)},
 \[
    \inttt{\{\moma\neq\bnu\}}\:
        \nm\big(\ttt,g(\ttt))\:\mu({\mathrm d}\ttt)=
    \zat{\vte}{\nabla\dufu_{\nm}(\vte)}-
    \inttt{\{\moma\neq\bnu\}}\:
        \conm\big(\ttt,\zat{\vte}{\moma(\ttt)}\big)\:\mu({\mathrm d}\ttt)
 \]
 where the right-hand side is finite. It follows that $\fuc_{\nm}(g)$ is finite
 and the assertion follows.
\end{Proof}

\begin{corollary}\label{C:dufuint}
    If $\dufu_{\nm}$ is finite on an open set then $\vafu_{\nm}$ is proper,
    and $\dufu_{\nm}$ is lsc.
\end{corollary}

\begin{Proof}
 By Proposition~\ref{P:properH},  $\vafu_{\nm}\not\equiv +\infty$ .
 Therefore, $\dufu_{\nm}=\vafu^*_{\nm}$ by Theorem~\ref{T:repre}.
 The assumption implies that $\dufu_\nm$ is proper hence so is also
 $\vafu_{\nm}$. As $\dufu_\nm$ is a convex conjugate, it is lsc.
\end{Proof}

\bigskip

 The hypothesis in Proposition~\ref{P:properH} is equivalent to
 assuming that $\dufu_{\nm}$ is proper and its effective domain
 has nonempty interior. To conclude $\vafu_{\nm}\not\equiv +\infty$,
 neither of the two  assumptions can be omitted, see
 Examples \ref{Ex:vafu.dual} and~\ref{Ex:noPCQ.DCQ}.

 The following lemma addresses, for later reference, essential smoothness of $\dufu_{\nm}$.

\begin{lemma}\label{L:essmooth}
  The function $\dufu_{\nm}$ is essentially smooth if and only if it is finite
  on an open set and the subdifferential $\partial \dufu_{\nm}(\vte)$ is empty
  for those $\vte$ in $\dom{\dufu_{\nm}}$ that are not in its interior. Here, the
  condition on emptiness is equivalent to $\dufu'_{\nm}(\vte;\,\teob-\vte)=-\infty$
  for each $\teob$ in the interior of $\dom{\dufu_{\nm}}$.
\end{lemma}

\begin{Proof}
 If $\dufu_{\nm}$ is finite on an open set then it is lsc by
 Corollary~\ref{C:dufuint}, and  differentiable in the interior
 of $\dom{\dufu_{\nm}}$ by Corollary~\ref{C:moma0}. Hence, the
 assertion follows from \cite[Theorem 26.1]{Rock} and the proof of
 \cite[Lemma 26.2]{Rock}.
\end{Proof}

%4444444444444444444444444444444444444444444444444444444444444444444444444444444444444444
\section{The constraint qualifications}\label{S:roleCQ}

 Most results of this section are well-known in more restrictive frameworks,
 typically for autonomous integrands which are essentially smooth
 or at least differentiable. The examples and Figure~\ref{F:1}
 in Section~\ref{S:examples} illustrate several situations
 encountered below. Lemma~\ref{L:duPyth} and
 Theorem~\ref{T:geprisol} are new results.

% FENCHEL and its integration

\paragraph{4.A.} Two simple lemmas are sent forward.

\begin{lemma}\label{L:duleqpri}
    $\vafu_{\nm}\geq\dufu_{\nm}^*$.
\end{lemma}

\noindent
 In particular, $\vafu_{\nm}=\dufu_{\nm}^*\equiv+\infty$
 if $\dufu_{\nm}$ attains the value $-\infty$.

\begin{Proof}
 By Fenchel inequality, for any $\vte\in\R^d$ and function $g$ on $\TT$
 \begin{equation}\label{E:Fen}
    \nm(\ttt,g(\ttt))+\conm\big(\ttt,\zat{\vte}{\moma(\ttt)}\big)
    \geq\zat{\vte}{\moma(\ttt)}g(\ttt)\,,\quad \ttt\in\TT\,.
 \end{equation}
 Integrating, for $a\in\R^d$
 \begin{equation}\label{E:Fen.int}
    \inttt{\TT}\:\big[\,\nm(\ttt,g(\ttt))+
            \conm\big(\ttt,\zat{\vte}{\moma(\ttt)}\big)\big]\:\mu({\mathrm d}\ttt)
                \geq\zat{\vte}{a}\,,\qquad g\in\calG_a^+\,.
 \end{equation}
 If $\dufu_{\nm}(\vte)=-\infty$ for some $\vte\in\R^d$ then this inequality
 implies $\fuc_\nm\equiv+\infty$. Otherwise, $\dufu_{\nm}$ is finite
 on its effective domain and $\fuc_\nm(g)\geq\zat{\vte}{a}- \dufu_{\nm}(\vte)$
 holds for every $g\in\calG_a$ and $\vte\in\R$.
\end{Proof}

\begin{lemma}\label{L:pri=du}
      If the \PCQ holds for $a\in\R^d$ then $\vafu_{\nm}(a)=\dufu_{\nm}^{*}(a)$
      and a dual solution for $a$ exists.
\end{lemma}

\begin{Proof}
 By Theorem~\ref{T:repre}, $\vafu_{\nm}^*=\dufu_{\nm}$ whence the equality
 is a consequence of the equality $\vafu_{\nm}(a)=\vafu_{\nm}^{**}(a)$, valid in \ri{\dom{\vafu_{\nm}}}.
 The existence of a dual solution follows from \cite[Theorems 23.4 and 23.5]{Rock}.
\end{Proof}

\begin{remark}\label{R:dualatt}\rm
 The \PCQ is also necessary for the existence of a dual solution if
 $\vafu_\nm\not\equiv +\infty$ and \nm is essentially smooth, see
 Corollary~\ref{C:dualatt}, but not in general, see Example~\ref{Ex:zero.density}.
\end{remark}

 In a `regular' situation the families $\calG_a$ and $\FF_\nm$ intersect.

\begin{lemma}\label{L:basic}
   If $a\in\R^d$ and $f_{\teob}\in\calG_a$ for some $\teob\in\Theta_{\nm}$ with
   $\dufu_{\nm}(\teob)$ finite then
   \[
    \vafu_{\nm}(a)
        =\fuc_\nm(f_{\teob})
        =\zat{\teob}{a}-\dufu_{\nm}(\teob)=\dufu_{\nm}^*(a)
        =\vafu_{\nm}^{**}(a)\,,
   \]
  the primal solution $g_a$ exists, $g_a=f_{\teob}$,
  $\teob$ is a dual solution for $a$, and $\calG_a\cap\FF_\nm=\{g_a\}$.
\end{lemma}

\begin{Proof}
 For $\vte\in\R^d$, $\ttt\in\TT$ and a function $g$, ineq.~\eqref{E:Fen} is
 tight if and only if $g(z)$ is equal to the derivative of $\conm(\ttt,\,.\,)$
 at $\zat{\vte}{\moma(\ttt)}$ \cite[Theorem 23.5]{Rock}. It follows that
 ineq.~\eqref{E:Fen.int} is tight if and only if $\vte\in\Theta_{\nm}$
 and $g=f_\vte$. This and finiteness of $\dufu_{\nm}(\teob)$ imply
 $\fuc_\nm(g)+\dufu_{\nm}(\teob)\geq\zat{\teob}{a}$, $g\in\calG_a$,
 with the equality if and only if $g=f_\teob$. By assumption,
 $f_{\teob}\in\calG_a$, and thus $\fuc_\nm(f_{\teob})$ is finite
 and equals $\zat{\teob}{a}-\dufu_{\nm}(\teob)$. Therefore,
 $\fuc_\nm(g)\geq\fuc_\nm(f_{\teob})$, $g\in\calG_a$. This proves
 that $\vafu_{\nm}(a)=\fuc_\nm(f_{\teob})$, the primal solution
 $g_a$ exists and equals $f_{\teob}$. By Theorem~\ref{T:repre},
 $\vafu_{\nm}^*=\dufu_{\nm}$, and thus $\dufu_{\nm}^*(a)=\vafu_{\nm}^{**}(a)$.
 Therefore, the two inequalities in the chain
 \[
    \vafu_{\nm}^{**}(a)\leq\vafu_{\nm}(a)=\zat{\teob}{a}-\dufu_{\nm}(\teob)\leq\dufu_{\nm}^*(a)
 \]
 are tight and $\teob$ is a dual solution for $a$.
\end{Proof}

\begin{corollary}\label{C:prisol.integ}
   A sufficient condition for $\vafu_{\nm}\not\equiv+\infty$ is
   $f_{\vte}\in\calG$ (integrability of $\moma f_\vte$) for some
   $\vte\in\Theta_{\nm}$ with $\dufu_{\nm}(\vte)$ finite.
\end{corollary}

\begin{remark}\rm
 The hypotheses in Lemma~\ref{L:basic} may hold also when the
 \PCQ is not valid for $a$, see Example~\ref{Ex:zero.density}.
 If, however, $\nm'_+(\cdot,0)\equiv-\infty$ and in particular
 if $\nm$ is essentially smooth, then all the functions $f_{\vte}$
 in $\FF_\nm$ are positive $\mu$-a.e., and the assumption
 $f_{\vte}\in\calG_a$ does imply $a\in\ri{\dom{\vafu_{\nm}}}$,
 by Corollary~\ref{C:zerooutside} and Lemma~\ref{L:suppg}. For
 $a\in\R^d$ satisfying the \PCQ, the assertion of Lemma~\ref{L:basic}
 admits a conversion, see Lemma~\ref{L:ainri}. The assumption that
 $\dufu_{\nm}(\vte)$ is finite is essential in Lemma~\ref{L:basic}
 and Corollary~\ref{C:prisol.integ}, see Example~\ref{Ex:noPCQ.DCQ}.
 That assumption is automatically satisfied when $\dufu_{\nm}$
 is proper, since $\vte\in\Theta_{\nm}\subseteq\dom{\dufu_{\nm}}$ implies
 that $\dufu_{\nm}(\vte)<+\infty$.
\end{remark}

% EFFECTIVE DUAL SOLUTION

\paragraph{4.B.} This subsection introduces effective dual solutions in general.

\begin{lemma}\label{L:dualsol}
   Under the \DCQ, if $\vte\in\R^d$ is a dual solution for $a\in\R^d$
   then $\vte\in\Theta_{\nm}$ and
   \begin{equation}\label{E:dirder}
        \inttt{\TT}\:\zat{\teob-\vte}{\moma}f_\vte\:{\mathrm d}\mu\geq\zat{\teob-\vte}{a}\,,
            \qquad \teob\in\dom{\dufu_{\nm}}\,,
   \end{equation}
   where the integrals are finite.
\end{lemma}

\begin{Proof}
 By the assumption, $\dufu_{\nm}^*(a)$ is finite and equals
 $\zat{\vte}{a}-\dufu_{\nm}(\vte)$ whence $\dufu_{\nm}(\vte)$
 is finite and $\dufu_{\nm}$ proper. If $\teob\in\dom{\dufu_{\nm}}$
 and $0<t<1$ then
 \[
    \zat{\vte}{a}-\dufu_{\nm}(\vte)=\dufu_{\nm}^*(a)
    \geq\zat{t\teob+(1-t)\vte}{a}-\dufu_{\nm}(t\teob+(1-t)\vte)
 \]
 by the definition of conjugation. This implies $\tfrac{1}{t}
 [\dufu_{\nm}(\vte+t(\teob-\vte))-\dufu_{\nm}(\vte)]\geq\zat{\teob-\vte}{a}$.
 Then, the limit of the left-hand side as $t\downarrow0$ is at least
 $\zat{\teob-\vte}{a}$. By Lemma~\ref{L:limint}, the limit is finite,
 equals the integral \eqref{E:integr}, and the set
 $\{\zat{\vte}{\moma}=\nm'(\cdot,+\infty)\,,
                    \,\zat{\teob}{\moma}\neq\zat{\vte}{\moma}\}$
 is $\mu$-negligible. Therefore, if $\zat{\teob}{\moma}<\nm'(\cdot,+\infty)\,[\mu]$
 then also $\zat{\vte}{\moma}<\nm'(\cdot,+\infty)\,[\mu]$, by~Lemma~\ref{L:domdufu}.
 By~\eqref{E:DCQ} and the assumption $\Theta_\nm\neq\pmn$, it follows that
 $\vte\in\Theta_\nm$. By Remark~\ref{R:limint}, the integral \eqref{E:integr}
 rewrites to $\int_{\TT}\:\zat{\teob-\vte}{\moma}f_\vte\:{\mathrm d}\mu$ and
 the assertion follows.
\end{Proof}

\begin{corollary}\label{C:nondep}
    Under the \DCQ, if $\vte,\teob$ are dual solutions for $a\in\R^d$
    then $f_\vte=f_\teob$ ($\mu$-a.e.) and
    \[
        \inttt{\TT}\:\zat{\teob-\vte}{\moma}f_\vte\:{\mathrm d}\mu=\zat{\teob-\vte}{a}\,.
    \]
\end{corollary}

\begin{Proof}
 Summing ineq.~\eqref{E:dirder} and its instance with $\vte$ and $\teob$ interchanged,
 \[
    \inttt{\TT}\:\big[\zat{\vte}{\moma}-\zat{\teob}{\moma}\big]
                                \big[f_\vte-f_\teob\big]\:{\mathrm d}\mu\leq0\,.
 \]
 The definition \eqref{E:fte} of $f_\vte$ and monotonicity of the functions
 $(\conm)'(\ttt,\cdot)$, $\ttt\in\TT$, imply that the product in the
 integral is nonnegative $\mu$-a.e. Hence, the product vanishes
 $\mu$-a.e.\ and $f_\vte=f_\teob$ holds by \eqref{E:fte}.
 In turn, ineq.~\eqref{E:dirder} is tight.
\end{Proof}

\begin{remark}\label{R:efduso}\rm
 Assuming the \DCQ and existence of dual solutions for $a\in\R^d$,
 Corollary~\ref{C:nondep} implies that each dual solution $\vte$ for
 $a$ gives rise to the same function $f_\vte$. The unique function
 defined thereby is denoted by $\efa$ and referred to as the
 \emph{effective dual solution} for $a$.
 In Subsection~1.E., $\efa$ appeared assuming additionally
 to the \DCQ also the \PCQ for $a$, which is a sufficient but not always
 necessary condition for existence of dual solutions, see Remark~\ref{R:dualatt}.
 Corollary~\ref{C:nondep} goes beyond the situation $a\in\ri{\dom{\vafu_\nm}}$,
 see Example~\ref{Ex:zero.density}.
\end{remark}

\begin{lemma}\label{L:ainri}
   Assuming the \PCQ holds for $a\in\R^d$, the primal solution $g_a$
   exists if and only if $\Theta_{\nm}\neq\pmn$ and the moment vector
   of $\efa$ exists and equals $a$. This takes place if and only if
   $\calG_a$ intersects $\FF_\nm$. In this case, $g_a=\efa$ and
   $\calG_a\cap\FF_\nm=\{g_a\}$.
\end{lemma}

\begin{Proof}
 By Lemma~\ref{L:pri=du}, a dual solution $\vte$ for $a$ exists
 and the primal value $\vafu_{\nm}(a)$ coincides with the dual one
 $\dufu_{\nm}^*(a)=\zat{\vte}{a}-\dufu_{\nm}(\vte)$.
 Hence, $\dufu_{\nm}$ is proper.

 If the primal solution $g_a\in\calG_a$ exists, thus
 $\fuc_\nm(g_a)=\vafu_\nm(a)$, then ineq.~\eqref{E:Fen.int}
 applies to $g_a$ and $\vte$. It rewrites to
 $\fuc_\nm(g_a)+\dufu_\nm(\vte)\geq\zat{\vte}{g_a}$.
 It follows that this inequality is tight whence
 ineq.~\eqref{E:Fen} is tight for $\mu$-a.a.\ $\ttt\in\TT$.
 In such a case, $g_a(z)$ is equal to the derivative
 of $\conm(\ttt,\,.\,)$ at $\zat{\vte}{\moma(\ttt)}$
 \cite[Theorem 23.5]{Rock}. Therefore, $\vte\in\Theta_\nm$
 and $g_a=f_\vte$. This implies that the moment vector of $\efa=f_\vte$
 exists and equals $a$, and also that $\calG_a$ and $\FF_\nm$
 intersect. Thus, both conditions for existence are~necessary.

 If $\Theta_{\nm}\neq\pmn$ then the dual solution $\vte$ belongs
 to this set by Lemma~\ref{L:dualsol}. If also the moment vector
 of the effective dual solution $\efa=f_\vte$ exists and equals
 $a$ then $\calG_a$ intersects~$\FF_\nm$ in $f_\vte$. More generally,
 if the intersection contains $f_\teob$ for some $\teob\in\Theta_{\nm}$
 then $\dufu_{\nm}(\teob)$ is finite because $\dufu_{\nm}$ is proper.
 Lemma~\ref{L:basic} implies that the primal solution $g_a$
 exists and equals $\efa$, and $\calG_a\cap\FF_\nm=\{g_a\}$.
 Thus, both conditions for existence are sufficient and the last
 assertion holds.
\end{Proof}

\begin{proposition}\label{P:dualsol}
   If the \DCQ holds and $\dom{\dufu_{\nm}}$ has nonempty interior
   then each effective dual solution belongs to $\calG$ (has a moment vector).
   If the \DCQ holds and the dual problem for $a\in\R^d$ has a solution
   in the interior of $\dom{\dufu_{\nm}}$ then the primal solution $g_a$
   exists and equals~$\efa$.
\end{proposition}

\begin{Proof}
 If for some $a\in\R^d$ a dual solution $\vte$ exists then
 $\dufu_{\nm}$ is proper. The \DCQ and Lemma~\ref{L:dualsol}
 imply $\vte\in\Theta_\nm$. Since the interior is nonempty, the set
 of $\teob-\vte\in\R^d$ with $\teob\in\dom{\dufu_{\nm}}$ has full
 dimension~$d$. This and finiteness of the integrals in \eqref{E:dirder}
 imply that $\moma f_\vte$ is integrable, thus the first assertion holds.
 If additionally the dual solution $\vte$ is in the interior then the inequalities
 in~\eqref{E:dirder} turn into equalities whence $\efa=f_\vte\in\calG_a$.
 By Lemma~\ref{L:basic}, the primal solution $g_a$ exists and equals $\efa$.
\end{Proof}

\begin{corollary}\label{C:prisolclassic}
   If the \DCQ holds and $\dufu_{\nm}$ is essentially smooth then the primal
   solution $g_a$ exists for $a\in\R^d$ whenever a dual solution does,
   and then $g_a=\efa$. In particular, the primal solution exists and
   equals the effective dual solution for each $a\in\ri{\dom{\vafu_{\nm}}}$.
\end{corollary}

\noindent
 For essential smoothness of $\dufu_{\nm}$, see Lemma~\ref{L:essmooth},
 a trivial sufficient condition is that $\dufu_{\nm}$ is proper with
 its effective domain open. When $\dufu_{\nm}$ is essentially smooth,
 Corollary~\ref{C:neededinBoston} gives a necessary and sufficient
 condition for the \DCQ. Note that the essential smoothness
 of $\dufu_\nm$ is not related to that of the integrand \nm.

\begin{Proof}
 Assuming a dual solution $\vte$ for $a\in\R^d$ exists, $f_\vte$ has a moment
 vector  by Proposition~\ref{P:dualsol}. Thus, for $\teob\in\Theta_\nm$
 the directional derivative  $\dufu'_{\nm}(\vte;\teob-\vte)$ is finite,
 see Lemma~\ref{L:limint}. In particular, this holds for $\teob$ in the
 interior of $\dom{\dufu_{\nm}}$, hence the essential smoothness of
 $\dufu_{\nm}$ implies by Lemma~\ref{L:essmooth} that $\vte\in\Theta_\nm$
 is not on the boundary of $\dom{\dufu_{\nm}}$. Having $\vte$ in the interior,
 the second part of Proposition~\ref{P:dualsol} gives that $g_a$ exists and
 equals $\efa$. The last assertion follows by Lemma~\ref{L:pri=du}.
\end{Proof}

% CORRECTION

\paragraph{4.C.} This subsection introduces the correction term for the
functions $g\in\calG^+$
 whose moment vector belongs to $\ri{\dom{\vafu_{\nm}}}$, provided the \DCQ
 holds and $\vafu_\nm>-\infty$.

 Recall the function $\lad$, $\ga\in\varGamma$, given by~\eqref{E:upsi}
 in Subsection~2.D. For $\vte\in\Theta_{\nm}$ let $\ladnin^{\;\vte}$
 denote the function given for $s\geq0$ and $\ttt\in\TT$
 by\index{$\lad$, $\ladnin$~~integrand behind the correction}\label{integrand behind the correction 2}
 \begin{equation}\label{E:upsi.vte}
   \ladnin^{\;\vte}(\ttt,s)\triangleq
                  \ladarg{\nm(\ttt,\cdot)}{s,\zat{\vte}{\moma(\ttt)}}
                     =\big[\nm'_{\sg{s-f_\vte(\ttt)}}
                         (\ttt,f_\vte(\ttt))-\zat{\vte}{\moma(\ttt)}\big]\,
                       \big[s-f_\vte(\ttt)\big]
 \end{equation}
 if $\zat{\vte}{\moma(\ttt)}<\nm'(\ttt,+\infty)$ and
 $\ladnin^{\;\vte}(\ttt,s)\triangleq0$ otherwise. Since
 $\vte\in\Theta_{\nm}$ the latter case is $\mu$-negligible.
 The function $\ladnin^{\;\vte}$ is nonnegative. It
 vanishes if $\nm$ is essentially smooth. If $\nm$ is
 differentiable then $\ladnin^{\;\vte}(\ttt,s)=
 |\nm'_+(\ttt,0)-\zat{\vte}{\moma(\ttt)}|_+\cdot s$.

 For a \calT-measurable function $g\ge 0$ let
 \begin{equation}\label{E:corre}
        \Dore_{\nm}^\vte(g)\triangleq
            \inttt{\TT}\:\ladnin^{\;\vte}(\ttt,g(\ttt))\:\mu({\mathrm d}\ttt)\,,
                \qquad \vte\in\Theta_{\nm}\,.
 \end{equation}

\begin{remark}\label{R:Dore}\rm
  The \calT-measurability of the nonnegative function in
  the integral  \eqref{E:corre} follows from the identity
  in the proof of Lemma~\ref{L:KEYID}. By \eqref{E:upsi.vte},
  $\Dore_{\nm}^\vte(f_\vte)=0$, and if $\nm$ is essentially
  smooth then $\Dore_{\nm}^\vte\equiv0$ on $\calG^+$.
  If $\nm$ is differentiable then for $g\ge 0$
  \begin{equation}\label{E:corrediff}
     \Dore_{\nm}^\vte(g)=\inttt{\TT}\:
             |\nm'_+(\cdot,0)-\zat{\vte}{\moma}|_+\cdot g\,{\mathrm d}\mu\,,
                                    \qquad \vte\in\Theta_{\nm}\,,
  \end{equation}
  hence in this case $\zat{\vte}{\moma}\ge\nm'_+(\cdot,0)\; [\mu]$ is a
  sufficient condition for $\Dore_{\nm}^\vte(g)=0$.
\end{remark}

\begin{lemma}\label{L:PCQ.Pyth}
   Assuming the \DCQ, if $\vte,\teob$ are dual solutions for $a\in\R^d$ then
   \[
    \Dore_{\nm}^\vte(g)=\Dore_{\nm}^\teob(g),\quad g\in\calG_a^+\,.
   \]
\end{lemma}

\begin{Proof}
 By Lemma~\ref{L:dualsol}, $\vte,\teob$ are in $\Theta_\nm$. By Corollary~\ref{C:nondep},
 two such solutions $\vte,\teob$ give $f_\vte=f_\teob\,[\mu]$ and
 $\int_{\TT}\:\zat{\teob-\vte}{\moma}f_\vte\:{\mathrm d}\mu=\zat{\teob-\vte}{a}$.
 Eq.~\eqref{E:upsi.vte} implies that
  \[
     \ladnin^{\;\vte}(\ttt,s)=\ladnin^{\;\teob}(\ttt,s)+\zat{\teob-\vte}
         {\moma(\ttt)}\,[s-f_\vte(\ttt)]\,,\qquad \text{for $s\geq0$ and $\mu$-a.a. $\ttt\in\TT$}\,.
  \]
 For $g\in\calG_a^+$ the assertion follows by substituting $s=g(\ttt)$ and integrating.
\end{Proof}

\bigskip

 The \emph{correction functional} $\Cot_\nm$ alluded to
 in~\eqref{E:Pyth} is defined, temporarily, for the\index{$\Cot_{\nm}$~~correction functional}\label{correction functional}
 functions $g\in\calG^{+}$ whose moment vectors belong to
 $\ri{\dom{\vafu_{\nm}}}$, by
 \begin{equation}\label{E:Corr}
    \Cot_\nm(g)\triangleq \Dore_{\nm}^\vte(g)\quad
            \text{where $\vte\in\Theta_{\nm}$ is any dual solution for}\;
                \mbox{$\int_{\TT}$}\:\moma g\,{\mathrm d}\mu\,,
 \end{equation}
 provided that $\vafu_\nm$ is proper and the \DCQ holds. Here,
 a dual solution exists by Lemma~\ref{L:pri=du} and the definition
 does not depend on its choice by Lemma~\ref{L:PCQ.Pyth}.
 This definition of the correction functional
 is further extended in Section~\ref{S:main}.

% PYTH

\paragraph{4.D.} The key lemma of this subsection is formulated as follows.

\begin{lemma}\label{L:KEYID}
  For $a\in\R^d$, $g\in\calG_a^+$ and $\vte\in\Theta_{\nm}$
  \[
   \fuc_{\nm}(g)=\zat{\vte}{a}-\dufu_{\nm}(\vte)+
               \Bre{g,f_\vte}+\Dore_{\nm}^\vte(g)\,.
  \]
\end{lemma}

\begin{Proof}
 For $\ttt\in\TT$, Lemma~\ref{L:keyid} is applied to $\nm(\ttt,\cdot)$ in the role
 of $\ga$, with $s=g(\ttt)$ and $r=\zat{\vte}{\moma(\ttt)}$. It follows
 that if $\zat{\vte}{\moma(\ttt)}<\nm'(\ttt,+\infty)$, which holds for $\mu$-a.e.\
 $\ttt\in\TT$ by~\eqref{E:DCQ}, then
 \[
    \nm(\ttt,g(\ttt)) +\conm(\ttt,\zat{\vte}{\moma(\ttt)})
    =\zat{\vte}{\moma(\ttt)} g(\ttt)+\trnnin(\ttt,g(\ttt),f_\vte(\ttt))+
       \ladnin^{\;\vte}(\ttt,g(\ttt))\,.
 \]
 The assertion is obtained by integration since $\dufu_{\nm}(\vte)<+\infty$.
\end{Proof}

\bigskip

 Below, the generalized Pythagorean identity \eqref{E:Pyth} is formulated,
 under\index{generalized Pythagorean identity}\label{generalized Pythagorean identity 2}
 restrictive assumptions alleviated later in Theorem~\ref{T:main}.

\begin{lemma}\label{L:duPyth}
  Assuming the \PCQ for $a\in\R^d$ and the \DCQ,
  \[
        \fuc_{\nm}(g)=\vafu_\nm(a)+\Bre{g,\efa}+\,\Cot_\nm(g),
            \qquad g\in\calG_{a}^+\,.
  \]
\end{lemma}

\begin{Proof}
 By Lemmas~\ref{L:pri=du} and~\ref{L:dualsol}, a dual solution
 $\vte$ for $a$ exists, it belongs to~$\Theta_{\nm}$ and
 $\zat{\vte}{a}-\dufu_{\nm}(\vte)=\dufu_{\nm}^{*}(a)=\vafu_\nm(a)$.
 It suffices to invoke Lemma~\ref{L:KEYID},
 using Remark~\ref{R:efduso} and the definition~\eqref{E:Corr}.
\end{Proof}

\bigskip

 The main result of this section is based on the hypotheses that the duality gap
 between the primal and dual values is zero and the dual value is attained.
 By Lemma~\ref{L:pri=du}, the \PCQ is a sufficient condition for this. However,
 it is not necessary, see Example~\ref{Ex:zero.density}. In general, it is difficult
 to recognize whether the gap is zero and this problem is not addressed here.

\begin{theorem}\label{T:geprisol}
    For $a\in\R^d$ let the duality gap be zero and the dual value be attained.
    Then, the  generalized primal solution for $a$ exists if and only if
    the \DCQ holds, in which case $\gea=\efa$.
\end{theorem}

\begin{Proof}
 The hypotheses imply that $\vafu_{\nm}(a)=\dufu_{\nm}^*(a)$
 which is equal to $\zat{\vte}{a}-\dufu_{\nm}(\vte)$ for some dual solution
 $\vte\in\R^d$. Let $g_n$ be a sequence in $\calG_a^+$ such that $\fuc_{\nm}(g_n)$
 converges to the finite primal value $\vafu_{\nm}(a)$. Then the functions
 \[
    h_n\colon\ttt\mapsto\nm(\ttt,g_n(\ttt))+\conm\big(\ttt,\zat{\vte}{\moma(\ttt)}\big)
    -\zat{\vte}{\moma(\ttt)}g_n(\ttt)
 \]
 are nonnegative, $\calT$-measurable and their integrals
 $\fuc_{\nm}(g_n)+\dufu_{\nm}(\vte)-\zat{\vte}{a}$ go to zero.
 Then, going to a subsequence if necessary, $h_n\to 0,\;\mu$-a.e.
 If $\Theta_\nm=\pmn$ then Corollary~\ref{C:Thetaindom} implies that
 $\zat{\vte}{\moma(\ttt)}=\nm'(\ttt,+\infty)$ for $\ttt$ in a subset
 $Y\in\calT$ of $\asli$ of positive $\mu$-measure, and thus
 \[
    \big[\nm'(\ttt,+\infty)\,g_n(\ttt)-\nm(\ttt,g_n(\ttt))\big]\to
    \conm\big(\ttt,\nm'(\ttt,+\infty)\big)\,, \qquad
    \text{for $\mu$-a.e.\ $\ttt\in Y\pdm\asli$}\,.
 \]
 Since $Y\pdm\asli$ it follows that $g_n$ goes to $+\infty$ $\mu$-a.e.\
 on $Y$. Therefore, the sequence $g_n$ is not convergent locally in measure, and
 thus the generalized primal solution for $a$ does not exist.

 Assuming the \DCQ holds, the dual solution $\vte$ belongs to $\Theta_\nm$ by
 Lemma~\ref{L:dualsol}. Thus, \mbox{$f_\vte=\efa$} by Remark~\ref{R:efduso}.
 Lemma~\ref{L:KEYID} implies that $\fuc_{\nm}(g)\geq\vafu_{\nm}(a)+\Bre{g,\efa}$,
 \mbox{$g\in\calG_{a}^+$}. For any sequence $g_n$ in $\calG_a^+$ with
 $\fuc_\nm(g_n)\to \vafu_{\nm}(a)$ necessarily $\Bre{g_n,\efa}\to0$.
 By Corollary~\ref{C:inmeas}, $g_n\rightsquigarrow \efa$. This proves
 that the generalized primal solution $\gea$ exists and equals~$\efa$.
\end{Proof}

\bigskip

 Let the Bregman closure of $\calG_a$ be defined as the set of \calT-measurable
 functions $h$ such that $\Bre{g_n,h}\to 0$ for some sequence $g_n$ in $\calG_a$.
 In the `irregular' situation when $\calG_a$ and $\FF_\nm$ are disjoint, $\FF_\nm$
 can still intersect the closure. For example, $\efa\in\FF_\nm$ belongs to the
 closure if the duality gap is zero and the dual value for $a$ is attained, by
 the last part of the above proof. The following assertion provides the converse under
 a regularity condition. Let $\Theta_\nm^+$\index{$\Theta_\nm^+$~~special subset of $\Theta_\nm$}\label{special subset of theta}
 denote the set of those $\vte\in\Theta_\nm$ for which $\zat{\vte}{\moma}\ge\nm_+'(\cdot,0)$ $[\mu]$.
 If \nm is essentially smooth then $\Theta_\nm^+=\Theta_\nm$. In general, $\Theta_\nm^+$
 cannot be replaced by $\Theta_\nm$ in Proposition~\ref{P:Breclo},
 see Example~\ref{E:int.big}.

\begin{proposition}\label{P:Breclo}
   Let $\nm$ be differentiable and $a\in\R^d$. If $\vte\in\Theta_\nm^+$
   has $\dufu_{\nm}(\vte)$ finite and $f_\vte$ belongs to the Bregman
   closure of $\calG_a$ then $\vafu_\nm(a)=\dufu_\nm^*(a)$, $\vte$
   is a dual solution for $a$ and $f_\vte=\efa$.
\end{proposition}

\begin{Proof}
 The hypotheses on $\nm$ and $\vte$ imply that $\Dore_{\nm}^\vte$ vanishes on $\calG^+$,
 on account of \eqref{E:corrediff} in Remark~\ref{R:Dore}. By assumption, there
 exists a sequence $g_n$ in $\calG_a^+$ with $\Bre{g_n,f_\vte}\to 0$. It follows from
 Lemma~\ref{L:KEYID} that
 \[
    \vafu_\nm(a)\le\lim_{n\to\infty}\fuc_\nm(g_n)
        =\zat{\vte}{a}-\dufu_\nm(\vte)\le\dufu_\nm^*(a)\,.
 \]
 Lemma~\ref{L:duleqpri} implies that the above inequalities are tight
 and the assertions follow.
\end{Proof}

%555555555555555555555555555555555555555555555555555555555555555555555555555555
\section{Conic cores}\label{S:cnc}

 A set $C$ in~$\R^d$ is a \emph{cone} if it contains the origin $\bnu$ and
 $tx\in C$ whenever $t>0$ and $x\in C$. The convex/conic hull of~$C$
 is denoted by $\conv{C}/\cone{C}$.\index{$\cnc{Q}$~~conic core of a measure $Q$}\label{conic core of a measure}

 In this section $Q$ typically denotes a $\sigma$-finite Borel measure
 on $\R^d$. A Borel subset of~$\R^d$ is $Q$-\emph{full} if its complement
 has $Q$-measure zero. The intersection of all closed $Q$-full sets in~$\R^d$
 is the \emph{support} $\s{Q}$ of $Q$ and the intersection of all convex,
 closed and $Q$-full sets is the \emph{convex support} $\cs{Q}$ of $Q$.

 The \emph{convex core} $\cc{Q}$ of a probability measure (pm) $Q$ was introduced
 in \cite{Csi.Ma.cc} as the intersection of all $Q$-full convex Borel sets
 in $\R^d$. The concept extends naturally to the $\sigma$-finite  measures
 \cite{Csi.Ma.gmle} since $\cc{Q}$ does not change when $Q$ is replaced
 by a finite measure equivalent to $Q$.  An equivalent definition involves
 the means of probability measures dominated by $Q$, namely by
 \cite[Theorem~3]{Csi.Ma.cc},\label{probmeas}
 \begin{equation}\label{E:cc.repr}
     \cc{Q}=\big\{\mbox{$\int_{\R^d}$}\:x\, P({\mathrm d}x)\colon
            \text{$P$ is a pm with mean and $P\ll Q$}\big\}\,.
 \end{equation}

\begin{definition}\label{D:cnc,cns}
{\rm   The \emph{conic core} $\cnc{Q}$ of a $\sigma$-finite Borel measure~$Q$
   on $\R^d$ is the intersection of the convex, Borel and $Q$-full cones.
   The \emph{conic support} $\cns{Q}$ is the intersection of the convex,
   closed and $Q$-full cones.
}
\end{definition}

\begin{remark}\label{R:zero}\rm
 The conic core is a convex cone, not necessarily $Q$-full.  The conic support
 is a convex, closed and $Q$-full cone. Both are nonempty since they contain
 the origin; they are equal to the singleton $\{\bnu\}$ if and only if
 $\mu(\R^d\setminus\{\bnu\})=0$. The conic core and support do not depend
 on the weight assigned by $Q$ to $\{\bnu\}$. Thus, in Definition~\ref{D:cnc,cns},
 one can admit infinite $Q$-mass at $\bnu$, going slightly beyond
 $\sigma$-finiteness.
\end{remark}

 Some properties of conic cores can be derived also from known facts
 on the convex cores, but a direct self-contained approach
 is preferred in this section.

\begin{lemma}\label{L:csccs}
   $\cl{\cnc{Q}}\,{=}\,\cns{Q}$ \,and\,
   $\ri{\cnc{Q}}\,{=}\,\ri{\cns{Q}}\,{=}\,\ri{\cone{\conv{\s{Q}}}}$.
\end{lemma}

\begin{Proof}
 By definition, $\cnc{Q}\pdm \cns{Q}$ whence $\cl{\cnc{Q}}\pdm \cns{Q}$
 using that $\cns{Q}$ is closed. There is no loss of generality in assuming
 that the dimension of $\{\bnu\}\cup\s{Q}$ is $d$. If $K$ is a convex, Borel
 and $Q$-full cone then $\cl{K}$ is a convex, closed and $Q$-full cone, and
 thus $\cns{Q}\pdm\cl{K}$. By the assumption on dimension, $\ri{\cns{Q}}\pdm K$,
 and, in turn, $\ri{\cns{Q}}$ is contained in $\cnc{Q}$. Hence, $\cns{Q}$
 is contained in $\cl{\cnc{Q}}$. The second assertion is a consequence
 of the first one and $\cns{Q}=\cl{\cone{\conv{\s{Q}}}}$.
\end{Proof}

\bigskip

 A \emph{supporting hyperplane} to a convex cone $K$ is a hyperplane $H$
 intersecting $K$ such that one of the closed half-spaces bordered by $H$
 contains $K$. A \emph{nontrivial} supporting hyperplane does not
 contain $K$. Any supporting hyperplane to a convex cone $K$ contains
 the origin. Thus, there exists $\vte\in\R^d$ nonzero such that
 $H=\{x\colon\zat{\vte}{x}=0\}$ and $K\pdm H\cup H_<$ where
 $H_<=\{x\colon\zat{\vte}{x}<0\}$.

 The restriction of a $\sigma$-finite measure~$Q$ to a Borel set $A\pdm\R^d$
 is denoted by $Q^A$. It is given by $Q^A(B)=Q(A\cap B)$ for every
 $B\subseteq\R^d$ Borel.

\begin{lemma}\label{L:supphyp}
    If $H$ is a supporting hyperplane of \cns{Q} then
    $\cnc{Q}\cap H=\cnc{Q^H}$ and $Q(H\setminus\cl{\cnc{Q}\cap H})=0$.
\end{lemma}

\begin{Proof}
 The hyperplane $H$ contains the origin and $\cns{Q}\pdm H\cup H_<$ as above.
 If $K$ is any convex, Borel and $Q^H$-full cone then $K\cap H$ has the same
 properties. Then, $(K\cap H)\cup H_<$ is a convex, Borel cone which is
 $Q$-full by $\cns{Q}\pdm H\cup H_<$. Hence, $\cnc{Q}\pdm(K\cap H)\cup H_<$
 and, intersecting with~$H$, $\cnc{Q}\cap H\pdm K$. This implies that
 $\cnc{Q}\cap H\subseteq\cnc{Q^H}$. The opposite inclusion holds because
 $\cnc{Q}$ and $H$ contain $\cnc{Q^H}$, by definitions. The first assertion
 and Lemma~\ref{L:csccs} imply that $\cl{\cnc{Q}\cap H}$ equals $\cns{Q^H}$.
 Then, the second assertion follows since this set is $Q^H$-full.
\end{Proof}

\bigskip

 A \emph{face} of a convex set $C\pdm\R^d$ is a nonempty convex subset
 $F\subseteq C$ such that every closed line segment in $C$ with a relative
 interior point in $F$ is contained in $F$. The face is proper if $F\not=C$.
 The relative interiors \riF of the faces $F$ partition the set $C$
 \cite[Theorem~18.2]{Rock}. A face of a convex cone is a convex cone.
 The \emph{smallest face} of a convex cone (the intersection of all faces)
 is either the singleton $\{\bnu\}$ or a linear subspace of $\R^d$.

\begin{lemma}\label{L:face}
   If $F$ is a face of \cnc{Q} then $\cnc{Q^{\clF}}=F$.
\end{lemma}

\begin{Proof}
 Induction on the dimension of \cnc{Q} is employed. If $F=\cnc{Q}$ then the
 assertion follows by Lemma~\ref{L:csccs} using that $\cns{Q}$ is $Q$-full.
 Otherwise, $F$ is a proper face and there exists a nontrivial supporting
 hyperplane $H$ to \cnc{Q} containing $F$ \cite[Theorem~11.6]{Rock}. Then,
 $\cnc{Q}\cap H$ is a proper face of \cnc{Q} containing $F$. Lemma~\ref{L:supphyp}
 implies that $F$ is a face of \cnc{Q^H}. As \cnc{Q^H} has smaller dimension
 than \cnc{Q}, by induction, $\cnc{(Q^H)^{\clF}}=F$.
\end{Proof}

\begin{corollary}\label{C:este}
   $Q(\clF)>0$ for each face $F$ of \cnc{Q}, except perhaps for $F=\{\bnu\}$.
\end{corollary}

\begin{lemma}\label{L:inteincnc}
   If the integral $\int_{\R^d}\:x\,Q({\mathrm d}x)$ exists then it belongs to $\ri{\cnc{Q}}$.
\end{lemma}

\begin{Proof}
 Let $H$ be a supporting hyperplane of \cns{Q}, thus
 $\cns{Q}\pdm H\cup H_<$ where $H$ and $H_<$ are
 parameterized by $\vte$ as above. Denoting the
 integral by~$a$, $\zat{\vte}{a}$ equals
 $\int_{H_<}\:\zat{\vte}{x}\,Q({\mathrm d}x)$.
 If $Q(H_<)=0$ then $\zat{\vte}{a}=0$ whence $a\in H$. Otherwise, the
 supporting hyperplane $H$ is nontrivial and  $\zat{\vte}{a}<0$, since
 $\zat{\vte}{x}<0$ for $x\in H_<$. Thus, $a\in H_<$. It follows that $a$
 belongs to the intersection of all closed halfspaces $ H\cup H_<$
 as above, which equals $\cns{Q}$, but to none of the nontrivial
 supporting hyperplanes of $\cns{Q}$. Therefore, $a\in\ri{\cns{Q}}$
 and the assertion follows by Lemma~\ref{L:csccs}.
\end{Proof}

\begin{corollary}\label{C:intincnc}
   If $P\ll Q$ and the integral $\mbox{$\int_{\R^d}$}\:x\, P({\mathrm d}x)$ exists
   then it belongs to \cnc{Q}.
\end{corollary}

\begin{Proof}
 The integral belongs to $\ri{\cnc{P}}\pdm\cnc{P}\pdm\cnc{Q}$, where the
 latter inclusion follows from $P\ll Q$, by the definition of the conic core.
\end{Proof}

\begin{lemma}\label{L:interepr}
   Each $a\in\ri{\cnc{Q}}$ can be represented as $\mbox{$\int_{\R^d}$}\:x\,P({\mathrm d}x)$
   where $P$ is a finite measure that is dominated by $Q$, has compact
   support, and its $Q$-density takes a finite number of values.
\end{lemma}

\begin{Proof}
 Let $C_Q$ denote the set of points that can be represented
 as the above integral with $P$ having the stated properties.
 By Corollary~\ref{C:intincnc}, $C_Q$ is a convex subcone of $\cnc{Q}$.
 Then, it suffices to show that $\ri{\cnc{Q}}\subseteq \cl{C_Q}$,
 since this implies the assertion $\ri{\cnc{Q}}\subseteq C_Q$.

 By Lemma~\ref{L:csccs}, each $a\in\ri{\cnc{Q}}$ can be represented
 as $\sum_{y\in Y}\;t_y\,y$ where $Y$ is a finite subset
 of \s{Q} and all $t_y$ are positive. Since $Q$  is $\sigma$-finite,
 for any $\vare>0$ and $y\in Y$ there exists a Borel subset
 $A_{\vare,y}$ of the $\vare$-ball $B_y(\vare)$ around $y$
 such that $Q(A_{\vare,y})$ is positive and finite.
 Let $y_\vare=\mbox{$\int_{A_{\vare,y}}$}\:x\,Q({\mathrm d}x)/Q(A_{\vare,y})$.
 Then, each $y_\vare$ belongs to $C_Q$ and $\norm{y_\vare-y}\le\vare$
 because $y_\vare$ is the mean of a pm concentrated on $B_y(\vare)$.
 Therefore, the point $\sum_{y\in Y}\; t_y\,y_\vare$ of $C_Q$ is
 arbitrarily close to $a$ if $\vare$ is sufficiently small.
 It follows that $a\in\cl{C_Q}$.
 \index{$\norm{\cdot}$ Euclidean norm}\label{Euclidean norm}
\end{Proof}

\begin{theorem}\label{T:cnc}
    The conic core $\cnc{Q}$ consists of the integrals
    $\mbox{$\int_{\R^d}$}\:x\,P({\mathrm d}x)$ where $P$ runs
    over all finite measures dominated by $Q$.
\end{theorem}

\begin{Proof}
  One inclusion follows from Corollary~\ref{C:intincnc}. If $a\in\cnc{Q}$
  then $a\in\riF$ for a face $F$ of \cnc{Q}. By Lemma~\ref{L:face},
  $a\in\ri{\cnc{Q^{\clF}}}$. By Lemma~\ref{L:interepr},
  $a=\mbox{$\int_{\R^d}$}\:x\,P({\mathrm d}x)$ for a finite measure
  $P$ dominated by $Q^{\clF}$, and thus by $Q$.
\end{Proof}

\begin{remark}\label{R:cnc.restricted}\rm
 The measures $P$ in Theorem~\ref{T:cnc} can be also restricted
 as in Lemma~\ref{L:interepr}.
\end{remark}

\begin{corollary}\label{C:cc+cnc}
   $\cnc{Q}=\cone{\cc{Q}}$.
\end{corollary}

\begin{Proof}
 The equality follows from \eqref{E:cc.repr}, which is \cite[Theorem~3]{Csi.Ma.cc},
 and Theorem~\ref{T:cnc}.\rule{0cm}{0cm}\quad\rule{3cm}{0cm}\quad\rule{0cm}{0cm}
\end{Proof}

\begin{remark}\label{R:hyperpla}\rm
  The faces of \cnc{Q} and \cc{Q} are not related to each other in general.
  However, if $Q$ is concentrated on a hyperplane that does not contain
  the origin then there is a bijection between the families of faces
  of \cc{Q} and \cnc{Q}, up to the face $\{\bnu\}$ of the latter:
  the faces of \cnc{Q} are the conic hulls of the faces of $\cc{Q}$.
\end{remark}

\begin{remark}\label{R:faces.compar}\rm
  The number of faces of any convex core is at most countable \cite[Theorem~3]{Csi.Ma.cc}.
  This remains true also for the conic cores. In fact, it suffices to prove that if $Q$
  is a pm then $\cnc{Q}=\cc{R}$ for $R=\sum_{n\geq0} Q^{[n]}2^{-n}$ where $Q^{[n]}$ is
  the image of $Q$ under the scaling $x\mapsto nx$. If $a\in\cnc{Q}$ then
  $a=t\int_{\R^d}\:x\,P({\mathrm d}x)$ for $t\geq0$ and a pm $P\ll Q$, by Theorem~\ref{T:cnc}.
  Then $a$ is the mean of a convex combination of $P^{[0]}$ and $P^{[n]}$, $n\geq t$.
  Since $R$ dominates these pm's, $a\in\cc{R}$ by \eqref{E:cc.repr}. In the opposite
  direction, any convex, Borel and $Q$-full cone is also $Q^{[n]}$-full whence $R$-full.
  Therefore, $\cnc{Q}\supseteq\cnc{R}\supseteq\cc{R}$ by definitions.
\end{remark}

%666666666666666666666666666666666666666666666666666666666666666666666666666666666666666
\section{The effective domain of the value function}\label{S:dom}

 Recall that the $\moma$-cone \cnmoma{\mu} of $\mu$ consists of the moment
 vectors $\int_{\TT}\:\moma g \,{\mathrm d}\mu$ of the functions $g\in\calG^+$. The
 $\moma$-cone contains the effective domain of $\vafu_{\nm}$ for each $\nm\in B$.
 In this section, a geometric description of this domain is presented that relies
 upon results on conic cores  from Section~\ref{S:cnc} and pays special attention
 to the relative boundary.

 Let $\mu_\moma$ denote the $\moma$-image of $\mu$. The intuitive meaning of the
 following lemma is that the $\moma$-cone of $\mu$ is equal to the conic core of
 $\mu_\moma$. Conic cores, however, have been defined only for measures on $\R^d$
 which are $\sigma$-finite on $\R^d\sm\{\bnu\}$. As the measure $\mu_\moma$ may
 fail to satisfy this condition, an auxiliary measure $\nu$ is invoked.

\begin{lemma}\label{L:conemoma}
    If $\nu$ is a measure equivalent to $\mu$ and the image $\nu_\moma$
    is $\sigma$-finite on $\R^d\sm\{\bnu\}$ then $\cnmoma{\mu}=\cnc{\nu_\moma}$.
\end{lemma}

\begin{Proof}
 To prove that $\cnc{\nu_\moma}\subseteq\cnmoma{\mu}$, it can be
 assumed that $\nu$ is finite because measures which are
 $\sigma$-finite and equivalent on $\R^d\sm\{\bnu\}$ have the same conic
 core, see~Remark~\ref{R:zero}. Let $h$ be a positive $\mu$-density
 of $\nu$. By Theorem~\ref{T:cnc}, any $a\in\cnc{\nu_\moma}$ can be
 written as  $\mbox{$\int_{\R^d}$}\:x f(x)\,\nu_\moma({\mathrm d}x)$ where
 $f\geq0$ is Borel. If $g(\ttt)= f(\moma (\ttt)) h(\ttt) $ then
 $\mbox{$\int_{\TT}$}\:\moma g \,{\mathrm d}\mu=\mbox{$\int_{\TT}$}\:\moma
 f(\moma) \,{\mathrm d}\nu=a$ which implies $a\in\cnmoma{\mu}$.

 In the opposite direction, suppose $a=\int_{\TT}\:\moma g \,{\mathrm d}\mu$
 for $g\in\calG^+$. There is no loss of generality in assuming that
 $g$ vanishes on the set $\{\moma=\bnu\}$. Denote by $\lambda$ the measure
 with $\mu$-density~$g$. Then, $\{\moma=0\}$ is $\lambda$-negligible
 and $\{\moma\neq\bnu\}$ partitions into at most countably many sets $A_n\in\calT$
 with $\lambda(A_n)$ finite. Let $Q_n$ be the $\moma$-image of $\lambda^{A_n}$
 and $Q$ denote the sum of the measures $Q_n$. By the assumption on $g$,
 $Q(\{\bnu\})=0$. Since
 \[
    +\infty>\mbox{$\int_{\TT}$}\:\norm{\moma}\,{\mathrm d}\lambda
        =\mbox{$\sum_n$}\;\mbox{$\int_{\TT}$}\:\norm{\moma}\,{\mathrm d}\lambda^{A_n}
        =\mbox{$\sum_n$}\;\mbox{$\int_{\R^d}$}\:\norm{x}\,Q_n({\mathrm d}x)\,,
 \]
 the complement of any ball around the origin has finite $Q$-measure.
 Therefore, $Q$ is $\sigma$-finite and $a=\mbox{$\int_{\R^d}$}\:x \,Q({\mathrm d}x)$.
 Since $Q_n\ll\nu_\moma$, it follows that $a\in\cnc{\nu_\moma}$,
 using Corollary~\ref{C:intincnc}.
\end{Proof}

\begin{corollary}\label{C:zerooutside}
  The set $\{\moma\not\in\cl{\cnmoma{\mu}}\}$ is $\mu$-negligible.
\end{corollary}

\begin{Proof}
 By Lemma~\ref{L:conemoma}, it $\nu$ is finite and equivalent to $\mu$
 then $\cl{\cnmoma{\mu}}$ is equal to $\cl{\cnc{\nu_\moma}}$ which is
 $\nu_\moma$-full by Lemma~\ref{L:csccs}. This implies that
 $\moma^{-1}(\cl{\cnc{\nu_\moma}})$ is $\nu$-full and the
 assertion follows.
\end{Proof}

\begin{remark}\rm
  The $\moma$-cone \cnmoma{\mu} can be equivalently defined to consist of the
  moment vectors $\int_{\TT}\:\moma g \,{\mathrm d}\mu$ of the $\mu$-integrable (rather
  than all) functions $g$ from $\calG^+$. This follows from Lemma~\ref{L:conemoma}
  and the first part of its proof. In fact, $f$ can be taken $\nu_\moma$-integrable
  by Theorem~\ref{T:cnc}, and then $g=f(\moma)\cdot h$ is $\mu$-integrable.
\end{remark}

\begin{lemma}\label{L:facemoma}
   If $F$ is a face of \cnmoma{\mu} then $\cnmoma{\mu^{\moma^{-1}(\clF)}}=F$.
\end{lemma}

\begin{Proof}
 This follows by Lemma~\ref{L:face} and Lemma~\ref{L:conemoma}.
\end{Proof}

\begin{lemma}\label{L:suppg}
   The moment vector $\int_{\TT}\:\moma g \,{\mathrm d}\mu$ of a function $g\in\calG^+$
   belongs to a face $F$ of \cnmoma{\mu} if and only if
   $g$ vanishes $\mu$-a.e.\ on $\{\moma\notin\cl{F}\}$.
\end{lemma}

\begin{Proof}
 Since any face contains the origin, there is no loss of generality in assuming
 that $g=0$ on $\{\moma=\bnu\}$. Let $a=\int_{\TT}\:\moma g \,{\mathrm d}\mu$, $\nu$
 be a finite  measure equivalent to $\mu$, and $\lambda$ denote the measure
 with $\mu$-density~$g$. Arguing as in the second part of the proof
 of Lemma~\ref{L:conemoma}, $a=\mbox{$\int_{\R^d}$}\:x \,Q({\mathrm d}x)$ for
 a $\sigma$-finite measure $Q=\sum_n\,(\lambda^{A_n})_\moma\ll\nu_\moma$.
 Then, $\cnc{Q}\pdm\cnc{\nu_\moma}=\cnmoma{\mu}$ by Lemma~\ref{L:conemoma},
 and $a\in\ri{\cnc{Q}}$ by Lemma~\ref{L:inteincnc}.

 It follows that if $a\in F$ then $\cnc{Q}\pdm F$. By Lemma~\ref{L:csccs},
 $\clF$ is $Q$-full.  Hence, $\moma^{-1}(\clF)$ is $\lambda^{A_n}$-full
 for all $n$, and thus  $\lambda$-full. This implies that $g$ vanishes
 $\mu$-a.e.\ on $\{\moma\notin\cl{F}\}$. In the opposite direction,
 the vanishing of $g$ implies that $a$ belongs to $\cnmoma{\mu^{\moma^{-1}(\clF)}}$
 which equals $F$, by Lemma~\ref{L:facemoma}.
\end{Proof}

\begin{lemma}\label{L:domvafu}
    If \dom{\vafu_{\nm}} is nonempty then
    $\ri{\dom{\vafu_{\nm}}}=\ri{\cnmoma{\mu}}$.
\end{lemma}

\begin{Proof}
 As $\dom{\vafu_{\nm}}\pdm\cnmoma{\mu}$ are convex sets, the assertion is
 a consequence of the inclusion $\cnmoma{\mu}\pdm\cl{\dom{\vafu_{\nm}}}$
 which is proved in two steps as follows.

 First, let $a\in\cnmoma{\mu}$ equal $\int_{\TT}\:\moma g \,{\mathrm d}\mu$
 for a function $g\in\calG$ that is everywhere positive. Since $\dom{\vafu_{\nm}}$
 is nonempty, $\fuc_\nm(h)<+\infty$ for some $h\in\calG^+$. Since $\mu$ is
 $\sigma$-finite, there exists a positive integrable function $f$ on \TT.
 Let $Y_n$ denote the set of those $\ttt\in\TT$ that satisfy the inequality
 $\nm(\ttt,g(\ttt))\leq\nm(\ttt,h(\ttt))+nf(\ttt)$. As $g$ and $f$ are positive
 and $\nm(\ttt,h(\ttt))<+\infty$ for $\mu$-a.a.\ $\ttt\in\TT$, the sequence
 $Y_n\in\calT$ increases to a $\mu$-full set. Let $g_n$ equal $g$ on $Y_n$
 and $h$ otherwise. It follows that $g_n\in\calG$, the moments
 $a_n=\int_{\TT}\:\moma g_n\,{\mathrm d}\mu$ converge to $a$, and
 \[
     \fuc_\nm(g_n)=\inttt{\,Y_n}\:\nm(\ttt,g(\ttt))\;\mu({\mathrm d}\ttt)
                +\inttt{\,\TT\sm Y_n}\:\nm(\ttt,h(\ttt))\;\mu({\mathrm d}\ttt)
                \leq\fuc_\nm(h)+n\inttt{\,Y_n}\:f\,{\mathrm d}\mu<+\infty\,.
 \]
 Hence, $a_n\in\dom{\vafu_{\nm}}$, and in turn $a\in\cl{\dom{\vafu_{\nm}}}$.

 Second, let $a$ be the moment vector $\int_{\TT}\:\moma g\,{\mathrm d}\mu$ of some function
 $g\in\calG^+$ that may vanish somewhere. The family $\calG$ contains a positive
 function $f$. As the function $g+\frac1n f$ is positive and belongs to~$\calG$,
 its moment $b_n=\int_{\TT}\:\moma (g+\frac1n f) \,{\mathrm d}\mu$ belongs to the closure
 of \dom{\vafu_{\nm}} by the previous part of the proof. Since $b_n\to a$, this
 completes the proof.
\end{Proof}

\begin{remark}\label{R:Slater}\rm
 The relative interior of $\cnmoma{\mu}$ is equal to the set $\cnmomapl{\mu}$
 of points that are representable as $\int_{\TT}\,\moma g\,{\mathrm d}\mu$ with strictly
 positive $g\in\calG^+$. Indeed, $\cnmomapl{\mu}$ is a convex subset of
 \ri{\cnmoma{\mu}} by Lemma~\ref{L:suppg}. Arguing as in the second part
 of the proof of Lemma~\ref{L:domvafu}, the closure of $\cnmomapl{\mu}$
 contains \cnmoma{\mu} whence $\ri{\cnmoma{\mu}}\subseteq\cnmomapl{\mu}$.
\end{remark}

\smallskip
 For a face $F$ of $\cnmoma{\mu}$, let\index{$\Crc{F,\nm}$~~shift in modified problems}\label{shift in modified problems 2}
 \[
    \Crc{F,\nm}\triangleq\inttt{\{\moma\notin\cl{F}\}}\:\nm(\cdot,0)\,{\mathrm d}\mu\,,
 \]
 and $\Fa{\nm}$\index{$\Fa{\nm}$~~family of faces of $\cnmoma{\mu}$}\label{family of faces of}
 denote the family of the faces $F$ such that $\Crc{F,\nm}<+\infty$.
 By Corollary~\ref{C:zerooutside}, $\cnmoma{\mu}\in\Fa{\nm}$.
 If $F\pdm G$ are faces of $\cnmoma{\mu}$ and $F$ belongs to $\Fa{\nm}$ then so
 does also $G$. In particular, $\Fa{\nm}$ contains all faces of $\cnmoma{\mu}$
 if and only if the smallest face of $\cnmoma{\mu}$ belongs to $\Fa{\nm}$.

\begin{theorem}\label{T:domain}
    If $\dom{\vafu_{\nm}}$ is nonempty then it is equal to
    $\bigcup_{F\in\Fa{\nm}}\,\riF$.
\end{theorem}

\begin{Proof}
 In this proof the notation, Definition~\ref{D:F-problems} and Lemma~\ref{L:modif}
 from the beginning of Section~\ref{S:main} are employed. Supposing $\dom{\vafu_{\nm}}
 \neq\pmn$, there exists $g\in\calG^+$ such that $\fuc_{\nm}(g)<+\infty$. If $F$ is
 a face of $\cnmoma{\mu}$ then $g\in\calG_{F}^+$ and $\fuc_{F,\nm}(g)<+\infty$.
 Denoting $\int_{\{\moma\in\cl{F}\}}\:\moma g\,{\mathrm d}\mu$ by $a$, the function~$g$ is
 in $\calG_{F,a}^+$, and $\vafu_{F,\nm}(a)<+\infty$, see Definition~\ref{D:F-problems}.
 Thus, $\dom{\vafu_{F,\nm}}$ is nonempty, and by Lemma~\ref{L:domvafu} it contains
 $\ri{\cnmoma{\mu^{\{\moma\in\cl{F}\}}}}$ that equals $\riF$ by Lemma~\ref{L:facemoma}.
 Hence $\vafu_{F,\nm}<+\infty$ on $\riF$. If $F\in\Fa{\nm}$ then \eqref{E:prevod3}
 from Lemma~\ref{L:modif} implies that $\vafu_{\nm}<+\infty$ on $\riF$. This
 proves that $\dom{\vafu_{\nm}}$ contains the union.

 Conversely, if $\dom{\vafu_{\nm}}$ intersects a face $F$ of $\cnmoma{\mu}$
 then $\fuc_{\nm}(g)<+\infty$ for some $g\in\calG$ with the moment vector
 $\int_{\TT}\:\moma g \,{\mathrm d}\mu \in F$. By Lemma~\ref{L:suppg},
 $g=0\,[\mu]$ on $\{\moma\notin\cl{F}\}$, and thus $\fuc_{\nm}(g)<+\infty$
 implies that $F\in\Fa{\nm}$. Since $\dom{\vafu_{\nm}}$ is
 a subset of \cnmoma{\mu}, it is contained in the union.
\end{Proof}

\begin{corollary}\label{C:domvafu}
   The effective domain of $\vafu_{\nm}$ is closed to positive multiples.
\end{corollary}

\begin{corollary}\label{C:dom.fini}
   A sufficient condition for $\dom{\vafu_{\nm}}=\cnmoma{\mu}$ is
   $\int_\TT\nm(\cdot,0)\,{\mathrm d}\mu<+\infty$. If $\{\bnu\}$ is a face of
   $\cnmoma{\mu}$ and $\{\moma=\bnu\}$ is $\mu$-negligible then this
   condition is necessary, as well. If the integral equals $-\infty$
   then $\vafu_{\nm}=-\infty$ on $\cnmoma{\mu}$.
\end{corollary}

\begin{Proof}
 The first assertion follows from Theorem~\ref{T:domain}, for the hypothesis
 implies that each face of $\cnmoma{\mu}$ belongs to $\Fa{\nm}$ and $\vafu_{\nm}$
 is not $+\infty$ at the origin $\bnu$. By Theorem~\ref{T:domain}, the equality
 $\dom{\vafu_{\nm}}=\cnmoma{\mu}$ implies that the smallest face of $\cnmoma{\mu}$
 belongs to $\Fa{\nm}$. If this smallest face is the singleton $\{\bnu\}$ and
 $\mu(\{\moma=\bnu\})=0$ then $\Crc{\{\bnu\},\nm}=\int_\TT\nm(\cdot,0)\,{\mathrm d}\mu$,
 and the second assertion follows. If the integral equals $-\infty$ then
 $\vafu_{\nm}(\bnu)=-\infty$, and the third assertion follows by convexity
 of $\vafu_{\nm}$.
\end{Proof}

\begin{remark}\label{R:momas*}\rm
 The hypotheses of the second assertion of Corollary~\ref{C:dom.fini},
 guaranteeing the necessity, are valid under the moment assumption.
 In fact, if \eqref{E:ma} holds with  $\teob\in\R^d$ then $\cnmoma{\mu}$
 is contained in the cone  $\{\bnu\}\cup\{x\in\R^d\colon\zat{\teob}{x}>0\}$,
 due to Lemma~\ref{L:conemoma}. Hence, $\{\bnu\}$ is the smallest face of
 $\cnmoma{\mu}$. The second hypothesis $\mu(\{\moma=\bnu\})=0$ follows
 directly from \eqref{E:ma}.
\end{remark}

\begin{corollary}\label{C:dom.open}
   If $\int_{\{\moma\notin H\}}\nm(\cdot,0)\,{\mathrm d}\mu=+\infty$ for each
   nontrivial supporting hyperplane $H$ of the cone $\cnmoma{\mu}$ then
   $\dom{\vafu_{\nm}}$ is either empty or equals $\ri{\cnmoma{\mu}}$.
   Moreover, this condition is necessary for
   $\dom{\vafu_{\nm}}=\ri{\cnmoma{\mu}}$.
\end{corollary}

\begin{Proof}
 Each proper face of $\cnmoma{\mu}$ is contained in a nontrivial supporting
 hyperplane, hence the hypothesis implies that no proper face of $\cnmoma{\mu}$
 belongs to $\Fa{\nm}$. By Theorem~\ref{T:domain}, the first assertion follows.
 Further, the intersection of $\cnmoma{\mu}$ with a nontrivial supporting
 hyperplane $H$ is a proper face $F$. If $\dom{\vafu_{\nm}}=\ri{\cnmoma{\mu}}$,
 Theorem~\ref{T:domain} implies that $F$ does not belong belong to $\Fa{\nm}$,
 thus $\Crc{F,\nm}=+\infty$. By Lemmas~\ref{L:supphyp} and~\ref{L:conemoma},
 this is equivalent to $\int_{\{\moma\notin H\}}\nm(\cdot,0)\,{\mathrm d}\mu=+\infty$.
\end{Proof}

\begin{remark}\label{R:dom.open}\rm
  If $\beta$ is autonomous, $\beta(\ttt,t)=\ga(t)$ for $\ttt\in\TT$ and $t\in\R$,
  then Corollary~\ref{C:dom.fini} states that $\dom{\vafu_{\ga}}$ coincides with
  \cnmoma{\mu} whenever $\ga(0)$ and $\mu$ are finite, or $\ga(0)\leq0$. If
  $\ga(0)=+\infty$ then Corollary~\ref{C:dom.open} gives that  $\dom{\vafu_{\ga}}$
  is either empty or equals $\ri{\cnmoma{\mu}}$, as observed
  in~\cite[Section~3]{Csi.Ma.minentfu}.
\end{remark}

%7777777777777777777777777777777777777777777777777777777777777777777777777777
\section{Dispensing with the PCQ in the primal problem}\label{S:main}

 In this section, the primal problem is studied for the vectors $a\in\R^d$
 with a finite value $\vafu_{\nm}(a)$. The \PCQ is not assumed. Recall
 that the primal and dual problems are constructed from three objects:
 the measure $\mu$, the moment mapping $\moma$ and the integrand $\nm\in B$.
 The notation has not made explicit the dependence on $\mu$ and~$\moma$.
 In this section,  $\mu$ will be replaced by its restriction to
 $\{\moma\in\clF\}$ where $F$ is a face of the $\moma$-cone
 \cnmoma{\mu}. To indicate this restriction, the letter $F$
 is added to indices.

 Correspondingly, for a face $F$ of \cnmoma{\mu} let $\calG_F$ denote
 the linear space  of the \calT-measurable functions $g\colon\TT\to\R$
 such that $\moma g$ is $\mu$-integrable on $\{\moma\in\clF\}$, and
 \[
    \calG_{F,a}\triangleq
        \big\{g\in\calG_F\colon \mbox{$\int_{\{\moma\in\clF\}}$}\,\moma g\,{\mathrm d}\mu=a\big\}\,,
            \qquad a\in\R^d\,.
 \]
 Let $\calG_{F}^+$/$\calG_{F,a}^+$ denote the set of nonnegative
 functions in  $\calG_{F}$/$\calG_{F,a}$.

\begin{definition}\label{D:F-problems}
{\rm   For a face $F$ of \cnmoma{\mu} and $a\in\R^d$, the minimization in
   \[
    \vafu_{F,\nm}(a)\triangleq{\inf}_{g\in\calG_{F,a}^+}\: \fuc_{F,\nm}(g)
        \quad \text{where}\quad
    \fuc_{F,\nm}(g)\triangleq \inttt{\{\moma\in\clF\}}\:\nm(\ttt,g(\ttt))\;\mu({\mathrm d}\ttt)
   \]
   is the \emph{$F$-primal problem} and the maximization in
   \[
    \dufu_{F,\nm}^*(a)\triangleq{\sup}_{\vte\in\R^d}\:[\zat{\vte}{a}-\dufu_{F,\nm}(\vte)]
          \;\;\text{where}\;\;
        \dufu_{F,\nm}(\vte)\triangleq\inttt{\{\moma\in\clF\}}\:
                \conm\big(\ttt,\zat{\vte}{\moma(\ttt)}\big)\:\mu({\mathrm d}\ttt)
   \]
   is the \emph{$F$-dual problem} for $a$. If $\vafu_{F,\nm}(a)$ is finite and the infimum
   is attained then the minimizers can be assumed to vanish outside $\{\moma\in\clF\}$.
   These minimizers define the $\mu$-unique \emph{$F$-primal solution} $g_{F,a}$ for $a$.
   The \emph{generalized $F$-primal solution} $\hat{g}_{F,a}$ is defined likewise.
}
\end{definition}

\begin{remark}\label{R:compar}\rm
  The $F$-primal/$F$-dual problem constructed from $\mu$, $\moma$ and $\nm$
  is identical to the primal/dual problem constructed from $\mu^{\{\moma\in\clF\}}$,
  $\moma$ and $\nm$. Note that if $F=\cnmoma{\mu}$ then $\mu$ does not change
  when restricted to $\{\moma\in\clF\}$, by Corollary~\ref{C:zerooutside}.
\end{remark}

 Two lemmas are sent forward.

\begin{lemma}\label{L:modif}
   Let $a$ be a point in a face $F$ of \cnmoma{\mu}. A function $g$ belongs to $\calG_a^+$
   if and only if it belongs to $\calG_{F,a}^+$ and vanishes $\mu$-a.e.\ on $\{\moma\notin\cl{F}\}$.
   Assuming that a left-hand side is not $+\infty$
   or that no term on a right-hand side is $+\infty$,
   \begin{align}\label{E:prevod2}
        \fuc_{\nm}(g)&=\Crc{F,\nm}+\fuc_{F,\nm}(g)\,,\qquad g\in\calG_{a}^+\,,\\
            \label{E:prevod3}
        \vafu_{\nm}(a)&=\Crc{F,\nm}+\vafu_{F,\nm}(a)\,.
   \end{align}
%   where $\Crc{F,\nm}\triangleq\int_{\{\moma\not\in\clF\}}\,\nm(\cdot,0)\,d\mu$.
\end{lemma}

\begin{Proof}
 The first assertion follows from Lemma~\ref{L:suppg}. Then, for $g$ in $\calG_{a}^+$
 \[
    \inttt{\TT}\:\nm(\ttt,g(\ttt))\:\mu({\mathrm d}\ttt)
        =\inttt{\{\moma\not\in\clF\}}\nm(\ttt,0)\,\mu({\mathrm d}\ttt)
            +\inttt{\{\moma\in\clF\}}\:\nm(\ttt,g(\ttt))\:\mu({\mathrm d}\ttt)\,,
 \]
 if the integral on the left differs from $+\infty$ or if neither integral
 on the right equals $+\infty$. Hence, \eqref{E:prevod2} holds, and
 the quantification there is equivalently over $g\in\calG_{F,a}^+$
 vanishing on $\{\moma\notin\cl{F}\}$. As $\vafu_{F,\nm}(a)$ equals the infimum
 of $\fuc_{F,\nm}(g)$ over such functions $g$, eq.~\eqref{E:prevod3} follows.
\end{Proof}

\begin{lemma}\label{L:modifdom}
   If $F$ is a face of \cnmoma{\mu}, and $\vafu_{\nm}(a)<+\infty$ for some
   $a\in\riF$, then $\ri{\dom{\vafu_{F,\nm}}}$ is equal to \riF.
\end{lemma}

\begin{Proof}
  The assumptions and Theorem~\ref{T:domain} imply that $\riF$ is contained in
  $\dom{\vafu_{\nm}}$. By eq.~\eqref{E:prevod3}, $\riF$ is contained in $\dom{\vafu_{F,\nm}}$.
  Since $\dom{\vafu_{F,\nm}}\pdm\cnmoma{\mu^{\moma^{-1}(\clF)}}=F$, using
  Lemma~\ref{L:facemoma}, the assertion follows.
\end{Proof}

\bigskip

 The set $\Theta_{F,\nm}$ consists of those $\vte\in\dom{\dufu_{F,\nm}}$
 for which the function $r\mapsto\conm(\ttt,r)$ is finite around
 $r=\zat{\vte}{\moma(\ttt)}$ when $\moma(\ttt)\in\clF$, for $\mu$-a.a.\
 $\ttt\in\TT$. If $\vte\in\Theta_{F,\nm}$ let
 \[
            f_{F,\vte}(\ttt)\triangleq
                \begin{cases}\displaystyle
                        \:(\conm)'(\ttt,\zat{\vte}{\moma(\ttt)})\,,
                                \quad&\displaystyle\text{if $\moma(\ttt)\in\clF$
                                            and the derivative exists,}\\
                        \:0\,,&\text{otherwise.}
                \end{cases}
 \]

\begin{remark}\label{R:efa.modif}\rm
  The assumption $\Theta_{F,\nm}\neq\pmn$ plays the role of \DCQ
  in the $F$-dual problems and is implied by  the \DCQ for
  the original problem \eqref{E:duva}. Under this assumption and
  attainment in the $F$-dual problem for $a\in\R^d$, each $F$-dual
  solution $\vte$ belongs to $\Theta_{F,\nm}$ and gives rise to the
  same function $f_{F,\vte}$, arguing as in Remark~\ref{R:efduso}.
  This function is referred to as the \emph{effective $F$-dual
  solution} $g^*_{F,a}$ for $a$.
\end{remark}

 For $a\in\cnmoma{\mu}$ let $F(a)$ denote the unique face
 of $\cnmoma{\mu}$ whose relative interior contains~$a$.

\begin{theorem}\label{T:primal.dual}
    For $a\in\R^d$ such that $\vafu_{\nm}(a)$ is finite\index{$\Crc{F,\nm}$~~shift in modified problems}\label{shift in modified problems}

    (i) the $F(a)$-dual value $\dufu_{F(a),\nm}^*(a)$ is attained
        and $\vafu_{\nm}(a)=\Crc{F(a),\nm}+\dufu_{F(a),\nm}^*(a)$,

    (ii) the primal solution $g_a$ exists if and only if
        $\Theta_{F(a),\nm}\neq\pmn$ and the moment vector
        of the effective $F(a)$-dual solution $g^*_{F(a),a}$ exists
        and equals $a$, in which case $g_a=g^*_{F(a),a}$,

    (iii) the generalized primal solution exists if and only if
        $\Theta_{F(a),\nm}\neq\pmn$, in which case $\gea=g^*_{F(a),a}$.
\end{theorem}

\begin{Proof}
 By finiteness of $\vafu_{\nm}(a)$ and Lemma~\ref{L:modifdom},
 $\ri{\dom{\vafu_{F(a),\nm}}} =\ri{F(a)}$. Then, eq.~\eqref{E:prevod3}
 implies that $\vafu_{\nm}(a)$ equals $\Crc{F(a),\nm}+\vafu_{F(a),\nm}(a)$
 where both quantities are finite. Since $\vafu_{F(a),\nm}(a)$
 is finite, the \PCQ in the $F(a)$-primal problem for $a$ holds.
 By Lemma~\ref{L:pri=du}, $\vafu_{F(a),\nm}(a)=\dufu_{F(a),\nm}^*(a)$ and
 an $F(a)$-dual solution for $a$ exists. These observations imply~\emph{(i)}.

 By Lemma~\ref{L:modif}, $g_a$ exists if and only if
 the $F(a)$-primal solution $g_{F(a),a}$ does, in which case
 they coincide. Knowing that the \PCQ holds in the $F(a)$-primal problem
 for~$a$, the latter existence is equivalent by Lemma~\ref{L:ainri}
 to $\Theta_{F(a),\nm}\neq\pmn$ and $g^*_{F(a),a}\in\calG_{F(a),a}$,
 in which case $g_{F(a),a}=g^*_{F(a),a}$ $\mu$-a.e.\
 on $\{\moma\in\cl{F(a)}\}$. The incidence means
 $g^*_{F(a),a}\in\calG_a$. These observations imply \emph{(ii)}.

 By Remark~\ref{R:efa.modif} and Theorem~\ref{T:geprisol}, $\Theta_{F(a),\nm}\neq\pmn$
 is equivalent to existence of the generalized $F(a)$-primal solution $\hat{g}_{F(a),a}$.
 In this case, $\hat{g}_{F(a),a}=g^*_{F(a),a}$. Lemma~\ref{L:modif} implies that
 $\hat{g}_{F(a),a}$ exists if and only if $\gea$ does, in which case they coincide.
 Hence \emph{(iii)} follows.
\end{Proof}

\begin{corollary}\label{C:uvid}
    Existence of the primal solution $g_a$ implies that the
    generalized primal solution $\gea$ exists and equals $g_a$.
\end{corollary}

 Theorem~\ref{T:primal.dual} makes sense also when the value function
 $\vafu_{\nm}$ equals $-\infty$ at some point, thus the \PCQ holds
 for no $a$, as in Example~\ref{Ex:-infty}. There, $\vafu_{\nm}^*$
 is identically $+\infty$ and the dual values equal $-\infty$, so that
 the dual problems \eqref{E:duva} bear no information on the primal ones.
 However, $\vafu_{\nm}$ can be yet finite at some point $a$ and, due to
 Theorem~\ref{T:primal.dual}, the $F(a)$-dual problem of
 Definition~\ref{D:F-problems} provides complete understanding
 of the primal problem for this $a$.

\begin{definition}\label{D:exn}\index{$\exn{\FF_{\nm}}$~~extension of $\FF_{\nm}$}\label{extension of}\rm
   The \emph{extension} \exn{\FF_{\nm}} of the family $\FF_{\nm}$ is defined
   as union of the families $\FF_{F,\nm}=\{f_{F,\vte}\colon\vte\in\Theta_{F,\nm}\}$
   over the faces $F$ of $\cnmoma{\mu}$.
\end{definition}

 The necessary and sufficient condition for existence of a primal solution can be
 reformulated by means of the extension, without ever mentioning convex duality.

\begin{corollary}\label{C:primal.dual}
    Let $a\in\R^d$ and $\vafu_{\nm}(a)$ be finite.
    The families $\calG_a$ and \exn{\FF_{\nm}} intersect if and only if
    the primal solution for $a$ exists, in which case the intersection
    equals $\{g_a\}$.
\end{corollary}

\begin{Proof}
 By Theorem~\ref{T:primal.dual}\emph{(ii)}, if the primal solution $g_a$
 for $a$ exists then $\Theta_{F(a),\nm}\neq\pmn$, the effective $F(a)$-dual
 solution $g^*_{F(a),a}$ is defined, and $g_a=g^*_{F(a),a}$ belongs to
 $\calG_a\cap\FF_{F(a),\nm}$, contained in $\calG_a\cap\exn{\FF_{\nm}}$.

 In the opposite direction, if $\calG_a\cap\exn{\FF_{\nm}}$ contains a function
 $f_{G,\vte}$, where $G$ is a face of \cnmoma{\mu} and $\vte\in\Theta_{G,\nm}$,
 then $f_{G,\vte}\in\calG_a$ implies $a\in G$, by Lemma~\ref{L:suppg}. Therefore,
 $F(a)\pdm G$. It follows from $\Theta_{G,\nm}\pdm\Theta_{F(a),\nm}$ that
 $\vte\in\Theta_{F(a),\nm}$. By Lemma~\ref{L:suppg}, $f_{G,\vte}$ equals $f_{F(a),\vte}$.
 Hence, $\calG_a$ intersects~$\FF_{F(a),\nm}$.  By Lemma~\ref{L:ainri},
 the $F(a)$-primal solution $g_{F(a),a}$ exists and equals $f_{F(a),\vte}$.
 Therefore, $g_a$ exists by Lemma~\ref{L:modif}, and $f_{G,\vte}=g_a$.
 Thus, $\calG_a\cap\exn{\FF_{\nm}}$ equals $\{g_a\}$.
\end{Proof}

\bigskip

 Corollary~\ref{C:primal.dual} practically amounts to solving the equation
 $\int_\TT\,\moma f_{F,\vte}\,{\mathrm d}\mu=a$ over the faces $F$ of \cnmoma{\mu}
 and $\vte\in\Theta_{F,\nm}$, which is within the framework of the
 last inference principle of Subsection~1.B.

\smallskip
 The correction functional $\Cot_{\nm}$ has been temporarily defined
 in eq.~\eqref{E:Corr} under certain conditions which are now relaxed,
 adapting the former definition to the $F$-problems. Analogously to
 \eqref{E:corre}, for $\vte\in\Theta_{F,\nm}$ and $g\geq0$
 $\calT$-measurable let
 \[
    \Dore_{F,\nm}^\vte(g)\triangleq
    \inttt{\{\moma\in\cl{F}\}}\:\big[\nm'_{\sg{s-f_{F,\vte}(\ttt)}}
                   (\ttt,f_{F,\vte}(\ttt))-\zat{\vte}{\moma(\ttt)}\big]\,
                  \big[g(\ttt)-f_{F,\vte}(\ttt)\big]\:\mu({\mathrm d}\ttt)\,.
 \]
 In turn, for any function $g\in\calG^{+}$ with the moment vector
 $\int_{\TT}\,\moma g\,{\mathrm d}\mu$ denoted by $a$, let
 \begin{eqnarray}\label{E:Corr2}
       \Cot_{\nm}(g)\triangleq\Dore_{F(a),\nm}^\vte(g) \quad
        \text{where $\vte\in\Theta_{F(a),\nm}$ is any $F(a)$-dual solution for $a$}\,,
 \end{eqnarray}
 provided that $\vafu_{\nm}(a)$ is finite and $\Theta_{F(a),\nm}\neq\pmn$.
 Recalling that in the $F(a)$-problem the \PCQ holds for $a$ by the finiteness,
 the correction functional is thereby well defined as it has been in eq.~\eqref{E:Corr}.
 By~\eqref{E:corrediff}, if $\nm$ is differentiable then for $g$ and $\vte$ as above
 \begin{equation}\label{E:corrediffgen}
    \Cot_{\nm}(g)=\inttt{\TT}\,
            |\nm'_+(\cdot,0)-\zat{\vte}{\moma}|_+\,g\,{\mathrm d}\mu
 \end{equation}
 where the integral is indeed over $\TT$ because $g\in\calG_a^+$ vanishes
 on $\{\moma\notin\cl{F(a)}\}$.

\begin{theorem}\label{T:main}\index{generalized Pythagorean identity}\label{generalized Pythagorean identity 3}
   For every $a\in\R^d$ with $\vafu_{\nm}(a)$ finite and
   $\Theta_{F(a),\nm}\neq\pmn$, there exists a ($\mu$-a.e.)
   unique $\calT$-measurable function $\pya$ such that
   \begin{equation}\label{E:Pyth.gen}
      \fuc_{\nm}(g)=\vafu_{\nm}(a)+\Bre{g,\pya}+\Cot_{\nm}(g)\,,
                    \qquad g\in\calG_{a}^+\,.
   \end{equation}
   This function $\pya$ equals the effective solution $g^*_{F(a),a}$
   of the $F(a)$-dual problem.
\end{theorem}

\begin{Proof}
 If $\pya$ satisfying \eqref{E:Pyth.gen} exists then its uniqueness
 follows by considering minimizing sequences that necessarily converge
 to $\pya$ locally in measure, similarly to arguments
 at the end of the proof of Theorem~\ref{T:geprisol}.

 It suffices to prove \eqref{E:Pyth.gen} for $g^*_{F(a),a}$ in the role
 of $\pya$. As in the proof of Theorem~\ref{T:primal.dual}, the first
 two hypotheses imply that the \PCQ holds in the $F(a)$-primal problem
 for $a$. Since $\Theta_{F(a),\nm}\neq\pmn$, the \DCQ holds in the
 $F(a)$-dual problems by Remark~\ref{R:efa.modif}. Hence,
 Lemma~\ref{L:duPyth} implies
   \[
        \fuc_{F(a),\nm}(g)=
          \vafu_{F(a),\nm}(a)+\BreZ{g,g^*_{F(a),a}}{F(a)}+\Cot_{F(a),\nm}(g)\,,
            \qquad g\in\calG_{F(a),a}^+\,.
   \]
 Since the functions $g\in\calG_{a}^+$ and $g^*_{F(a),a}$ vanish on
 $\{\moma\not\in\cl{F(a)}\}$, see Lemma~\ref{L:suppg} and Remark~\ref{R:efa.modif},
 the above Bregman distance equals $\Bre{g,g^*_{F(a),a}}$. By the definition~\eqref{E:Corr2},
 $\Cot_{F(a),\nm}(g)=\Cot_{\nm}(g)$ for $g\in\calG_{a}^+$. The assertion
 follows by Lemma~\ref{L:modif} knowing that $\Crc{F(a),\nm}$ is finite.
\end{Proof}

\bigskip

 Comparing Theorems~\ref{T:primal.dual}\emph{(iii)} and~\ref{T:main},
 the generalized primal solution $\gea$ exists if and only if
 eq.~\eqref{E:Pyth.gen} is available, in which case $\pya=\gea$. If a
 primal solution $g_a$ exists then the hypothesis of Theorem~\ref{T:main}
 holds by Theorem~\ref{T:primal.dual}\emph{(ii)}, eq.~\eqref{E:Pyth.gen}
 takes the form
 \[
    \fuc_{\nm}(g)=\fuc_{\nm}(g_a)+\Bre{g,g_a}+\Cot_{\nm}(g)\,,
                   \qquad g\in\calG_{a}^+\,,
 \]
 and $g_a=\pya$. This implies again Corollary~\ref{C:uvid}.

%888888888888888888888888888888888888888888888888888888888888888888888888888888888888
\section{Bregman projections}\label{S:BREGproj}

 For any integrand $\nm\in B$ and $\calT$-measurable function $h$,
 the mapping\index{$\nmh$~~integrand underlying $\Bre{\cdot,h}$}\label{integrand underlying Bre}
 \[
    (\ttt,t)\mapsto\trnarg{\nm}{t,h(\ttt)}\,,
            \qquad\ttt\in\TT,\,t\in\R\,,
 \]
 is denoted by \nmh. It is a normal integrand, see Lemma~\ref{L:nindelta} and
 \cite[Proposition 14.45(c)]{RW}. It is always \emph{assumed} that $h\geq0$,
 and $h(\ttt)>0$ whenever $\nm'_+(\ttt,0)=-\infty$, $\ttt\in\TT$. Then
 for $t\ge 0$
 \[
    \nmh(\ttt,t)=\nm(\ttt,t)-\nm(\ttt,h(\ttt))-
        \nm'_{\sg{t-h(\ttt)}}(\ttt,h(\ttt))[t-h(\ttt)]\,,
 \]
 and $\nmh\in B$, by Lemma~\ref{L:gat}. Since $\Bre{g,h}=\fuc_{\nmh}(g)$
 for $\calT$-measurable functions $g$ on $\TT$, the Bregman distance $\Bre{g,h}$
 as a function of $g$ is an integral functional of the form~\eqref{E:fuc}.

 In this section, the results on the problem \eqref{E:priva} are specialized
 to the minimization in
 \begin{equation}\label{E:privaB}
    \vafu_{\nmh}(a)={\inf}_{g\in\calG^+_a}\: \Bre{g,h}\,,
                \qquad  a\in\R^d\,.
 \end{equation}
 A (generalized) primal solution of this problem is renamed to a
 \emph{(generalized) Bregman projection} of $h$ to $\calG^+_a$ or to $\calG_a$.

 The dual problem to \eqref{E:privaB} features the function $\dufu_{\nmh}$
 that is equal at $\vte\in\R^d$ to the $\mu$-integral of
 \begin{equation}\label{E:dufuninh}
     \conmh(\cdot,\zat{\vte}{\moma})=\conm(\cdot,\zat{\vte}{\moma}+
        \nm'_{\sg{\zat{\vte}{\moma}}}(\cdot,h))-
        \conm(\cdot,\nm'_{\sg{\zat{\vte}{\moma}}}(\cdot,h))
 \end{equation}
 using Lemma~\ref{L:gat}. In particular, $\dufu_{\nmh}(\vte)=0$ at $\vte=\bnu$,
 by the assumption on $h$. In \eqref{E:dufuninh}, the missing arguments
 $\ttt\in\TT$ of functions are the same, for example the left hand side denotes
 the function $\ttt\mapsto\conmh(\ttt,\zat{\vte}{\moma(\ttt)})$. This convention
 is applied below without any further comments. By  Lemma~\ref{L:gat},
 $\nmh'(\cdot,+\infty)$ equals $\nm'(\cdot,+\infty)-\nm'_{+}(\cdot,h)$.
 Referring to \eqref{E:DCQ}, the crucial set $\Theta_{\nmh}$
 consists of those $\vte\in\dom{\dufu_{\nmh}}$ that satisfy
 \begin{equation}\label{E:DCQBre}
   \zat{\vte}{\moma(\ttt)}< \nm'(\ttt,+\infty)-\nm'_+(\ttt,h(\ttt))
        \qquad \text{for $\mu$-a.a.\ $\ttt\in\TT$.}
 \end{equation}
 Since $\nm\in B$, the difference is positive whence $\vte=\bnu$ always belongs
 to $\Theta_{\nmh}$. For $\vte\in\Theta_{\nmh}$ the functions given by
 \begin{equation}\label{E:Bref}
    f_{\nmh,\vte}=\begin{cases}\displaystyle(\conm)'\big(\cdot\,,\zat{\vte}{\moma}
            +\nm'_{\sg{\zat{\vte}{\moma}}}(\cdot\,,h)\big)\,,
              \quad&\displaystyle\text{when the ineq.\ in \eqref{E:DCQBre} holds,}\\
                                        0\,,&\text{otherwise,}
                           \end{cases}
 \end{equation}
 form the family $\FF_{\nmh}$, see~\eqref{E:fte} and Lemma~\ref{L:gat}.
 The family contains the function $h$, parameterized
 by $\vte=\bnu$, see Lemma~\ref{L:gamma0}.

\begin{remark}\label{R:plusnull}\rm
   Equations~\eqref{E:dufuninh} and \eqref{E:Bref} admit simplifications
   on the set $\{h=0\}$. Namely, if $h=0$, it is possible to write
   $\nm'_+(\cdot,h)$ instead of $\nm'_{\sg{\zat{\vte}{\moma}}}(\cdot,h)$
   even if $\zat{\vte}{\moma}<0$, due to the fact that
   $\conm(\cdot, r)=-\nm(\cdot,0)$ and $(\conm)'(\cdot, r)=0$ for
   all $r\le\nm'_+(\cdot,0)$. In particular, if $\nm$ is differentiable
   and conventionally $\nm'(\cdot, 0)=\nm_+'(\cdot,0)$, the indices
   $\sg{\zat{\vte}{\moma}}$ can be omitted in \eqref{E:dufuninh}
   and \eqref{E:Bref}.
\end{remark}

 Since the integrand $\nmh$ is nonnegative, the \PCQ of the problem \eqref{E:privaB}
 for $a\in\R^d$ reduces to $a\in\ri{\dom{\vafu_{\nmh}}}$. Assuming $\vafu_{\nmh}\not\equiv+\infty$,
 thus existence of $g\in\calG^+$  with $\Bre{g,h}$ finite, the relative interiors
 of $\dom{\vafu_{\nmh}}$ and \cnmoma{\mu} coincide, by Lemma~\ref{L:domvafu}.
 Then, the \PCQ is equivalent to $a\in\ri{\cnmoma{\mu}}$, not depending on $h$.

 Theorems~\ref{T:primal.dual} and \ref{T:main} can be reformulated as follows.
 In these reformulations, in addition to restricting $\mu$, the integrand $\nm$
 is replaced by $\nmh$, as indicated in indices. Accordingly, $(F,\nmh)$-problems,
 $(F,\nmh)$-solutions, etc., come into play.

 Recall the running assumption on $h\geq0$, thus finiteness of $\nm'_+(\ttt,h(\ttt))$,
 $\ttt\in\TT$.

 \begin{theorem}\label{T:primal.dualB}
    For every $a\in\dom{\vafu_{\nmh}}$

    (i) the $(F(a),\nmh)$-dual value is attained and
            $\vafu_{\nmh}(a)=\Crc{F(a),\nmh}+\dufu_{F(a),\nmh}^*(a)$,

    (ii) the Bregman projection $g_{\nmh,a}$ of $h$ to $\calG_a$ exists if and only if
            the moment vector of the effective $(F(a),\nmh)$-dual solution $g^*_{F(a),\nmh,a}$
            exists and equals $a$, in which case $g_{\nmh,a}=g^*_{F(a),\nmh,a}$,

    (iii) the generalized Bregman projection $\geah$ of $h$ to $\calG_a$ exists
            and equals $g^*_{F(a),\nmh,a}$.
\end{theorem}

\begin{theorem}\label{T:mainB}
   For every $a\in\dom{\vafu_{\nmh}}$ there exists a unique
   $\calT$-measurable function $\pyah$ such that
   \begin{equation}\label{E:Pyth.genB}
      \Bre{g,h}=\vafu_{\nmh}(a)+\Breh{g,\pyah}+\Cot_{\nmh}(g)\,,
                    \qquad g\in\calG_{a}^+\,.
   \end{equation}
   This function $\pyah$ equals the effective dual solution
   $g^*_{F(a),\nmh,a}$ of the $(F(a),\nmh)$-dual problem.
\end{theorem}

 As a consequence, the generalized Bregman projection $\geah$ of $h$ to $\calG_a$
 equals $\pyah$. The genuine Bregman projection $g_{\nmh,a}$ exists if and only
 if $\geah\in\calG_a$, in which case they coincide and \eqref{E:Pyth.genB} reduces to
 \begin{equation}\label{E:Pyth.genB2}
      \Bre{g,h}=\Breh{g,g_{\nmh,a}}+\Bre{g_{\nmh,a},h}+\Cot_{\nmh}(g)\,,
                    \qquad g\in\calG_{a}^+\,.
 \end{equation}
 A new feature of eqs.~\eqref{E:Pyth.genB} and \eqref{E:Pyth.genB2}
 is the presence of two kinds of Bregman distances, the original
 one based on \nm, and another one based on \nmh. The following
 lemma presents a comparison.

 \begin{lemma}\label{L:twoBre}
    For any nonnegative $\calT$-measurable functions $g,\tilde g$,
    \[
    \Bre{g,\tilde g}=\Breh{g,\tilde g}+
        \inttt{\{\tilde g\neq h\}}\:
           \big[\nm'_{\sg{g-h}}(\ttt,h(\ttt))-
     \nm'_{\sg{\tilde g-h}}(\ttt,h(\ttt))\big][g(\ttt)-h(\ttt)]\,\mu({\mathrm d}\ttt).
    \]
    The integral is nonnegative and vanishes if $\nm(\ttt,\cdot)$
    is differentiable at $t=h(\ttt)$ for $\mu$-a.a.\ $\ttt\in\TT$
    with $h(\ttt)>0$.
\end{lemma}

\begin{Proof}
 Applying Lemma~\ref{L:compB} to $\ga=\nm(\ttt,\cdot)$, $s=g(\ttt)$, $t=h(\ttt)$
 and $r=\tilde{g}(\ttt)$, the above identity follows by integration and implies
 the remaining assertions.
\end{Proof}

\bigskip

 On account of Lemma~\ref{L:twoBre}, if \Breh{g,g_{\nmh,a}} were replaced by
 \Bre{g,g_{\nmh,a}} in eq.~\eqref{E:Pyth.genB2} then an (explicitly specified)
 nonnegative term had to be subtracted on the right-hand side. This term is not
 necessarily canceled by the correction term $\Cot_{\nmh}(g)$ and it may happen,
 see Example~\ref{Ex:Brepro}, that although the Bregman projection exists,
 the inequality $\Bre{g,h}\geq\Bre{g,g_{\nmh,a}}+\Bre{g_{\nmh,a},h}$
 does not hold for some $g\in\calG_{a}^+$.

 On the other hand, if $\nm(\ttt,\cdot)$, $\ttt\in\TT$, is differentiable at each
 positive number then the two kinds of Bregman distances coincide and
 in~eqs.~\eqref{E:Pyth.genB} and \eqref{E:Pyth.genB2} the nuisance of Bregman
 distance based on \nmh disappears. By Theorem~\ref{T:mainB},
 \eqref{E:corrediffgen} and Lemma~\ref{L:gat}, for $a\in\dom{\vafu_{\nmh}}$
 and $g\in\calG_{a}^+$ (see Remark~\ref{R:plusnull} for $\nm'(\ttt,0)$)
 \begin{equation}\label{E:Pyth.genB3}
      \Bre{g,h}=\vafu_{\nmh}(a)+\Bre{g,\pyah}
        +\inttt{\TT}\:
            |\nm'(\ttt,0)-\nm'(\ttt,h(\ttt))-\zat{\vte}{\moma(\ttt)}|_+\cdot g(\ttt)\,\mu({\mathrm d}\ttt),
 \end{equation}
 where $\vte\in\Theta_{F(a),\nmh}$ is any solution of the $(F(a),\nmh)$-dual problem.
 The integral accounts for the lack of essential smoothness of the functions
 $\nm(\ttt,\cdot)$ at $0$.

\smallskip
 The results below deal with the special situation when the function $h$ projected
 to~$\calG_a$ belongs to the family $\FF_{\nm}$. In this situation, regularity assumptions
 enable to relate directly the Bregman distance minimization \eqref{E:privaB},
 the original primal problem \eqref{E:priva}, and even its dual \eqref{E:duva}.

\begin{lemma}\label{L:esssmooBre}
   Let $\nm$ be essentially smooth and $a\in\R^d$.

   (i) If $\teob\in\Theta_{\nm}$ and $\dufu_{\nm}(\teob)$ is finite
   then the primal problem \eqref{E:priva} for $a$ is equivalent to
   minimization of $\Bre{g,f_\teob}$ subject to $g\in\calG_a^+$.

   (ii) If $g\in\calG^+_a$, $\fuc_{\nm}(g)$ is finite, and the \DCQ
   holds then the dual problem \eqref{E:duva} for $a$ is equivalent
   to minimization of $\Bre{g,f_\teob}$ subject to $\teob\in\Theta_{\nm}$.
\end{lemma}

\begin{Proof}
 By essential smoothness, Lemma~\ref{L:KEYID} implies that for
 $g\in\calG_a^+$ and  $\teob\in\Theta_{\nm}$
 \begin{equation}\label{E:spex}
   \fuc_{\nm}(g)=\zat{\teob}{a}-\dufu_{\nm}(\teob)+\Bre{g,f_\teob}
 \end{equation}
 because the term $\Dore_{\nm}^\teob(g)$ vanishes. Then, \emph{(i)} follows and
 likewise \emph{(ii)}, by Lemmas~\ref{L:Theta} and \ref{L:dualsol}.
\end{Proof}

\bigskip

 Lemma~\ref{L:esssmooBre} is well known, and so are also the following results
 when \nm is essentially smooth. Below, however, only the differentiability of
 this integrand is required, thus the correction term need not vanish. The role
 of $h$ is played by $f_\teob$ with $\teob\in\Theta_\nm^+$. Recall that
 $\Theta_\nm^+$ consists of those $\teob\in\Theta_\nm$ for which
 $\zat{\teob}{\moma}\geq \nm'(\cdot, 0)\;[\mu]$. When $\nm$ is not essentially
 smooth, this is, in general, a proper subset of $\Theta_\nm$, and the results
 below need not hold for all $\teob\in\Theta_\nm$, see Example~\ref{E:threeelem}.

\begin{theorem}\label{T:CLASSIC}
   Suppose $\nm$ is differentiable and $\Theta_{\nm}^+\neq\pmn$.

   (i)
   If $\vafu_\nm\not\equiv+\infty$ then $\vafu_{\nm}$ and $\dufu_{\nm}$
   are proper. If, in addition, $\teob\in\Theta_{\nm}^+$ then
   $\dom{\vafu_{[\nm f_\teob]}}$ equals $\dom{\vafu_{\nm}}$ and
   for $a$ in this domain
   \begin{equation}\label{E:porov}
     \fuc_{\nm}(g)-\vafu_{\nm}(a)=\Bre{g,f_\teob}-\vafu_{[\nm f_\teob]}(a)
            =\Bre{g,\gea}+\Cot_\nm (g)\,,
                \quad\quad g\in\calG^+_a\,,\:\teob\in\Theta_\nm^+\,.
   \end{equation}
   Further, for such $a$ the generalized Bregman projection
   $\hat{g}_{[\nm f_\teob],a}$ exists and equals $\gea$, and the condition
   $\gea\in\calG^+_a$ is necessary and sufficient both for the existence
   of the primal solution $g_a$ and of the Bregman projection $g_{[\nm f_\teob],a}$,
   in which case both are equal to $\gea$.

   (ii)
   If $\vafu_\nm\equiv+\infty$ then $\dom{\vafu_{[\nm f_\teob]}}=\dom{\vafu_{\nm}}=\pmn$
   for $\teob\in\Theta_{\nm}^+$ with $\dufu_{\nm}(\teob)$ finite.
\end{theorem}

\noindent
 For the last assertion, the finiteness hypothesis is essential, see
 Example~\ref{Ex:noPCQ.DCQ}.

\begin{Proof}
 \emph{(i)}~
 By the hypotheses and Theorem~\ref{T:repre}, $\vafu_{\nm}$ and
 $\dufu_{\nm}=\vafu^*_{\nm}$ have nonempty effective domains, hence
 both are proper. Lemma~\ref{L:KEYID} implies eq.~\eqref{E:spex} also
 under the current hypotheses, because if $\nm$ is differentiable and
 $\teob\in\Theta_{\nm}^+$ then $\Dore_{\nm}^\teob(g)=0$ due to
 eq.~\eqref{E:corrediff}. It follows minimizing in eq.~\eqref{E:spex}
 over $g\in\calG_a^+$, or trivially if $\calG^+_a=\pmn$, that
 \[
    \vafu_{\nm}(a)=\zat{\teob}{a}-\dufu_{\nm}(\teob)+\vafu_{[\nm f_\teob]}(a)\,,
        \quad\quad a\in\R^d\,,\;\teob\in\Theta_\nm^+\,,
 \]
 where $\dufu_{\nm}(\teob)$ is finite since $\dufu_{\nm}$ is proper and
 $\teob\in\dom{\dufu_{\nm}}$. This identity implies the claimed equality
 of domains. In case $a\in\dom{\vafu_{\nm}}$ it also implies, by subtraction
 from eq.~\eqref{E:spex}, the first equality in~\eqref{E:porov}. The second
 equality in~\eqref{E:porov} follows from~\eqref{E:Pyth.gen}, as $\tilde g_a$
 in~\eqref{E:Pyth.gen} equals the generalized primal solution $\gea$,
 see also Theorem~\ref{T:primal.dual}(iii).
 The assertions about (generalized) Bregman projections immediately follow
 from the first equality in eq.~\eqref{E:porov}.

 \emph{(ii)}~ The above proof of
 $\dom{\vafu_{\nm}}=\dom{\vafu_{[\nm f_\teob]}}$ goes through also when
 $\vafu_{\nm}\equiv +\infty$ but $\dufu_{\nm}(\teob)$ is finite, since
 the hypothesis $\vafu_{\nm}\not\equiv +\infty$ has been used only to
 guarantee that finiteness.%, when $h=f_{\vte}$, $\teob\in\Theta_{\nm}^+$.
\end{Proof}

\bigskip

 The following lemma relates the dual functions $\dufu_\nm$ and $\dufu_{\nmh}$,
 as well as the families $\FF_\nm$ and $\FF_{\nmh}$, corresponding to the problems
 \eqref{E:priva} and \eqref{E:privaB}.

\begin{lemma}\label{L:classic}
   If $\nm$ is differentiable then for $\teob\in\Theta_{\nm}^+$
   with $\dufu_{\nm}(\teob)$ finite
   \begin{equation}\label{E:duid}
         \dufu_{[\nm f_\teob]}(\vte)=\dufu_{\nm}(\vte+\teob)-\dufu_{\nm}(\teob)\,,
                    \qquad  \vte\in\R^d\,,
   \end{equation}
   and $\vte\in\Theta_{[\nm f_\teob]}$ is equivalent to $\vte+\teob\in\Theta_{\nm}$,
   in which case $f_{[\nm f_\teob],\vte} =f_{\vte+\teob}$. The families $\FF_{\nm}$
   and $\FF_{[\nm f_\teob]}$ coincide, and so do their extensions $\exn{\FF_{\nm}}$
   and \exn{\FF_{[\nm f_\teob]}}.
\end{lemma}

\begin{Proof}
 The value $\dufu_{[\nm f_\teob]}(\vte)$ is obtained by integrating the function
 $[\nm f_\teob]^*(\cdot,\zat{\vte}{\moma})$ which is equal to
 $\nm^*(\cdot, \zat{\vte}{\moma}+\nm'(\cdot,f_{\teob}))-\nm^*(\cdot,\nm'(\cdot,f_{\teob}))$
 by~\eqref{E:dufuninh} and Remark~\ref{R:plusnull}. Here,
 $\nm'(\cdot,f_{\teob})=\zat{\teob}{\moma}$ by the definition~\eqref{E:fte} of $f_{\teob}$,
 differentiability, Lemma~\ref{L:gamma}\emph{(i)} and the assumption $\teob\in\Theta_{\nm}^+$.
 Therefore,  $\dufu_{[\nm f_\teob]}(\vte)$ is obtained by integrating the difference
 $\nm^*(\cdot, \zat{\vte}{\moma}+\zat{\teob}{\moma})-\nm^*(\cdot, \zat{\teob}{\moma})$.
 This proves~\eqref{E:duid}, using that $\dufu_{\nm}(\teob)$ is finite.

 By~\eqref{E:DCQ}, $\vte\in\Theta_{[\nm f_\teob]}$ is equivalent to
 $\vte\!\in\!\dom{\dufu_{[\nm f_\teob]}}$ and
 $\zat{\vte}{\moma}\!<\![\nm f_\teob]'(\cdot,+\infty)\,[\mu]$.
 This takes place if and only if $\vte+\teob\in\dom{\dufu_{\nm}}$ and
 $\zat{\vte}{\moma}<\nm'(\cdot,+\infty)-\nm'(\cdot,f_\teob)\,[\mu]$,
 by~\eqref{E:duid} and Lemma~\ref{L:gat}. Hence, the second assertion
 is proved. By \eqref{E:Bref} and Remark~\ref{R:plusnull},
 $f_{[\nm f_\teob],\vte}$ is equal to
 $(\conm)'(\cdot,\zat{\vte}{\moma}+\nm'(\cdot,f_{\teob}))$
 where $ \nm'(\cdot,f_{\teob})=\zat{\teob}{\moma}$ as above.
 This implies that $f_{[\nm f_\teob],\vte} =f_{\vte+\teob}$
 and $\FF_{\nm}=\FF_{[\nm f_\teob]}$. The same proof
 works for the families built upon faces $F$ of $\cnmoma{\mu}$
 whose unions define the extensions of $\FF_{\nm}$ and
 $\FF_{[\nm f_\teob]}$, using that the hypothesis
 $\teob\in\Theta_{\nm}^+$ implies $\teob\in \Theta_{F,\nm}^+$
 for each face $F$ of $\cnmoma{\mu}$. Thus, also these extensions
 coincide.
\end{Proof}

\begin{corollary}\label{C:dduBre}
   If $\nm$ is differentiable then for $\teob\in\Theta_{\nm}^+$
   with $\dufu_{\nm}(\teob)$ finite
   \[
     \dufu^*_{[\nm f_\teob]}(a)=\dufu^*_{\nm}(a)+\dufu_{\nm}(\teob)-\zat{\teob}{a}\,,
                    \qquad  a\in\R^d\,.
   \]
\end{corollary}

 In above results, the (generalized) Bregman projections of $f_{\teob}$,
 to $\calG_a$, $\teob\in\Theta_{\nm}^+$, are related to the original primal
 and dual problems, \eqref{E:priva} and \eqref{E:duva}, not depending on
 $\teob$. This section is concluded by analogous results for the Bregman
 projection problem~\eqref{E:privaB} with the function $h$ arbitrary,
 subject to the running assumption. They feature the set $\Theta_{\nmh}^+$
 consisting of those $\teob\in\Theta_{\nmh}$ for which
 \[
    \zat{\teob}{\moma}\ge\nmh'(\cdot,0)=\nm'(\cdot,0)-\nm'(\cdot,h)\;\; [\mu]\,,
 \]
 using Lemma~\ref{L:gat}. To simplify the notation in Theorem~\ref{T:classic},
 the function $f_{\nmh,\teob}$, see eq.~\eqref{E:Bref} and Remark~\ref{R:plusnull},
 will be denoted by $h_{\teob}$. Note that $\Theta_{\nmh}^+$ contains the origin
 $\bnu$ and $h=h_{\bnu}$.

\begin{theorem}\label{T:classic}
   Suppose $\nm$ is differentiable.

   For $\teob\in\Theta_{\nmh}^+$ with $\dufu_{\nmh}(\teob)$ finite,
   $\Theta_{[\nm h_\teob]}$ coincides with $\Theta_{\nmh}-\teob$. For $\vte$ in
   that set $f_{[\nm h_\teob],\vte}$ equals $f_{\nmh,\vte+\teob}=h_{\vte+\teob}$.
   Further, $\FF_{\nmh}$ coincides with $\FF_{[\nm h_\teob]}$ and so do their extensions.

   If $\vafu_{\nmh}\not\equiv+\infty$ then
   $\dom{\vafu_{\nmh}}=\dom{\vafu_{[\nm h_\teob]}}$ for $\teob\in\Theta_{\nmh}^+$,
   and for \mbox{$a\in\dom{\vafu_{\nmh}}$}
   \[
     \Bre{g,h}-\vafu_{\nmh}(a)\,{=}\,\Bre{g,h_\teob}-\vafu_{[\nm h_\teob]}(a)
         \,{=}\,\Bre{g,\hat{g}_{\nmh,a}}+\Cot_{\nmh}(g)\,,
        \quad g\in\calG_a^+\,,\;\teob\in\Theta_{\nmh}^+\,,
   \]
   where $\Cot_{\nmh}(g)$ equals the integral in eq.~\eqref{E:Pyth.genB3}.
   Each $h_{\teob}$ above, including $h_{\bnu}=h$, has the same generalized
   Bregman projection to $\calG_a$. This generalized projection belongs to
   $\calG_a$ if and only if all the projections exist, in which case they coincide.
\end{theorem}

\begin{Proof}
 By the assumption, the integrand $\nmh$ is differentiable. By Lemma~\ref{L:compB}
 and the differentiability of $\nm$, $[\nmh h_{\teob}]=[\nm h_{\teob}]$. It suffices
 to apply Lemma~\ref{L:classic} and Theorem~\ref{T:CLASSIC} to the integrand
 $\nmh$ in the role of $\nm$, which gives rise to the same Bregman distance
 as $\nm$, by Lemma~\ref{L:twoBre}.
\end{Proof}

%9999999999999999999999999999999999999999999999999999999999999999999999999999999999999
\section{Generalized solutions of the dual problem}\label{S:dual}

 The results of Section~\ref{S:main} addressed the primal problem
 in absence of the \PCQ. This section approaches the dual problem
 \eqref{E:duva} in this very respect.

 When the value function $\vafu_{\nm}$ is proper, the bi\-conjugate
 $\vafu^{**}_{\nm}$ is equal to the lsc envelope of~$\vafu_{\nm}$,
 and $\vafu^{**}_{\nm}=\dufu_{\nm}^*$ by Theorem~\ref{T:repre}.
 In this case, $\dom{\dufu^*_{\nm}}$ contains $\dom{\vafu_{\nm}}$
 and is contained in its closure which is equal to the closure of
 $\cnmoma{\mu}$, by Lemma~\ref{L:domvafu}. In general, both inclusions
 can be strict. Proposition~\ref{P:dualdom} gives a sufficient
 condition for equality in the first one. First,
 dual attainment is briefly addressed.

\begin{lemma}\label{L:suppHyp}
   Let $H=\{x\colon \zat{\teob}{x}=0\}$ be a hyperplane such that
   $H_<=\{x\colon \zat{\teob}{x}<0\}$ contains $\ri{\cnmoma{\mu}}$.
   Then, for $a\in H$ and $\vte\in\R^d$ with $\dufu_{\nm}(\vte)$ finite
    \[
        \dufu_\nm^*(a)\ge\zat{\vte}{a}-\dufu_{\nm}(\vte)+\inttt{\{\moma\in
               H_<\}}\:\left[\conm\big(\ttt,\zat{\vte}{\moma(\ttt)}\big)-
                \conm(\ttt,-\infty)\right]\,\mu({\mathrm d}\ttt)\,.
    \]
\end{lemma}

\begin{Proof}
 By Corollary~\ref{C:zerooutside}, $\moma\in\cl{\cnmoma{\mu}}$
 $\mu$-a.e. The closure is contained in $H\cup H_<$ by assumption.
 Hence, using that $\zat{\teob}{\moma}$ vanishes on the set $\{\moma\in H\}$
 and  $\dufu_{\nm}(\vte)$ is finite,
 \[
    \dufu_{\nm}(\vte+t\teob)=
        \inttt{\{\moma\in H\}}\:\conm\big(\ttt,\zat{\vte}{\moma(\ttt)}\big)\,\mu({\mathrm d}\ttt)
        \;+
        \inttt{\{\moma\in H_<\}}\:\conm\big(\ttt,\zat{\vte+t\teob}{\moma(\ttt)}\big)\,\mu({\mathrm d}\ttt)\,.
 \]
 Here, the first integral is finite and equals
 $\dufu_{\nm}(\vte)-\int_{\{\moma\in H_<\}}\:
       \conm\big(\ttt,\zat{\vte}{\moma(\ttt)}\big)\,\mu({\mathrm d}\ttt)$.
 When $t\to+\infty$, the second one converges to the integral over
 $\{\moma\in H_<\}$ of $\conm(\ttt,-\infty)$, by monotone convergence.
 Since $\dufu_\nm^*(a)$ is lower bounded by
 $\lim_{t\to+\infty}\left[\zat{\vte+t\teob}{a}-\dufu_{\nm}(\vte+t\teob)\right]$
 and $\zat{\teob}{a}=0$, the assertion follows.
\end{Proof}

\begin{proposition}\label{P:dualatt}
   If $\vafu_\nm\not\equiv +\infty$ and a dual solution $\vte$ for some
   $a\in\R^d$  exists, then either the \PCQ holds for $a$ or else $H$ and
   $H_<$ exist as in Lemma~\ref{L:suppHyp} such that $a\in H$ and
   $\zat{\vte}{\moma}\le\nm'_+(\cdot,0)$ $\mu$-a.e.\ on $\{\moma\in H_<\}$.
\end{proposition}

\begin{Proof}
 The hypotheses imply that $\vafu_\nm$ is proper and $a\in\cl{\cnmoma{\mu}}$.
 If the \PCQ for~$a$ fails then $a\not\in\ri{\cnmoma{\mu}}$, see Lemma~\ref{L:domvafu}.
 Therefore, there exists a hyperplane $H=\{x\colon \zat{\teob}{x}=0\}$
 containing $a$ and the origin such that \ri{\cnmoma{\mu}} is contained
 in $H_<$ as in Lemma~\ref{L:suppHyp} \cite[Theorem 11.2]{Rock}.
 Since $\vte$ is a dual solution for $a$, the nonnegative difference
 in the integral of Lemma~\ref{L:suppHyp} equals zero $\mu$-a.e. on
 the set $\{\moma\in H_<\}$. As that difference vanishes if and only
 if $\zat{\vte}{\moma(\ttt)}\le\nm'_+(\ttt,0)$, this completes the proof.
\end{Proof}

\begin{corollary}\label{C:dualatt}
   If $\vafu_\nm\not\equiv +\infty$ and $\nm'_+(\cdot,0)=-\infty\,[\mu]$,
   in particular if $\nm$ is essentially smooth, then the \PCQ for~$a$
   is necessary and sufficient for existence of a dual solution for~$a$.
\end{corollary}

\begin{Proof}
 By Lemma~\ref{L:pri=du}, sufficiency holds for any $\nm\in B$. Necessity
 under the additional hypothesis follows from Proposition~\ref{P:dualatt}, as
 the hypothesis on \nm rules out the second contingency there.
\end{Proof}

\begin{proposition}\label{P:dualdom}
    If $\vafu_{\nm}$ is proper and its effective domain is equal to
    $\ri{\cnmoma{\mu}}$ then \dom{\dufu^*_{\nm}} coincides
    with \dom{\vafu_{\nm}}.
\end{proposition}

\begin{Proof}
 By the facts sent forward at the beginning of this section, it suffices
 to show that if $a$ belongs to the closure but not to the relative interior
 of \cnmoma{\mu} then $\dufu_{\nm}^*(a)$ equals $+\infty$. Let $H$,
 $H_<$ and $a\in H$ be as in Lemma~\ref{L:suppHyp}. Since $\vafu_{\nm}$
 is proper there exists $\vte$ with $\dufu_{\nm}(\vte)$ finite. The hypothesis
 $\dom{\vafu_\nm}=\ri{\cnmoma{\mu}}$ implies by Corollary~\ref{C:dom.open}
 that the integral of $ \nm(\cdot,0)$ over $H_<$ equals $+\infty$.
 Since $\conm(\ttt,-\infty)=-\nm(\ttt,0)$, Lemma~\ref{L:suppHyp}
 implies $\dufu_\nm^*(a)=+\infty$.
\end{Proof}

\bigskip

 Analogously to the generalized primal solutions, a generalized dual
 solution is introduced for each $a\in\R^d$ with finite dual value
 $\dufu_{\nm}^*(a)$, attained or not. More precisely, this concept
 generalizes that of the effective dual solution rather than that of a
 dual solution proper. It requires the \DCQ which is assumed throughout
 the remaining part of this section. By Lemma~\ref{L:Theta}, there exist
 sequences $\vte_n\in\Theta_{\nm}$ with $\zat{\vte_n}{a}-\dufu_{\nm}(\vte_n)$
 tending to $\dufu_{\nm}^{*}(a)$. A \calT-measurable function
 $h_a$\index{$h_a$~~generalized dual solution for $a$}\label{generalized dual solution for}
 is a \emph{generalized dual solution} for $a$ if for each sequence
 $\vte_n$ as above the functions $f_{\vte_n}$ converge to $h_a$
 locally in measure. Existence of generalized dual solutions follows
 from the main result of this section.

\begin{theorem}\label{T:GMLE}
   Assuming the \DCQ, for every $a\in\R^d$ with $\dufu_{\nm}^{*}(a)$ finite
   there exists a~unique \calT-measurable function $h_a$ such that
   \begin{equation}\label{E:GMLE}
       \dufu_{\nm}^{*}(a)-\big[\zat{\vte}{a}-\dufu_{\nm}(\vte)\big]
            \ge\Bre{h_a,f_\vte}\,,  \qquad \vte\in\Theta_{\nm}\,.
   \end{equation}
\end{theorem}

 If the effective dual solution $g_a^*=f_\vte$ exists, where $\vte\in\Theta_{\nm}$
 is a dual solution for $a$, then $g_a^*=h_a$ by ineq.~\eqref{E:GMLE}.
 If $\vte_n\in\Theta_{\nm}$ is a maximizing sequence for
 $\zat{\vte}{a}-\dufu_{\nm}(\vte)$, ineq.~\eqref{E:GMLE} implies that
 the Bregman distances $\Bre{h_a,f_{\vte_n}}$ tend to zero, and then
 $f_{\vte_n}\rightsquigarrow h_a$ by Corollary~\ref{C:inmeas}.
 Thus, $h_a$ is the generalized dual solution for $a$. This establishes
 also the uniqueness in Theorem~\ref{T:GMLE}.

 First, a special case of Theorem~\ref{T:GMLE} is established, for its
 simplicity and independent interest. In this case, the generalized primal
 and dual solutions for~$a$ coincide.

\begin{proposition}\label{P:GMLE}
   Assuming the \DCQ, for $a\in\R^d$ with $\dufu_{\nm}^{*}(a)$ finite
   and zero duality gap, ineq.~\eqref{E:GMLE} holds with $h_a=\tilde g_a$,
   see Theorem~\ref{T:main}.
\end{proposition}

\begin{Proof}
 By assumptions, $\vafu_\nm(a)$ is finite. Let $g_n$ be a sequence in $\calG_a^+$
 with $\fuc_\nm(g_n)$ converging to $\vafu_\nm(a)$. Limiting along $g_n$
 in Theorem~\ref{T:main}, $\Bre{g_n,\pya}\to0$ whence $g_n\rightsquigarrow\pya$,
 by Corollary~\ref{C:inmeas}. Lemma~\ref{L:KEYID} implies
 \[
        \fuc_\nm(g)\geq \zat{\vte}{a}-\dufu_{\nm}(\vte)+\Bre{g,f_\vte}\,,
            \quad g\in\calG_a^+\,,\;\vte\in\Theta_{\nm}\,.
 \]
 Limiting here along $g_n$, $\vafu_\nm(a)\geq \zat{\vte}{a}-\dufu_{\nm}(\vte)+\Bre{\pya,f_\vte}$
 for $\vte\in\Theta_{\nm}$, by Lemma~\ref{L:Brelsc}. This and the hypothesis
 $\vafu_\nm(a)=\dufu_\nm^*(a)$ imply that $\Bre{\pya,f_{\vte_n}}\to0$ for each
 sequence $\vte_n$ in $\Theta_{\nm}$ with $\zat{\vte_n}{a}-\dufu_{\nm}(\vte_n)$
 converging to $\dufu_\nm^*(a)$. The assertion $f_{\vte_n}\rightsquigarrow\pya$
 follows by Corollary~\ref{C:inmeas}.
\end{Proof}

\bigskip

 Example~\ref{Ex:gap} illustrates a situation when the \PCQ fails,
 the primal and dual values are finite but different and the primal
 solution $g_a$ is different from $h_a$. Additionally, $\mu$ is finite
 and $\moma$ bounded.

 To prove Theorem ~\ref{T:GMLE} in general, the following lemmas and
 corollary are needed. The inequality below compares the Jensen
 difference of $\dufu_{\nm}$ with Bregman distances.

\begin{lemma}\label{L:JBr}
   If $\dufu_\nm$ is proper then for $\teob_1,\teob_2$
   in $\Theta_{\nm}$ and $0<t<1$
   \[\begin{split}
       t \dufu_{\nm}(\teob_1) +&(1-t) \dufu_{\nm}(\teob_2)
                - \dufu_{\nm}(t\teob_1+(1-t)\teob_2)\\
           &\geq t \Bre{f_{t\teob_1+(1-t)\teob_2},f_{\teob_1}}
             +(1-t) \Bre{f_{t\teob_1+(1-t)\teob_2},f_{\teob_2}}\,.
             \end{split}
   \]
\end{lemma}

\begin{Proof}
 The left-hand side is equal to
   \[
       t [\dufu_{\nm}(\teob_1) - \dufu_{\nm}(t\teob_1+(1-t)\teob_2)]
       +(1-t) [\dufu_{\nm}(\teob_2) - \dufu_{\nm}(t\teob_1+(1-t)\teob_2)]
   \]
 where all values are finite. The left bracket takes the form
 \[\begin{split}
&\inttt{\TT}\:\Big[\conm\big(\ttt,\zat{\teob_1}{\moma(\ttt)}\big)\:\mu({\mathrm
d}\ttt)
              -\conm\big(\ttt,\zat{t\teob_1+(1\!-\!t)\teob_2}{\moma(\ttt)}\big)\Big]\:\mu({\mathrm d}\ttt)\\
             &=\inttt{\TT}\:\Big[\trnarg{\conm(\ttt,\cdot)}{\zat{\teob_1}{\moma(\ttt)},\zat{t\teob_1+(1\!-\!t)\teob_2}{\moma(\ttt)}}\\
             &\qquad+(1\!-\!t)\zat{\teob_1-\teob_2}{\moma(\ttt)}f_{t\teob_1+(1\!-\!t)\teob_2}(\ttt)\Big]\:\mu({\mathrm d}\ttt)
    \end{split}
 \]
 and the right bracket
 \[
 \inttt{\TT}\:\Big[\trnarg{\conm(\ttt,\cdot)}{\zat{\teob_2}{\moma(\ttt)},\zat{t\teob_1+(1\!-\!t)\teob_2}{\moma(\ttt)}}
             +t\zat{\teob_2-\teob_1}{\moma(\ttt)}f_{t\teob_1+(1\!-\!t)\teob_2}(\ttt)\Big]\:\mu({\mathrm d}\ttt)\,.
 \]
 Then, the left-hand side of the inequality rewrites to
 \[\begin{split}
 \inttt{\TT}\:\Big[&t\trnarg{\conm(\ttt,\cdot)}{\zat{\teob_1}{\moma(\ttt)},\zat{t\teob_1+(1\!-\!t)\teob_2}{\moma(\ttt)}}\\
                    &+ (1-t)\trnarg{\conm(\ttt,\cdot)}{\zat{\teob_2}{\moma(\ttt)},\zat{t\teob_1+(1\!-\!t)\teob_2}{\moma(\ttt)}}
                        \Big]\:\mu({\mathrm d}\ttt)\,.
 \end{split}\]
 The assertion follows by the consequence
 $\trndu{r_2,r_1}\geq\trnga{\cgd(r_1),\cgd(r_2)}$ of Lemma~\ref{L:Bld},
 where $u(r)={\coga}'(r)$.
\end{Proof}

\begin{corollary}\label{C:JBr}
 If $\dufu^*_{\nm}(a)$ is finite, $\teob_1,\teob_2\in\Theta_{\nm}$ and $0<t<1$ then
   \[\begin{split}
       t \big[\dufu^*_{\nm}(a)-[\zat{\teob_1}{a}-\dufu_{\nm}(\teob_1)]\big]
       +&(1-t) \big[\dufu^*_{\nm}(a)-[\zat{\teob_2}{a}-\dufu_{\nm}(\teob_2)]\big]
                \\
           \geq t& \Bre{f_{t\teob_1+(1-t)\teob_2},f_{\teob_1}}
             +(1-t) \Bre{f_{t\teob_1+(1-t)\teob_2},f_{\teob_2}}\,.
             \end{split}
   \]
\end{corollary}

\begin{Proof}
 The Jensen difference in  Lemma~\ref{L:JBr} is equal to
 \[\begin{split}
     \big[\zat{t\teob_1+(1-t)\teob_2}{a}-&\dufu_{\nm}(t\teob_1+(1-t)\teob_2)\big]\\
        &
       -t \big[\zat{\teob_1}{a}-\dufu_{\nm}(\teob_1)\big]
       -(1-t) \big[\zat{\teob_2}{a}-\dufu_{\nm}(\teob_2)\big]
 \end{split}\]
 where the first bracket is dominated by $\dufu^*(a)$.
\end{Proof}

\begin{lemma}\label{L:inmeas2}
    Let $C\in\calT$ have finite $\mu$-measure, let $L$, $\xi$, and $\delta$
    be positive numbers. Then there exists $K>L$ such that if $\Bre{h,g}\leq\delta$
    for some  nonnegative \calT-measurable functions $g,h$ then
         \[
            \mu(C\cap\{g>K\})<\xi+\mu(C\cap\{h>L\})\,.
         \]
\end{lemma}

\begin{Proof}
 Let $M=2\delta/\xi$. By monotonicity, for $K>L$
 \[\begin{split}
    \Bre{g,h}\geq\inttt{\{g> K,h\leq L\}}\:\trnarg{\nm}{\ttt,g(\ttt),h(\ttt)} \;\mu({\mathrm d}\ttt)
             \geq\inttt{\{g> K,h\leq L\}}\:\trnarg{\nm}{\ttt,K,L} \;\mu({\mathrm d}\ttt)\\
             \geq M\cdot\mu\big(C\cap\{\trnarg{\nm}{\cdot,K,L}\geq M\}\cap\{g> K,h\leq L\}\big)
 \end{split}\]
 whence
 \[
    \mu(C\cap\{g>K\})\leq
        \tfrac{1}{M}\Bre{g,h}+\mu(C\cap\{\trnarg{\nm}{\cdot,K,L}< M\})+\mu(C\cap\{h>L\})\,.
 \]
 Since $\trnarg{\nm}{\ttt,K,L}\uparrow+\infty$ if $K\uparrow+\infty$
 due to strict convexity, there exists $K>L$ such that
 $\mu\big(C\cap\{\trnarg{\nm}{\cdot,K,L}< M\}\big)<\tfrac12\xi$.
 With this $K$ the assertion follows by the choice of~$M$,
 implying $\tfrac{1}{M}\Bre{h,g}\leq\tfrac12\xi$ whenever
 $\Bre{h,g}\leq\delta$.
\end{Proof}

\bigskip

\noindent
P\,r\,o\,o\,f\ o\,f\ T\,h\,e\,o\,r\,e\,m~\ref{T:GMLE}.
 By assumptions and Lemma~\ref{L:Theta}, $\dufu_{\nm}$ is proper
 and there exists a sequence $\te_n$ in $\Theta_{\nm}$ such that
 $\zat{\te_n}{a}- \dufu_{\nm}(\te_n)$ converges to $\dufu^*(a)$.

 Let $\vte\in\Theta_{\nm}$. Applying Corollary~\ref{C:JBr}
 to $\teob_1=\vte$, $\teob_2=\te_n$ and $0<t_n<1$ yields
 \begin{equation}\label{E:*}\begin{split}
    \big[\dufu^*_{\nm}(a)-[\zat{\vte}{a}-\dufu_{\nm}(\vte)]\big]
                    +&\tfrac{1-t_n}{t_n}
       \big[\dufu_{\nm}^*(a)-[\zat{\te_n}{a}-\dufu_{\nm}(\te_n)]\big] \\
           \geq & \Bre{f_{t_n\vte+(1-t_n)\te_n},f_{\vte}}
             +\tfrac{1-t_n}{t_n}  \Bre{f_{t_n\vte+(1-t_n)\te_n},f_{\te_n}}\,.
             \end{split}
 \end{equation}
 Let $t_n\to 0$ sufficiently slowly to make the second term
 on the left hand side go to zero. Then, \eqref{E:*} implies that
 the sequence $\Bre{f_{t_n\vte+(1-t_n)\te_n},f_{\vte}}$ is bounded,
 the sequence $\Bre{f_{t_n\vte+(1-t_n)\te_n},f_{\te_n}}$
 tends to zero and
 \begin{equation}\label{E:*x}
    \dufu^*_{\nm}(a)-[\zat{\vte}{a}-\dufu_{\nm}(\vte)]
          \geq {\liminf}_{n\to\infty}\:\Bre{f_{t_n\vte+(1-t_n)\te_n},f_{\vte}}\,,
          \quad \vte\in\Theta_{\nm}\,.
 \end{equation}
 By Lemma~\ref{L:Brelsc}, it suffices to prove that the sequence
 $f_{t_n\vte+(1-t_n)\te_n}$ converges locally in measure.

 Let $C\in\mathcal Z$ have finite $\mu$-measure and $\xi>0$.
 Then, $\mu(C\cap\{f_{\vte}>L\})<\xi$ for some $L>0$. Since
 $\Bre{f_{t_n\vte+(1-t_n)\te_n},f_{\vte}}$ is bounded,
 by Lemma~\ref{L:inmeas2} there exists $K>L$ such that
 \begin{equation}\label{E:*iv}
    \mu(C\cap\{f_{t_n\vte+(1-t_n)\te_n}>K\})<\xi+\mu(C\cap\{f_{\vte}>L\})<2\xi
 \end{equation}
 for all $n$. Since $\Bre{f_{t_n\vte+(1-t_n)\te_n},f_{\te_n}}\to 0$,
 Lemma~\ref{L:inmeas} and \eqref{E:*iv} imply that
 for any $\vare>0$
 \begin{equation}\label{E:***}
    \mu(C\cap\{|f_{t_n\vte+(1-t_n)\te_n}-f_{\te_n}|>\vare\})<3\xi\,,
        \qquad\text{eventually in $n$.}
 \end{equation}
 Combining \eqref{E:*iv} and \eqref{E:***},
 \begin{equation}\label{E:unibound}
    \mu(C\cap\{f_{\te_n}>K+\vare\})< 5\xi\,,
        \qquad\text{eventually in $n$.}
 \end{equation}

 By Corollary~\ref{C:JBr} applied to $\teob_1=\te_m$, $\teob_2=\te_n$
 and $t=\frac12$,
 \[
    \Bre{f_{(\te_n+\te_m)/2},f_{\te_m}} +
        \Bre{f_{(\te_n+\te_m)/2},f_{\te_n}}\to0\qquad\text{as $n,m\to\infty$.}
 \]
 This convergence, Lemma~\ref{L:inmeas} and \eqref{E:unibound}
 imply that for any $\vare>0$ and $\xi>0$
 \begin{equation}\label{E:*v}
    \mu(C\cap\{|f_{\te_n}-f_{\te_m}|>\vare\})<\xi+\mu(C\cap\{f_{\te_n}>K+\vare\})<6\xi
 \end{equation}
 provided $n,m$ are sufficiently large. Thus, the sequence $f_{\te_n}$ is
 Cauchy, locally in $\mu$-measure. Hence, $f_{\te_n}\rightsquigarrow h$
 for some $\calT$-measurable nonnegative function $h$. This and \eqref{E:***}
 imply that also $f_{t_n\vte+(1-t_n)\te_n}\rightsquigarrow h$,
 needed to complete the proof. \hfill$\square$

\begin{remark}\label{R:ga=ha}\rm
 In Theorem~\ref{T:GMLE}, it can happen that $\vafu_{\nm}(a)$ is not finite,
 even $\vafu_{\nm}\equiv +\infty$ is allowed. Assuming the \DCQ, if
 $\vafu_{\nm}(a)$ is finite then Lemma~\ref{L:KEYID} implies for $g_n\in\calG_a^+$
 with $\fuc_{\nm}(g_n)\to\vafu_{\nm}(a)$ and $\vte_n\in\Theta_{\nm}$ with
 $\zat{\vte_n}{a}-\dufu_{\nm}(\vte_n)\to\dufu^*_{\nm}(a)$ that
 \[
    \vafu_{\nm}(a)-\dufu_{\nm}^*(a)\ge\limsup\nolimits_{n\to \infty}\;\Bre{g_n,f_{\vte_n}}\,.
        %\ge\inf_{g\in\calG_a, \vte\in\Theta_{\nm}}\;\Bre{g,f_\vte},
 \]
 Here, if $\nm$ is essentially smooth then the equality takes place and
 the limit exists. Since $g_n\rightsquigarrow \tilde g_a=\gea$ by Theorem~\ref{T:main},
 and $f_{\vte_n}\rightsquigarrow h_a$ by Theorem~\ref{T:GMLE}, it follows by
 Lemma~\ref{L:Brelsc} that the duality gap $\vafu_{\nm}(a)-\dufu_{\nm}^*(a)$
 majorizes the Bregman distance $\Bre{\gea,h_a}$ of the generalized solutions.
 Conditions making this bound tight remain elusive.
\end{remark}

%10 10 10 10 10 10 10 10 10 10 10 10 10 10 10 10 10 10 10 10 10 10 10 10 10 10
\section{Examples}\label{S:examples}

 This section demonstrates that `irregular' behavior
 may occur in the primal and dual problems, even in the autonomous case
 $\nm(\ttt, t)=\ga(t)$, $\ttt\in\TT$, $t\in\R$, with $\ga$ differentiable.
 In this situation, $\ga\in\varGamma$ replaces $\nm$ in notations like
 $\fuc_\nm(g)$, $\vafu_{\nm}$, etc. Fig.~\ref{F:1} summarizes some
 properties of the first eight examples. It contains a region marked
 by $\pmn$, see Lemma~\ref{L:ainri}.  Example~\ref{Ex:moma=0} describes
 the situation when the moment mapping vanishes identically, showing how
 the problem of unconstrained minimization of convex integral functionals
 fits into our framework. The remaining three examples illustrate Bregman
 projections and closure.

\begin{figure}[h!]
%\begin{center}
\hspace*{-0.8cm}%
\setlength{\unitlength}{1mm}
\scalebox{2}{\begin{picture}(70,25)
  \put(5,15){\line(1,0){60}}
  \put(33,8){\circle{14}}
  \put(40,15){\circle{14}}
  \put(42,7){\circle{12}}
\end{picture}}\hspace*{-140mm}\setlength{\unitlength}{2mm}%
\begin{picture}(70,25)
  \put(13,13){\makebox(0,0){\footnotesize $\vafu_{\nm}\not\equiv+\infty$}}
  \put(13,17){\makebox(0,0){\footnotesize $\vafu_{\nm}\equiv+\infty$}}
  \put(25,4){\makebox(0,0){\PCQ}}
  \put(48,21){\makebox(0,0){\DCQ}}
  \put(58,3){\makebox(0,0){\footnotesize the primal solution exists}}
  \put(40,18){\makebox(0,0){\bf\tiny\ref{Ex:noPCQ.DCQ}}}
  \put(30,18){\makebox(0,0){\bf\tiny\ref{Ex:vafu.dual}}}
  \put(31,8){\makebox(0,0){\bf\tiny\ref{Ex:emqua}}}
  \put(35.5,12.2){\makebox(0,0){\bf\tiny\ref{Ex:Burg}}}
  \put(34.8,13.6){\makebox(0,0){\bf\tiny\ref{Ex:momefa}}}
  \put(44,5){\makebox(0,0){\bf\tiny\ref{Ex:-infty}}}
  \put(42,12){\makebox(0,0){\bf\tiny\ref{Ex:zero.density}}}
  \put(42,10.5){\makebox(0,0){\bf\tiny\ref{Ex:gap}}}
  \put(37.5,6){\makebox(0,0){\bf\small$\pmn$}}
  \put(38.2,9.7){\makebox(0,0){$\heartsuit$}}
\end{picture}
%\end{center}
 \caption{}\label{F:1}
\end{figure}

\begin{example}\label{Ex:vafu.dual}\rm
 Let $\mu$ be the counting measure and $\moma$ the identity mapping on the set
 $\TT$ of integers. The functional $\fuc_\ga$ is considered with $\ga\in\varGamma$
 given by $\ga(t)=(2t)^{-1}$ for $t>0$. Then, $\fuc_\ga(g)<+\infty$ implies
 that $g$ is positive and $\sum_{\ttt\in\TT}\:1/g(\ttt)$ converges.
 In this case, $g(z)\geq1$ if $|\ttt|$ is sufficiently large whence both
 the positive and negative parts of $\int_\TT\moma g\,{\mathrm d}\mu$ are infinite.
 Therefore, $\vafu_{\ga}\equiv+\infty$ and $\vafu_{\ga}^*\equiv-\infty$.
 The conjugate of $\ga$ is given by $\coga(r)=-\sqrt{-2r}$
 for $r\leq0$ and $\coga(r)=+\infty$ otherwise. Then,
 $\dufu_{\ga}(\vte)=\int_\TT\coga(\vte\moma)\,{\mathrm d}\mu$ equals $0$
 for $\vte=0$ and $+\infty$, otherwise, both the positive and
 negative parts being infinite. In particular, $\dufu_{\ga}$ is proper,
 but its effective domain has empty interior. The \DCQ fails,
 $\Theta_{\ga}=\pmn$. For $a=0$ the dual value is $0$,
 different from the primal one $+\infty$.
\end{example}

\begin{example}\label{Ex:noPCQ.DCQ}\rm
 Let $\mu$ be the counting measure and $\moma$ the identity mapping on
 the set $\TT=\{1,2,\ldots\}$. Let $\ga(t)=t\ln t -t+1$, $t\ge 0$, so that
 $\ga^*(r)=e^r-1$, $r\in\R$. Then, $\dufu_{\ga}(\vte)=\sum_{\ttt\in\TT}\: \coga(\vte\ttt)
 =\vte\cdot(+\infty)$, $\vte\in\R$, and $\dom{\dufu_{\ga}}$ coincides with
 $\Theta_{\ga}=(-\infty,0]$. Although $\dom{\dufu_{\ga}}$ has nonempty interior and
 the \DCQ holds, $\vafu_{\ga}\equiv +\infty$ for otherwise $\dufu_{\ga}$ could not take
 the value $-\infty$. For $\vte<0$ the function $f_{\vte}$ belongs to $\calG$,
 i.\,e., the moment $\int_\TT\moma f_\vte\,{\mathrm d}\mu=\sum_{\ttt=1}^{\infty}\:\ttt e^{\vte\ttt}$
 exists, while $\dufu_{\ga}(\vte)=-\infty$. For $a$ equal to that moment,
 $\inf_{g\in\calG^+_a}\Brega{g,f_{\vte}}=0$, hence $\dom{\vafu_{[\ga f_\vte]}}$
 contains $a$, see~\eqref{E:privaB}. This shows that in the last assertion of
Theorem~\ref{T:CLASSIC} the finiteness assumption is essential.
\end{example}

\begin{example}\label{Ex:momefa}\rm
 Let $\mu$ be the Borel measure on $\TT=\R$ given by ${\mathrm d}\mu=\frac{{\mathrm d}\ttt}{1+\ttt^2}$,
 $\moma$ the identity mapping on $\TT$ and $\ga(t)=t\ln t$, $t\geq0$. Since
 $\mu(\TT)=\pi$ and $\ga\geq\ga(\frac{1}{e})=-\frac{1}{e}$, the value function
 $\vafu_\ga$ is lower bounded by $-\frac{\pi}{e}$. The functional $\fuc_\ga$
 is finite for functions $g\geq0$ on $\TT$ that are nonzero and bounded on bounded sets.
 Then the \PCQ holds for every $a\in\R$. In dual problems, $\coga(r)=e^{r-1}$,
 $r\in\R$, and $\dufu_{\ga}(\vte)=\int_\R\,e^{\vte\ttt-1}\,\mu({\mathrm d}\ttt)$
 is equal to~$\frac{\pi}{e}$ for $\vte=0$ and $+\infty$, otherwise. Therefore,
 the \DCQ holds, $\Theta_{\ga}=\{0\}$ and $\vafu_\ga=\dufu_{\ga}^*\equiv-\frac{\pi}{e}$.
 For each $a\in\R$ the effective dual solution $\efa$ is identically equal
 to $\frac{1}{e}$. The moment of \efa does not exist whence the primal problem
 has no solution. Nevertheless, the generalized primal solution $\gea$ exists
 and equals $\efa$. A modification of this example with information theoretical
 interpretation appeared in~\cite[Example~1]{Csi.Ma.minentrev}.
\end{example}

\begin{example}\label{Ex:zero.density}\rm
 On $\TT=\{0,1\}$, let $\mu$ be the counting measure and $\moma$
 the identity mapping. Let $\ga(t)$ equal $\frac12t^2$ for $t\geq0$.
 Then, $\vafu_{\ga}=\ga$ has the effective domain $[0,+\infty)$. Since
 $\coga(r)$ equals $\frac12r^2$ for $r\geq0$ and $0$ otherwise, the
 \DCQ holds and $\Theta_{\ga}=\R$. For $a=0$ not enjoying the \PCQ,
 each $\vte\leq0$ is a dual solution and the effective dual solution
 $\efa$ is identically equal to $0$. It belongs to $\calG_a$, and
 thus coincides with the primal solution $g_a$.
\end{example}

\begin{example}\label{Ex:emqua}\rm
 On $\TT=\{-1,1\}$, let $\mu$ be the counting measure and $\moma$
 the identity mapping. Let $\ga$ be the same as in Example~\ref{Ex:vafu.dual}.
 In the primal problem
 \[
    \vafu_{\ga}(a)=\inf\Big\{\tfrac{1}{2g(1)}+\tfrac{1}{2g(-1)}\colon
                        g(1),g(-1)>0\,,\,g(1)-g(-1)=a\Big\}=0\,,\quad a\in\R\,,
 \]
 the infimum is not attained. Though the \PCQ holds for each $a\in\R$,
 no primal solution exists. No generalized primal solution exists either,
 since any minimizing sequence $g_n$ in the primal problem for $a$ satisfies
 $g_n(1)\to+\infty$ and $g_n(1)-g_n(-1)=a$, and thus $g_n$ cannot converge
 in any standard sense. Since $\dom{\vafu_\ga}=\R$, the effective domain
 of $\dufu_\ga=\vafu_{\ga}^*$ equals the singleton $\{0\}$. For $\vte=0$
 the function $\ga^*$ is not finite around $\zat{\vte}{\moma(\ttt)}=0$
 whence the \DCQ does not hold.
\end{example}

\begin{example}\label{Ex:Burg}\rm(a modification of the example in \cite[p.~263]{Bo.Le})
 Let $\mu$ be the Borel measure on $\TT=[0,1]$ with ${\mathrm d}\mu=2\ttt\,{\mathrm d}\ttt$,
 $\moma(\ttt)=(1,\ttt)$, $\ttt\in\TT$, and $\ga( t)=-\ln t$, $t>0$. By
 Theorem~\ref{T:domain}, $\dom{\vafu_{\ga}}$ consists of the pairs $a=(a_1,a_2)$
 such that $0<a_2<a_1$. Since $\ga^*(r)$ equals $-1-\ln(-r)$ for $r<0$ and
 $+\infty$ otherwise, for $\vte=(\vte_1,\vte_2)\in\R^2$
 \[
     \dufu_{\ga}(\vte)=
        \begin{cases}
                \;-1-\mbox{\large$\int_{0}^{1}$}\:2z\,\ln(-\vte_1-\vte_2z)\,{\mathrm d}z\,,
                    \quad &\vte_1\leq0\,,\;\vte_1+\vte_2<0\,,\\
                \;+\infty\,, &\text{otherwise,}
        \end{cases}
 \]
 where the integral is finite. The set $\Theta_\ga$ is equal to $\dom{\dufu_\ga}$,
 given by the two above inequalities. The dual problem takes the form
 \[
     \dufu_{\ga}^{*}(a)=\sup\nolimits_{(\vte_1,\vte_2)\in\Theta_\ga}\;
            \Big[\vte_1 a_1 + \vte_2 a_2 +1+\mbox{$\int_{0}^{1}$}\:
            2z\,\ln(-\vte_1-\vte_2z)\,{\mathrm d}z\Big]\,,  \qquad  a\in\R^2\,,
 \]
 and $\dom{\dufu_{\ga}^{*}}$ equals $\dom{\vafu_{\ga}}$ determined above,
 by Proposition~\ref{P:dualdom}. In the interior of~$\Theta_{\ga}$, the derivative
 $\vte_1\tfrac{\partial}{\partial\vte_1}+\vte_2\tfrac{\partial}{\partial\vte_2}$
 of the above bracket is equal to $\vte_1 a_1 + \vte_2 a_2+1$. Assuming
 $a\in\dom{\dufu_{\ga}^{*}}$, the derivative vanishes if and only if $\vte$
 belongs to the relatively  open segment between $\te=(0,-\frac{1}{a_2})$ and
 $(-\frac{1}{a_1-a_2},\frac{1}{a_1-a_2})$. The supremum can be restricted to this
 segment, parallel to the direction $(-a_2,a_1)$. The directional derivative of
 the bracket at $\te$ in this direction is equal to $a_2(2a_2-a_1)$, by a direct
 computation. Therefore, if $2a_2\leq a_1$ then $\te$ is the unique dual solution
 for $a$, not depending on $a_1$. In this case, the dual value $\dufu_{\ga}^{*}(a)$
 is equal to $-\tfrac{1}{2}-\ln a_2$ and the effective dual solution is
 $\efa\colon z\mapsto\frac{a_2}{z}$. When even $2a_2<a_1$, the first coordinate
 $\int_{0}^{1}\: 1\cdot\frac{a_2}{z}\cdot 2z\,{\mathrm d}z$ of the moment vector
 of \efa is less than $a_1$. Thus, the primal solution for~$a$ does not exist.
 Nevertheless, $\efa$ provides the generalized primal solution $\gea$.
\end{example}

\begin{example}\label{Ex:-infty}\rm
 Let $\TT=\R^2$, $\mu$ be the sum of the Lebesgue measure on the horizontal axis
 and the unit point mass at $(0,1)$, $\moma(\ttt)=(1,\ttt_1,\ttt_2)$ for
 $\ttt=(\ttt_1,\ttt_2)\in\TT$, and $\ga(t)=t\ln t$, $t\geq0$. The $\moma$-cone
 of $\mu$ is the convex hull of two of its proper faces $F=\{(t,0,t)\colon t\geq0\}$
 and $\{(0,0,0)\}\cup\{(t,r,0)\colon t>0\}$. Then, $\vafu_{\ga}(a)=\ga(t)$
 for $a=(t,0,t)\in F$ and $\vafu_{\ga}\equiv-\infty$ on $\cnmoma{\mu}\sm F$,
 using the fact that the Shannon functional can be explicitly evaluated
 at the Gaussian densities. Hence, $\dufu_{\ga}\equiv+\infty$ and the \DCQ
 fails. Nevertheless, by Definition~\ref{D:F-problems}, for $a\in\ri{F}$
 \[
    \dufu_{F,\ga}(\vte)=\inttt{\{\ttt_1=0\,,\,\ttt_2=1\}}
            e^{\vte_0+\vte_1\ttt_1+\vte_2\ttt_2-1}\,\mu({\mathrm d}\ttt_1,{\mathrm d}\ttt_2)
        =e^{\vte_0+\vte_2-1}\,,\qquad \vte=(\vte_0,\vte_1,\vte_2)\,,
 \]
 whence $\Theta_{F,\ga}=\R^3$. The $F$-dual problem for $a=(t,0,t)$, $t>0$, has
 many solutions, e.g.\ $(1+\ln t,0,0)$, and $\dufu_{F,\ga}^*(a)=\ga(t)$.
 By Theorem~\ref{T:primal.dual}, the primal solution $g_a$ exists and equals
 $g^*_{F,a}$, a function equal to $t$ at $(0,1)\in\TT$ and zero otherwise.
\end{example}

\begin{example}\label{Ex:gap}\rm
 On $\TT=[0,1]^2$, let $\mu$ be the sum of the Lebesgue measure
 and the unit masses at $(0,\frac13)$ and $(0,\frac23)$.
 Let $\moma(\ttt)=(1,\ttt_1,\ttt_2)$, $\ttt=(\ttt_1,\ttt_2)\in\TT$, and
 $\ga(t)=-2\sqrt{t}$ for $t\geq0$ whence $\ga^*(r)=-r^{-1}$ for $r<0$.
 If $a=(3,0,1)$ then $\vafu_\ga(a)=\fuc_\ga(g_a)=-2\sqrt{3}$
 where $g_a\in\calG_a$ is equal to~$3$ at~$(0,\frac13)$ and to~$0$ otherwise.
 If $\vte=(\vte_0,\vte_1,\vte_2)\in\R^3$ has $\vte_0$, $\vte_0+\vte_1$,
 $\vte_0+\vte_2$ and  $\vte_0+\vte_1+\vte_2$ negative then
 \[
    \dufu_\ga(\vte)=
     -\zlo{3}{3\vte_0+\vte_2}-\zlo{3}{3\vte_0+2\vte_2}
     -\inttt{\,[0,1]^2}\:\zlo{{\mathrm d}\ttt_1\,{\mathrm d}\ttt_2}{\vte_0+\vte_1\ttt_1+\vte_2\ttt_2}
 \]
 where the integral is finite. Otherwise, $\dufu_\ga(\vte)=+\infty$. Hence,
 the \DCQ holds. The maximization in the dual problem for $a=(3,0,1)$ includes
 the limiting $\vte_1\downarrow-\infty$, thus
 \[
   \dufu_\ga^{*}(a)=\sup\nolimits_{\,\vte_0<0,\;\vte_0+\vte_2<0}\;
                        \Big[\,3\vte_0+\vte_2
                        +\zlo{3}{3\vte_0+\vte_2}+\zlo{3}{3\vte_0+2\vte_2}\,\Big]\,.
 \]
 The bracket is increasing when $(\vte_0,\vte_2)$ moves in the
 direction $(1,-3)$ which implies that $\dufu_\ga^{*}(a)$ is equal
 to $\max_{\vte_2<0}\:[\vte_2+\frac{9}{2\vte_2}]=-3\sqrt{2}$.
 Hence, the primal value is strictly greater than the dual one.
 The sequence $-(\frac1n,n,\frac{3}{\sqrt{2}})$ is maximizing
 in the dual problem. By Theorem~\ref{T:GMLE}, the generalized dual
 solution $h_a$ is the limit in measure of the sequence of functions
 $(\frac1n+n\ttt_1+\frac{3}{\sqrt{2}}\ttt_2)^{-1}$.
 Since $h_a$ is equal to $\sqrt{2}$ at $(0,\frac13)$
 it differs from the primal solution $g_a$.
\end{example}

\begin{example}\label{Ex:moma=0}\rm
 Let $\moma\equiv\bnu$. Then $\calG_{a}$ consists of all $\calT$-measurable
 functions if $a=\bnu$, and is empty otherwise. Thus,
 $\dom{\vafu_\nm}\pdm\{\bnu\}$ for each $\nm\in B$. For the equality, i.\,e.,
 for the existence of a measurable function $g$ with $\fuc_{\nm}(g)<+\infty$,
 the obvious necessary condition $\int_{\TT}\,\inf_{t}\nm(\cdot,t)\,{\mathrm d}\mu<+\infty$
 is sufficient, as well, by Lemma~\ref{L:pomstand}, also implying that
 $\vafu_\nm(\bnu)$ equals that integral. Further,
 $\dufu_{\nm}(\vte)=\int_{\TT}\,\conm(\cdot,0)\,{\mathrm d}\mu$ for each $\vte\in\R^d$,
 whence $\dufu_\nm^*(\bnu)=-\int_{\TT}\,\conm(\cdot,0)\,{\mathrm d}\mu$. If the
 integral is finite then each $\vte\in\R^d$ is a dual solution for $a=\bnu$.
 Since $\conm(\cdot,0)$ is equal to $-\inf_{t}\nm(\cdot,t)$, in case
 $\vafu_\nm(\bnu)<+\infty$ the primal and dual values are equal,
 $\vafu_\nm(\bnu)=\dufu_\nm^*(\bnu)$. The latter may fail if
 $\vafu_\nm(\bnu)=+\infty$, for the adopted convention admits both
 $\pm\inf_{t}\nm(\cdot,t)$ to have integral $+\infty$.
 The finiteness of $\vafu_\nm(\bnu)$ is equivalent to the \PCQ for $a=\bnu$,
 in which case the primal solution for $a=\bnu$ exists if and only if
 $\inf_{t}\nm(\cdot,t)$ is attained $\mu$-a.e.\ \cite[Theorem 14.60]{RW}.
 This is equivalent to $\nm'(\cdot,+\infty)>0\,[\mu]$, thus
 $\bnu\in\Theta_\nm$ by \eqref{E:DCQ}, hence the mentioned result is contained
 in Corollary~\ref{C:prisolclassic}. By Theorem~\ref{T:geprisol},
 if the \DCQ fails then no generalized primal solution exists, either.
\end{example}

\begin{example}\label{Ex:Brepro}\rm
 Let $\ga\in\varGamma$ be differentiable except at $t=1$,
 $\mu$ be a pm on $(\TT,\calT)$, $d=1$ and $\moma\equiv1$.
 Then $\calG_a$ consists of the $\calT$-measurable functions
 whose $\mu$-integral equals~$a\in\R$. The Bregman projection
 of $h\equiv1$ to $\calG_a$ features the integrand
 \[
    [\ga h](s)=\trnarg{\ga}{s,1}
            =\ga(s)-\ga(1)-\ga'_{\sg{s-1}}(1)[s-1]\,,
                \qquad s\ge 0\,.
 \]
 For any $a>0$ the minimum subject to $g\in\calG_a$
 of $\Brega{g,h}=\int_\TT\,\trnarg{\ga}{g,1}\,{\mathrm d}\mu$ is attained when
 $g\equiv a$, by Jensen inequality. In other words, the Bregman projection
 of $h$ to $\calG_a$ exists and $g_{[\ga h],a}\equiv a$. For $0<a<1$ and
 $g\in\calG_a$ with $\int_\TT\,\ga(g)\,{\mathrm d}\mu$ finite,
 \[\begin{split}
    & \Brega{g,g_{[\ga h],a}}+\Brega{g_{[\ga h],a},h}\\
        &=\inttt{\TT}\:\big[\ga(g)-\ga(a)-\ga'(a)[g-a]\big]{\mathrm d}\mu
         +\ga(a)-\ga(1)-\ga'_-(1)[a-1]\\
         &=\inttt{\TT}\,\ga(g)\,{\mathrm d}\mu-\ga(1)-\ga'_-(1)[a-1]
                \end{split}
  \]
 while
 \[
    \Brega{g,h}=\inttt{\TT}\,\ga(g)\,{\mathrm d}\mu-\ga(1)-\ga'_-(1)[a-1]
             +\inttt{\{g>1\}}\:[\ga'_-(1)-\ga'_+(1)][g-1]{\mathrm d}\mu\,.
 \]
 This shows that $ \Brega{g,h}<\Brega{g,g_{[\ga h],a}}+\Brega{g_{[\ga h],a},h}$
 when $g\in\calG_a$ and the set $\{g>1\}$ is not $\mu$-negligible.
\end{example}

 \newcommand{\popa}[1]{|#1|_{+}}
\begin{example}\label{E:threeelem}\rm
 Let $\mu$ be the counting measure on $\TT=\{1,2,3\}$, and $\moma$ have the
 values $\moma(1)=(1,1)$, $\moma(2)=(1,-1)$ and $\moma(3)=(1,0)$. Functions
 $g$ on~$\TT$ are identified with points in~$\R^3$. Thus,
 \[
   \calG_a=\big\{(\tfrac{t+a_2}{2},\tfrac{t-a_2}{2},a_1-t)\colon t\in\R\big\}\,,
       \qquad a=(a_1,a_2)\in\R^2\,,
 \]
 and the \moma-cone $\cnmoma{\mu}\pdm\R^2$ is given by $|a_2|\leq a_1$.
 If $\ga(t)=\frac{\;t^2}{2}$, $t\geq0$, then $\Theta_\ga=\R^2$ and,
 using that $(\ga^*)'=\popa{\cdot}$, the family
 \[
   \FF_\ga=\big\{f_\vte=(\popa{\vte_1+\vte_2}\,,\popa{\vte_1-\vte_2}\,,\popa{\vte_1})\,
       \colon\, \vte=(\vte_1,\vte_2)\in\R^2\,\big\}
 \]
 is the union of three two-dimensional cones in~$\R^3$. The set
 $\Theta_\ga^+$ of $\vte\in\Theta_\ga$ with $\zat{\vte}{\moma}\geq\ga'(0)=0$
 is determined by $|\vte_2|\le\vte_1$ and coincides with one of the three cones.
 Theorem~\ref{T:CLASSIC} fails if $\Theta_\ga^+$ is replaced by the whole $\Theta_\ga$,
 noting that the Bregman distance ${\mathsfsl  B}_{\ga}$ of two points $g,h\in\R^3$
 equals their squared Euclidean distance divided by two. Lemma~\ref{L:classic} fails
 as well, e.\,g.,\ if $\teob=(0,1)\not\in\Theta_\ga^+$ then $f_\teob=(1,1,0)$,
 \[\begin{split}
    f_{\vte+\teob}&=(\popa{\vte_1+\vte_2+1}\,,\popa{\vte_1-\vte_2-1}\,,\popa{\vte_1})\,,\\
    f_{[\nm f_\teob],\vte}&=(\popa{\vte_1+\vte_2+1}\,,\popa{\vte_1-\vte_2+1}\,,\popa{\vte_1})\,,
 \end{split}\]
 by Lemma~\ref{L:gat}.
\end{example}

\begin{example}\label{E:int.big}\rm
 On $\TT=(1,+\infty)$ let $\mu$ be the pm with density $2\ttt^{-3}d\ttt$,
 $\moma(\ttt)=(1,\ttt)$, $\ttt\in\TT$, and $\ga(t)=\frac12 t^{2}$, $t\geq0$.
 Then, $\cnmoma{\mu}=\{(a_1,a_2)\colon a_2>a_1>0\}\cup\{(0,0)\}$ coincides
 with $\dom{\vafu_\ga}$. %The \PCQ holds for $a$ if and only if $a_2>a_1>0$.
 For $\vte=(\vte_1,\vte_2)\in\R^2$,
 \[
    \dufu_\ga(\vte)=
        \begin{cases}
                \;\mbox{\large$\int_{1}^{+\infty}$}\,|\vte_1+\vte_2 z|_+^2\;\frac{2}{\ttt^3}\,{\mathrm d}\ttt\,,
                    \quad &\displaystyle\vte_2\leq0\,,\\
                \;+\infty\,, &\text{otherwise.}
        \end{cases}
 \]
 Hence, $\Theta_\nm$ is also given by the above inequality. The function
 $f_\vte(z)=\popa{\vte_1+\vte_2 z}$, $\vte\in\Theta_\nm$, identically
 vanishes on $\TT$ if $\vte_1\le-\vte_2$. Otherwise, $\vte=t(1,-r)$
 with $t>0$ and $0\le r<1$, and $f_\vte$ has the moment vector
 $t\big((1-r)^2,\; 2(1-r+r\ln r)\big)$. By a straightforward calculation,
 the moment vectors of $f_\vte$, $\vte\in\Theta_\nm$, exhaust the
 subcone of $\cnmoma{\mu}$ given by $2a_1\ge a_2$. It follows that for
 $a$ in this cone the primal solution $g_a$ exists.

 In the case $a_2>2a_1>0$ Proposition~\ref{P:Breclo} is employed.
 Since $\Brega{g,h}=\frac{1}{2}\norm{g-h}_{L_2(\mu)}^2$ for $g,h$
 nonnegative measurable, Bregman closure equals {the} $L_2(\mu)$-closure.
 The $L_2(\mu)$-closure of $\calG_a^+$ with $0<a_1<a_2$
 contains each function $f_\vte$, $\vte\in \Theta_\ga$, whose moment
 vector $(b_1,b_2)$ satisfies $b_1=a_1,\,b_2<a_2$. Indeed, for $n\geq1$
 there exists $x_n>1$ such that
 \[
    (a_2-b_2)n=
    {\intt}_{\!\!\!\!\!\scriptscriptstyle1}^{\scriptscriptstyle\!\!x_n}\zlo{{\mathrm d}z}{z\ln z}-
        \zlo{b_2}{a_1}\,{\intt}_{\!\!\!\!\!\scriptscriptstyle1}^{\scriptscriptstyle\!\!x_n}\zlo{{\mathrm d}z}{z^2\ln z}\,,
 \]
 by continuity. Then, $x_n\uparrow+\infty$, the function
 \[
   \ttt\mapsto \Big[1-\zlo{1}{a_1 n}{\intt}_{\!\!\!\!\!\scriptscriptstyle1}^{\scriptscriptstyle\!\!x_n}\:
        \zlo{{\mathrm d}r}{r^2\ln r}\Big]\,f_\vte(\ttt)+\zlo{\ttt}{2n\ln\ttt}\egy_{(1,x_n)}(\ttt)
 \]
 belongs to $\calG_a$, and $L_2(\mu)$-converges to $f_\vte$ as $n\to +\infty$ because
 $\frac{z}{\ln z}\in L_2(\mu)$.

 The set $\Theta_\nm^+$ consists of those $\vte\in\Theta_\ga$ for which
 $\vte_1+\vte_2 z\ge 0,\;z\in \TT$. Thus, it is given by
 $\vte_1\ge 0$ and  $\vte_2=0$. For such $(\vte_1,\vte_2)$, the function
 $f_\vte$ equals identically~$\vte_1$ and has the moment $(\vte_1,2\vte_1)$.
 By the above result, if $a_2>2a_1>0$ then the $L_2(\mu)$-closure of $\calG_a$
 contains $f_\vte$ with $\vte_1=a_1$. Hence, Proposition~\ref{P:Breclo}
 implies that the effective dual solution $g_a^*$ equals the constant
 $a_1$. Note that the last assertion of Proposition~\ref{P:Breclo} fails
 if the restriction $\vte\in\Theta_\nm^+$ is dropped.
\end{example}

%11 11 11 11 11 11 11 11 11 11 11 11 11 11 11 11 11 11 11 11 11 11 11 11 11 11
\section{Relation of this work to previous ones}\label{S:relation}

 The subject addressed in this paper has one of its origins in
 the principle of maximum entropy (\ME) which comes from statistical physics
 and has been promoted as a general principle of inference primarily by
 Jaynes \cite{Jay} and Kullback \cite{Kull}. While \ME calls for maximizing
 Shannon entropy or for minimizing $I$-divergence (Kullback--Leibler
 distance~\cite{KulL}), maximization of Burg entropy~\cite{Burg0,Burg}
 and other `entropy functionals' is also widely used in sciences. These
 applications motivated the formulation of the general minimization
 problem in eq.~\eqref{E:priva} with autonomous integrands. It is for
 convenience that minimization of convex integral functionals is addressed,
 maximization of concave ones as Shannon or Burg entropy is covered by taking
 their negatives.

 The literature of the subject is extensive, some pointers are given here
 to works that have influenced ours. The integral functional~\eqref{E:fuc}
 with an autonomous integrand $\nm(\ttt,t)=\varPhi(t)$ is called $\varPhi$-entropy
 in~\cite{Bo.Le.conv1,Te.Va}, where for a growing number of constraints,
 convergence of the solutions of the problem~\eqref{E:priva} is established
 under suitable conditions.  Divergences of form~\eqref{E:gadiv} have been
 introduced by Csisz\'ar~\cite{Csi.fdiv,Csi.fdiv2}, called $f$-divergences,
 and by Ali and Silvey~\cite{Ali}. They, as well as $I$-divergence, were
 originally defined for probability densities $g$ and $h$ only. Substantial
 developments in their theory and applications are due to I.~Vajda, see
 e.\,g.~\cite{LieVa}. More recent references include
 \cite{Csi.Ga.Ga, Bro.Ke, Am.C}. Bregman distances were introduced
 in~\cite{Bre} as non-metric distances between vectors in~$\R^d$, associated
 with a convex function on~$\R^d$, for numerous applications in convex
 programming problems see the book~\cite{CZ}. The subclass of separable
 Bregman distances is the one whose infinite dimensional extension is used
 in this paper, see Remark~\ref{R:Bre}. Their statistical applications, initiated
 by Jones and Byrne~\cite{Jon.Byr}, are currently wide ranging, see e.\,g.\ Murata
 et al~\cite{MTKA}. For more general Bregman distances see e.\,g.~\cite{BBC,Fri.Sri.Gup}.
 The axiomatic study~\cite{Csi.ax} highlights (in the finite dimensional case)
 the distinguished role of $\ga$-divergences and separable Bregman distances,
 and primarily that of $I$-divergence. The problem of minimizing convex
 integral functionals arises also in large deviations theory, for interplay
 with this field see e.\,g.~\cite{Da-Cs.Ga,GG,Csi.Ga.Ga,Leo.domin}.
 Other fields could also be mentioned, such as in control theory
 where, typically, also derivatives of the unknown function are involved.
 Regarding a possible interplay see \cite{BLN}.

 A large body of the literature on the minimization problem of eq.~\eqref{E:priva}
 is application oriented and mathematically non-rigorous. For example,
 while the form of the solution like~\eqref{E:fte} is derived via
 Lagrange multipliers, little attention is payed to conditions under
 which a solution exists, and is indeed of this form when it does. Early
 rigorous results about $I$-divergence minimization and associated Pythagorean
 (in)equalities  were obtained by Chentsov~\cite{Chen} and Csisz\'ar~\cite{Csi.I-div}.
 Recent works typically employ convex duality, following the lead of
 Borwein and Lewis, see \cite{Bo.Le.dua,Bo.Le.conv1,Bo.Le}. Previously,
 convex duality had been applied to the $I$-divergence minimization problem
 in~\cite{Be-Ch}. Advanced tools from functional analysis appear indispensable
 to efficiently deal with the case, not treated here, when the range of the
 moment mapping is infinite dimensional, see L\'eonard \cite{Leo2}--\cite{Leo.domin}.
 L\'eonard's results are strong and general also when restricted to finite
 dimensional mappings.  Still, they appear to require assumptions on the integrand
 and the moment mapping not needed here, e.\,g., that $\nm(\ttt,t)$ is nonnegative
 and equals $0$ for some $t=t_{\ttt}$. Another tools are provided by
 differential geometry \cite{Amari, Am.C}, first applied in the \ME context
 by Chentsov~\cite{Chen}. They require strong regularity conditions but lead
 to impressive `geometric' results for example about Pythagorean identities,
 going beyond those obtainable otherwise.

 This paper generalizes the results obtained for the Shannon case
 in~\cite{Csi.Ma.proj}. Convex duality is used, as there, for the value
 function only, a convex conjugate of the integral functional is not
 needed. Accordingly, the functional is not restricted to a `good' space
 that has a manageable dual space, as frequently done in the literature.
 A key tool in~\cite{Csi.Ma.proj} has been the convex core of a measure on
 $\R^d$, introduced in~\cite{Csi.Ma.cc}. In the present generality, its
 role is played by the concept of conic core, introduced here.
 The framework is in several respects more general than usual:
 (i) non-differentiable integrands are allowed (ii) the
 integrands need not be autonomous (iii) there are no restrictions
 beyond measurability, neither on the functions $g$ over which the
 functional is minimized nor on the moment mapping \moma, other than
 that \moma has finite dimensional range. While feature (i) is not unique
 for this paper, in the literature often stronger assumptions are adopted
 on $\nm$ than differentiability (in addition to strict convexity which
 is assumed also here). Typical ones are essential smoothness plus
 cofiniteness, or the equivalent assumption that $\conm$ is strictly
 convex on $\R$, as in~\cite{MTKA}. Non-autonomous integrands (also admitted
 in~\cite{Leo,Leo.domin}) do not cause conceptual difficulties but do
 cause technical ones concerning measurability. These are handled here
 via the theory of normal integrands initiated in \cite{Rock.cofu, Rock.cofu.dual}
 and summarized in detail in the recent book~\cite{RW}. The latter is relied
 upon in the paper also elsewhere.

 The limitations of our framework are, in addition to restricting
 attention to moment mappings of finite dimensional range, that only
 equality constraints are considered, and the integrand value $\nm(\ttt,t)$
 has to be finite for $t>0$ and $+\infty$ for $t<0$. To consider only
 equality constraints does not seem a serious restriction. It should
 not be difficult to extend the results to constraints of the form that the
 moment vector belongs to a convex subset of $\R^d$, see e.\,g.~\cite{Csi.Ga.Ga}.
 Our restriction on the integrand is not needed for the mere extension of
 familiar results to the generality of (i)--(iii). It is, however, essential for
 the main results, viz.\ the geometric characterization of the effective domain
 of the value function and its implications that extend several results previously
 proved only under the \PCQ beyond that assumption. These main
 results are relevant in those cases when the effective domain includes
 a nontrivial boundary, which is typical when the underlying measure $\mu$
 has discrete components. In `classical' moment problems involving Lebesgue
 measure on $\R^k$ and moment mappings formed by polynomials or trigonometric
 polynomials, the need for going beyond the \PCQ does not arise (while
 the problem of nonexistence of primal solution does).

 \enlargethispage{.3cm} As genuine primal and dual solutions do not always exist, also generalized
 ones are studied which are universal limits of minimizing (maximizing) sequences.
 Another concept of generalized solution, not used in this paper, involves relaxation
 of the problem~\eqref{E:priva} to minimization over a larger space. The latter is
 typically obtained regarding the functional to be defined over a specific linear
 space with a manageable dual space, which paves a road to extend the functional
 to the second dual. In~\cite{Leo2,Leo3,Leo} this extension is to the topological
 dual of an Orlicz space. In~\cite{Leo.domin} it is also shown in considerable
 generality that the `absolutely continuous component' of this kind of generalized
 solution coincides with the generalized solution in our sense; a similar but more
 special result appeared previously in \cite{Csi.genproj}. Note that nonexistence
 of a minimizer in eq.~\eqref{E:priva} had emerged as a practical problem
 in the context of three dimensional density reconstruction via Burg entropy
 maximization, see references in \cite{Bo.Le} where the mathematical background
 of this phenomenon has been clarified.

 By key results of this paper, generalized primal and dual solutions
 in our sense exist, subject to the \DCQ, in all nontrivial cases, and
 their Bregman distance is a lower bound to the duality gap. The adopted concept
 of  generalized solution dates back to Topsoe~\cite{Top} who, for Shannon
 entropy maximization over a convex set of probability distributions,
 established a Pythagorean inequality involving a `center of attraction'
 perhaps not in that set. Actually, the existence of generalized $I$-projections
 is implicit already in \cite{Csi.I-div}. They were studied in detail in
 \cite{Csi.Sanov,Csi.Ma.proj}; in~\cite{Csi.Sanov} also their relevance for
 large deviations theory is demonstrated. Generalized minimizers of integral
 functionals with (differentiable) autonomous integrands $\ga\in\varGamma$,
 for arbitrary convex sets of functions, were introduced and corresponding
 Pythagorean inequalities established in~\cite{Csi.genproj}. Generalized
 primal solutions for the minimization problem in eq.~\eqref{E:priva}, assuming the \PCQ and implicitly
 the \DCQ, are treated in \cite{Csi.Ga.Ga}. Note that the existence result
 \cite[Theorem~1(c)]{Csi.genproj} does not hold in full generality, its proof
 contains a gap of implicitly assuming that functions in a minimizing sequence
 always have bounded Bregman distance from some fixed function. For the minimization
 problem in eq.~\eqref{E:priva} the latter is true if the \DCQ holds, due to Lemma~\ref{L:KEYID},
 hence the existence of a generalized primal solution subject to the \DCQ does follow
 from \cite[Theorem~1(c)]{Csi.genproj}. In this respect, the new feature of the
 result here is that it explicitly specifies the generalized primal solution.

 Generalized dual solutions for the $I$-divergence case, viz.\ generalized
 maximum likelihood estimates, have been treated in \cite{Csi.Ma.infdim};
 stronger results including their explicit description appear in~\cite{Csi.Ma.gmle}.
 In the present paper the existence of generalized dual solution is proved,
 for any integral functional with integrand $\nm\in B$ subject to the \DCQ,
 via a nontrivial updating of the technique used in \cite{Csi.Ma.infdim},
 similar to that in \cite{Csi.genproj}, going back to \cite{Top}. This proof
 works also when the moment mapping has infinite dimensional range, but
 gives no indication how to construct the generalized dual solution. The
 latter remains an open problem. Note that when the duality gap is zero,
 a direct proof shows that, subject to the \DCQ, the generalized dual
 solution exists and coincides with the generalized primal solution,
 which has been explicitly described.

 \enlargethispage{.3cm} Finally, let us comment on the family of functions $\FF_\nm=\{f_\vte\colon\vte\in\Theta_\nm\}$.
 It has been observed many times that $\FF_\nm$ plays the same role as a classical
 exponential family does in the Shannon case, see Appendix~B for formal details.
 For statistical applications see e.\,g.~\cite{MTKA} where our $\FF_\nm$ is
 called $U$-model (with $U=\conm$). For the classical theory of exponential
 families see~\cite{Chen,BaNi}. For a very general concept see~\cite{Da.Gr}.
 Apparently, it has not been pointed out before that some key properties
 of exponential families extend to $\FF_\nm$ only when $\nm$ is essentially
 smooth. If $\nm$ is merely differentiable, those properties extend only
 partially, to functions in $\FF_\nm$ parameterized by $\vte$ in a specific
 set $\Theta_\nm^+\subseteq \Theta_\nm$, see Proposition~\ref{P:Breclo}
 and Theorems~\ref{T:CLASSIC},~\ref{T:classic}.   \pagebreak
%\\[1cm]

\begin{large} {APPENDICES} \end{large}
\appendix

%AAAAAAAAAAAAAAAAAAAAAAAAAAAAAAAAAAAAAAAAAAAAAAAAAAAAAAAAAAAAAAAAAAAAAAAAAA
\section{Integral representation}\label{S:intrep}

 The proof of Theorem~\ref{T:repre} depends on an interchange of
 minimization and integration, see below. We are not aware of
 a reference that would guarantee admissibility of this interchange
 in the required generality, though in the special case of finite
 $\mu$ and $\mu$-integrable $\moma$,  \cite[Theorem 14.60]{RW}
 suffices. An extension of the latter, Theorem~\ref{T:infint},
 will be proved below and applied to cover the general case.

 A linear space $\calH$ of real $\calT$-measurable functions
 is \emph{decomposable} w.r.t.~$\mu$ \cite[Definition~14.59]{RW}
 if $g\egy_{\TT\sm A}+h\egy_{A}$ belongs to~\calH whenever $g\in\calH$,
 $A\in\calT$ has finite $\mu$-measure and $h$ is bounded \calT-measurable.
 For $\mu$ finite, \calH is decomposable if and only if it contains all
 bounded $\calT$-measurable functions. A weaker notion is introduced
 as follows.

\begin{definition}\label{D:sigma}\rm
 A space $\calH$ is \emph{$\sigma$-decomposable} w.r.t.~$\mu$
 if $Z$ can be covered by a countable family of sets $Z_n\in\calT$
 with $\mu(Z_n)$ finite such that $g\egy_{Z\sm A}+h\egy_{A}\in\calH$
 whenever $g\in\calH$, $A\in\calT$ is contained in some
 $Z_n$ and $h$ is bounded  \calT-measurable.
\end{definition}

\begin{remark}\label{R:nlog}\rm
 There is no loss of generality in assuming that $Z_n\pdm Z_{n+1}$
 in Definition~\ref{D:sigma}.
\end{remark}

\begin{remark}\label{L:sideco}\rm
 If $\moma\colon\TT\to\R^d$ is any moment mapping then the space $\calG$, consisting of
 those functions $g$ for which the moment vector $\int_\TT\: \moma g\,{\mathrm d}\mu$ exists,
 is $\sigma$-decomposable w.r.t.~$\mu$. Namely, by $\sigma$-finiteness, $\TT$ can be
 covered by sets $Y_n\in\calT$ of finite measure, and the countable family
 $Z_{n,m}=Y_n\cap\{\norm{\moma}\leq m\}$ indexed by $n,m$ is suitable.
 In fact, if $A\in\calT$ is contained in some $Z_{n,m}$ and $h$ is bounded
 \calT-measurable then  $g\in\calG$ implies existence of
 $\int_{\TT\sm A}\: \moma g\,{\mathrm d}\mu+\int_{A}\: \moma h\,{\mathrm d}\mu$.
 The space \calG need not be decomposable w.r.t.~$\mu$ even if $\mu$ is finite.
 For example, if $\TT=\R$,
 ${\mathrm d}\mu=\frac{{\mathrm d}\ttt}{1+\ttt^2}$ and $\moma(\ttt)=\ttt$ then
 $\calG$ does not contain the constant functions.
\end{remark}

 The following assertion on interchange of minimization and integration
 is an extension of \cite[Theorem 14.60]{RW} to the $\sigma$-decomposable
 spaces.

\begin{theorem} \label{T:infint}
   Let $\calH$ be a $\sigma$-decompos\-able linear space of $\calT$-measurable
   functions on a $\sigma$-finite measure space $(\TT,\calT,\mu)$,
   and let $\alpha\colon\TT\times\R\to[-\infty,+\infty]$ be a normal integrand
   such that the integral functional $\fuc_\alpha(g)\triangleq\int_{\TT}\:\alpha(\ttt,g(\ttt))\;\mu({\mathrm d}\ttt)$
   does not identically equal $+\infty$ for $g\in\calH$. Then
   \begin{equation}\label{E:RW}
       \inf_{g\in\calH} \inttt{\TT}\:\alpha(\ttt,g(\ttt))\;\mu({\mathrm d}\ttt)
        = \inttt{\TT}\:\inf_{t\in\R}\:\alpha(\ttt,t)\;\mu({\mathrm d}\ttt)\,.
   \end{equation}
\end{theorem}

\noindent
 The function on $\TT$ that is integrated on the right is denoted by
 $\alpha_{\text{inf}}$. It is \calT-measurable by \cite[Theorem~14.37]{RW}.
 The following lemma does not involve $\calH$.

\begin{lemma}\label{L:pomstand}
    If $\int_\TT\,\alpha_{\text{\emph{inf}}}\,{\mathrm d}\mu<t$ then
    \mbox{$\fuc_\alpha(h)<t$} for some real $\calT$-measurable function $h$.
\end{lemma}

\begin{Proof}
  This is proved neatly on the lines 7--17 of the proof of \cite[Theorem 14.60]{RW}.
\end{Proof}

\begin{lemma}\label{L:stand}
    If $\calH$ is $\sigma$-decomposable, $\fuc_\alpha(f)$ finite for some
    $f\in\calH$, and \mbox{$\int_\TT\,\alpha_{\text{\emph{inf}}}\,{\mathrm d}\mu<t$}
    then $\fuc_\alpha(g)<t$ for some $g\in\calH$.
\end{lemma}

\begin{Proof}
 By Lemma~\ref{L:pomstand}, $\fuc_\alpha(h)<t$ for some
 $\calT$-measurable $h$. Let $0<2\vare<t-\fuc_\alpha(h)$.
 By $\sigma$-decomposability and
 Remark~\ref{R:nlog}, $Z$ is covered by a countable increasing
 sequence $Z_n$ with $\mu(Z_n)$ finite that has the property
 from Definition~\ref{D:sigma}. Since $\fuc_\alpha(f)$ is finite,
 for $n$ sufficiently large
 \[
        \vare\geq\inttt{Z\sm Z_n}\:\alpha(\ttt,f(\ttt))\;\mu({\mathrm d}\ttt)\quad\text{and}\quad
        t-2\vare>\inttt{Z_n}\:\alpha(\ttt,h(\ttt))\;\mu({\mathrm d}\ttt)\,.
 \]
 For $m$ sufficiently large
 \[
        \vare\geq\inttt{Z_n\sm\{h\leq m\}}\:\alpha(\ttt,f(\ttt))\;\mu({\mathrm d}\ttt)\quad\text{and}\quad
        t-2\vare>\inttt{Z_n\cap\{h\leq m\}}\:\alpha(\ttt,h(\ttt))\;\mu({\mathrm d}\ttt)\,.
 \]
 Since $h$ is bounded on $A=Z_n\cap\{h\leq m\}\pdm Z_n$, the function $g=f\egy_{\TT\sm A}+h\egy_{A}$
 belongs to $\calH$ by $\sigma$-decomposability. Combining the above inequalities, $\fuc_\alpha(g)<t$.
\end{Proof}

\bigskip

\noindent
P\,r\,o\,o\,f\ o\,f\ T\,h\,e\,o\,r\,e\,m~\ref{T:infint}.
 The inequality $\geq$ in \eqref{E:RW} follows from
 $\alpha(\ttt,g(\ttt))\geq\alpha_{\text{inf}}(\ttt)$, $\ttt\in\TT$,
 by integration. By assumption, $\fuc_\alpha(f)<+\infty$ for some
 $f$ in $\calH$. If $\fuc_\alpha(f)=-\infty$ then \eqref{E:RW}
 has $-\infty$ on both sides. Otherwise, the inequality $\leq$
 in \eqref{E:RW} follows from Lemma~\ref{L:stand}.
\hfill$\square$

\bigskip

\noindent
P\,r\,o\,o\,f\ o\,f\ T\,h\,e\,o\,r\,e\,m~\ref{T:repre}.
 By definition,
 \[
    \vafu_{\nm}^*(\vte)=-{\inf}_{a\in\R^d}\:
                                \big[\,-\zat{\vte}{a}+{\inf}_{g\in\calG_{a}}\:\fuc_{\nm}(g)\,\big]\,,
                                        \qquad\vte\in\R^d\,.
 \]
 The expression on the right rewrites to
 \[
   -  {\inf}_{a\in\R^d}\:{\inf}_{g\in\calG_{a}}\:\inttt{\TT}\:
            \big[\,-\zat{\vte}{\moma(\ttt) g(\ttt)}+\nm(\ttt,g(\ttt))\,\big]\mu({\mathrm d}\ttt)
 \]
 where two infima reduce to one, over $g\in\calG$. By Remark~\ref{L:sideco},
 the linear space $\calG$ defined by any moment mapping is $\sigma$-decomposable.
 Hence, by Theorem~\ref{T:infint}, if $\vafu_{\nm}\not\equiv+\infty$ then
 \[
    \vafu_{\nm}^*(\vte)=-\inttt{\TT}\:
            {\inf}_{t\in\R}\big[\,-\zat{\vte}{\moma(\ttt)}t+\nm(\ttt,t)\,\big]\mu({\mathrm d}\ttt)\,,
                                        \qquad\vte\in\R^d\,,
 \]
 and thus $\vafu_{\nm}^*=\dufu_{\nm}$.
\hfill$\square$

\bigskip

 The assumption in Theorem~\ref{T:repre} that $\vafu_{\nm}$ is not identically
 $+\infty$ does matter, see Example~\ref{Ex:vafu.dual}.

%BBBBBBBBBBBBBBBBBBBBBBBBBBBBBBBBBBBBBBBBBBBBBBBBBBBBBBBBBBBBBBBBBBBBBBBBBBBBBBBBBBBBB
 \section{Restricted value function}\label{S:revafu}

 In this appendix, some details are discussed for the Shannon's integral
 functional defined by the autonomous integrand $\ga(t)=t\ln t$, $t>0$.
 If $g$ is nonnegative and $\int_\TT\,g\,{\mathrm d}\mu=t>0$ then
 $\fuc_{\ga}(g)=\ga(t)+t\fuc_{\ga}(g/t)$. Hence, the integral functional
 $\fuc_{\ga}$ is determined by its values on the probability densities
 $g$ w.r.t.~$\mu$. For such a density, $\fuc_\ga(g)$ is the negative entropy
 of the corresponding probability measure w.r.t.~$\mu$, or its
 $I$-divergence from $\mu$ when $\mu(\TT)=1$.

 Further, it is assumed that the moment mapping $\moma$ has first coordinate
 identically equal to~$1$, which is nonrestrictive in many applications, see
 Section~1.B.
 Let $\moma=(1,\psi)$ where $\psi\colon\TT\to\R^{d-1}$.
 For any vector $a\in\R^d$ with a positive first component, writing it as $(t,b)$
 with $t>0$ and $b\in\R^{d-1}$, $\vafu_{\ga}(a)=\ga(t)+ t\,\vafu_{\ga}(1,b/t)$.
 Hence, the value function $\vafu_{\ga}$ is uniquely determined by its restriction.
 Let $\revafu_{\ga}(b)\triangleq\vafu_{\ga}(1,b)$, $b\in\R^{d-1}$.
 The conjugate of the value function at any point $\vte=(r,\te)$, where $r\in\R$
 and $\te\in\R^{d-1}$, is
 \begin{equation}\label{E:pocit}
 \begin{split}
    \vafu^*_{\ga}(\vte)
        &={\sup}_{t>0,\,b\in\R^{d-1}}\:
            \big[rt+\zat{\te}{b}-\ga(t)- t\,\revafu_{\ga}(b/t)\big]\\
        &={\sup}_{t>0}\:
            \big[rt-\ga(t)+ t\revafu^*_{\ga}(\te) \big]
        =\ga^*(r+\revafu^*_{\ga}(\te))=\exp[r+\revafu^*_{\ga}(\te)-1]
 \end{split}\end{equation}
 using that the convex conjugate of $\ga$ is $\ga^*(r)=e^{r-1}$, $r\in\R$.
 By Theorem~\ref{T:repre}, knowing that $\fuc_{\ga}(g)=0$ for $g\equiv 0$,
 $\vafu^*$ admits the integral representation, and hence
 \[
   \revafu^*_{\ga}(\te)=1-r+\ln \vafu^*_{\ga}(r,\te)
    =\ln\inttt{\TT}\; e^{\zat{\te}{\psi}}\,{\mathrm d}\mu\,,\quad \te\in\R^{d-1}\,.
 \]
 This formula is well-known and has been a key tool when minimizing the negative
 Shannon entropy $\fuc_\ga(g)$ of a probability density $g$ subject to
 moment constraints, e.\,g.~in \cite{Csi.Ma.proj}.

 The set $\Theta_{\ga}$ equals $\dom{\vafu^*_{\ga}}$, consisting of all $\vte\in\R^d$
 with $e^{\zat{\vte}{\moma}}$ $\mu$-integrable. The functions \eqref{E:fte} of the
 family $\FF_{\ga}$ are given by $f_\vte=e^{r+\zat{\te}{\psi}-1}$ where $\vte=(r,\te)$.
 The family of $f_\vte$ that integrate to $1$ is known as the exponential family
 based on $\mu$ with canonical statistic $\psi$, see \cite{BaNi,Chen}.

 The original and restricted dual problems are also simply related, for
 $a=(t,b)$ with $t>0$
 \[\begin{split}
    \vafu^{**}_{\ga}(a)&={\sup}_{\te\in\R^{d-1}}
        \Big[\zat{\te}{b}+{\sup}_{r\in\R}
        \big[rt-\ga^*(r+\revafu^*_{\ga}(\te))\big]\Big]\\
        &={\sup}_{\te\in\R^{d-1}}
        \big[\zat{\te}{b}- t\revafu^*_{\ga}(\te)+\ga(t)\big]
        =\ga(t)+t\,\revafu^{**}_{\ga}(b/t)\,,
 \end{split}\]
 using \eqref{E:pocit}. For example, the duality gap of the original problem at $a$
 is $t$ times the duality gap $\revafu^{**}_{\ga}(b/t)-\revafu_{\ga}(b/t)$
 of the restricted problem at $b/t$.

 Note that other integral functionals do not admit such simple formulas that would relate
 unrestricted and restricted value functions, and their conjugates and biconjugates.

%CCCCCCCCCCCCCCCCCCCCCCCCCCCCCCCCCCCCCCCCCCCCCCCCCCCCCCCCCCCCCCCCCCCCCCCCCCCCCCCCCCC
\section{{\large$\gamma$}-divergences}\label{S:gadiv}

 Let $\ga\in\varGamma$ be nonnegative with $\ga(1)=0$. The $\ga$-divergence
 of a function $g:\TT\to [0,+\infty)$ from $h:\TT\to (0,+\infty)$, both
 $\calT$-measurable, is defined by\index{$\ga$-divergence}\label{ga divergence 2}
 \begin{equation}\label{E:gadiv}
      \gadiv{g}{h}\triangleq
      \inttt{\TT}\:h\,\ga(g/h)\;{\mathrm d}\mu\,.
 \end{equation}
 This divergence is nonnegative, and equals $0$ only if $g=h\;[\mu]$. If
 $\ga(t)=t\ln t-t+1$ then $\gadiv{g}{h}$ is equal  to the $I$-divergence
 of $g$ from $h$. As mentioned in subsection~1.B., the minimization of a
 $\ga$-divergence $\gadiv{g}{h}$ subject to $g\in\calG_{a}$, for fixed $h$,
 is a frequently studied instance of the minimization problem addressed
 in this paper. The integrand
 $\nm\colon(\ttt,t)\mapsto h(\ttt)\,\ga (t/h(\ttt))$, $\ttt\in\TT$,
 $t\in\R$, in \eqref{E:gadiv} is in general non-autonomous, but belongs to $B$
 and $\gadiv{g}{h}=\fuc_{\nm}(g)$. Lemma~\ref{L:change} below shows that
 a reformulation of the minimization with an autonomous integrand is possible.

 A general idea is to modify simultaneously a measure $\mu$ and integrand
 $\nm\in B$ to $\tilde\mu$ and $\tilde\nm$ given by
 \[
    {\mathrm d}\tilde\mu=h\,{\mathrm d}\mu \quad\text{~~and~~}\quad
    \tilde\nm(\ttt,t)=\nm(\ttt,t\,h(\ttt))/h(\ttt)\,,\quad \ttt\in\TT\,, t\in\R\,,
 \]
 where $h$ is a given positive \calT-measurable function.
 By \cite[Proposition~14.45]{RW}, $\tilde\nm\in B$. As earlier,
 the dependence on $\mu$ in $\fuc$, $\vafu$, $\calG_a$ is
 added to indices while the moment mapping $\moma$ is not changed.

\begin{lemma}\label{L:change}
 Given a positive \calT-measurable function $h$, let $\tilde g=g/h$
 for any \calT-measur\-able function $g$. Then
 $\fuc_{\mu,\nm}(g)=\fuc_{\tilde\mu,\tilde\nm}(\tilde g)$ and
 $\int_{\TT}\: g\moma\:{\mathrm d}\mu =\int_{\TT}\:\tilde g\moma\:{\mathrm d}\tilde\mu$
 if one of the integrals exists. Further, $g\in\calG_{\mu,a}$
 if and only if $\tilde g\in\calG_{\tilde\mu,a}$, and
 $\vafu_{\mu,\nm}=\vafu_{\tilde\mu,\tilde\nm}$.
\end{lemma}

\noindent A simple proof based on substitutions in integrals is omitted.

 When considering the minimization of $\gadiv{g}{h}$ subject to a
 moment constraint on~$g$, Lemma~\ref{L:change} applies with the
 very function $h$ from the divergence. The corresponding integrand
 $\tilde\nm$ is autonomous and coincides with $\ga$. Hence, the
 minimization of $\gadiv{\cdot}{h}$ over $\calG_{a}$ is equivalent
 to that of $\fuc_{\tilde\mu,\ga}$ over $\calG_{\tilde\mu,a}$,
 whence an autonomous integrand suffices.

\begin{remark}\rm
  The function $h>0$ in Lemma~\ref{L:change} can always be chosen to make the
  measure $\tilde\mu$ finite. Therefore, the finiteness of the underlying measure
  could have been assumed throughout this paper, without any loss of generality.
\end{remark}

\section*{ACKNOWLEDGEMENT}\label{acknow}
\small {This work was supported by the Hungarian National Foundation for Scientific
       Research under Grant K76088, and by Grant Agency of the Czech Republic
       under Grants P202/10/0618 and 201/08/0539.
       Preliminary results on this research were published in \emph{Proceedings ISIT}
       2008~\cite{Csi.Ma.minentfu} and \emph{Proceedings ITW} 2009~\cite{Csi.Ma.minmultfu}.
       An abridged version of the presented results appears in \emph{Proceedings ISIT}
       2012~\cite{Csi.Ma.minentrev}. }

%\makesubmdate
\bigskip
\hfill(\small Received December 3, 2011)

%BIBLIOGRAPHY

%\vspace*{1cm}

\bigskip\medskip

\small
I\,N\,D\,E\,X \\ \label{index}
\rule[2mm]{0.999\textwidth}{0.02mm}\\
\vspace*{-1mm}
\noindent
\begin{minipage}[t]{7cm}
\small
  $\triangleq$~~equal by definition,                    \pageref{equal by definition}
  \\[.1cm]
  $^*$~~convex conjugate,                               \pageref{conjugate of convex function},
                                                        \pageref{conjugate of convex function 2}
  \\[.1cm]
  $\zat{\cdot}{\cdot}$~~inner product in $\R^d$,        \pageref{inner product}
  \\[.1cm]
  $\norm{\cdot}$~~Euclidean norm in $\R^d$,             \pageref{Euclidean norm}
  \\[.1cm]
  $\rightsquigarrow$~~local convergence in measure,     \pageref{local convergence in measure},
                                                        \pageref{local convergence in measure 2}
  \\[.1cm]
  $\nm$~~integrand with $\nm(\ttt,\cdot)\in\varGamma$,  \pageref{integrand with nm}
  \\[.1cm]
  $B$~~class of integrands $\nm$,                       \pageref{the class of integrands}
  \\[.1cm]
  $\ga$~~convex function in class $\varGamma$,          \pageref{strictly convex lsc function in}
  \\[.1cm]
  $\varGamma$~~class of convex functions $\ga$,         \pageref{class of convex functions}
  \\[.1cm]
  $\trn$, $\trnnin$~~Bregman integrand,                 \pageref{Bregman integrand}
   \\[.1cm]
  $\Theta_{\nm}$, $\Theta_\nm^+$~~special subsets of \dom{\dufu_{\nm}},     \pageref{a subset of dom},
                                                                            \pageref{special subset of theta}
  \\[.1cm]
  $[\mu]$ ~~$\mu$-almost everywhere,                    \pageref{almost everywhere}
  \\[.1cm]
  $\lad$, $\ladnin$~~correction integrand,              \pageref{integrand behind the correction},
                                                        \pageref{integrand behind the correction 2}
  \\[.1cm]
  $\moma$~~moment mapping,                              \pageref{moment mapping}
  \\[.1cm]
  $\Crc{F,\nm}$~~integral assigned to face $F$,         \pageref{shift in modified problems 2},
                                                        \pageref{shift in modified problems}
  \\[.1cm]
  asymptotically linear function,                       \pageref{asymptotically linear function},
                                                        \pageref{asymptotically linear function 2}
%  \\[.1cm]
%  biconjugate,                                          \pageref{biconjugate}
  \\[.1cm]
  $\Bres$~~Bregman distance,                            \pageref{Bregman distance}
  \\[.1cm]
  $\Cot_{\nm}$~~correction functional,                  \pageref{correction functional}
    \\[.1cm]
  $\mathsfsl{cl}$~~closure,                             \pageref{closure of}
  \\[.1cm]
  $\cnc{Q}$~~conic core of measure $Q$ on $\R^d$,       \pageref{conic core of a measure}
  \\[.1cm]
  $\cnmoma{\mu}$~~$\moma$-cone of measure $\mu$ on $\TT$,       \pageref{cone of}
  \\[.1cm]
  cofinite function,                                    \pageref{cofinite function},
                                                        \pageref{cofinite function 2}
  \\[.1cm]
  \DCQ~~dual constraint qualification,                  \pageref{dual constraint qualification}
  \\[.1cm]
  $\mathsfsl{dom}$~~effective domain,                   \pageref{effective domain}
\end{minipage}
\begin{minipage}[t]{7cm}\small
  $\exn{\FF_{\nm}}$~~extension of $\FF_{\nm}$,          \pageref{extension of}
  \\[.1cm]
  $\FF_{\nm}$~~family of functions $f_\vte$,            \pageref{family of functions}
    \\[.1cm]
  $f_\vte$~~function in $\FF_{\nm}$,                    \pageref{family of functions}
  \\[.1cm]
  $\Fa{\nm}$~~special family of faces of $\cnmoma{\mu}$,       \pageref{family of faces of}
  \\[.1cm]
  $\calG$~~class of functions with a moment,            \pageref{class of functions with a moment}
  \\[.1cm]
  $\calG_a$~~class of functions with the moment $a$,    \pageref{class of functions with the moment a}
  \\[.1cm]
  $\calG^+$, $\calG^+_a$~~subclasses of nonneg.\ functions,  \pageref{class of nonnegative functions in}
  \\[.1cm]
  $g_a$~~primal solution for $a$,                       \pageref{primal solution for a}
  \\[.1cm]
  $\gea$~~generalized primal solution for $a$,    \pageref{generalized primal solution for a}
  \\[.1cm]
  $\efa$~~effective dual solution for $a$,    \pageref{effective dual solution}
  \\[.1cm]
  generalized Pythagorean identity,  \pageref{generalized Pythagorean identity},
                                     \pageref{generalized Pythagorean identity 2},
                                     \pageref{generalized Pythagorean identity 3}
  \\[.1cm]
  $h_a$~~generalized dual solution for $a$,   \pageref{generalized dual solution for}
  \\[.1cm]
  $\fuc_{\nm}$~~integral functional, \pageref{integral functional}
  \\[.1cm]
  $\vafu_{\nm}$~~value function,     \pageref{nm value function}
%  \\[.1cm]
%  $\ga$-divergence,                  \pageref{ga divergence},
%                                     \pageref{ga divergence 2}
%  \\[.1cm]
%  $\nmh$~~integrand underlying $\Bre{\cdot,h}$,     \pageref{integrand underlying Bre}
  \\[.1cm]
  $\dufu_{\nm}$~~function in dual problem,          \pageref{dual value for}
  \\[.1cm]
  lsc~~lower semicontinuous,                        \pageref{lower semicontinuous, lsc}
  \\[.1cm]
  moment assumption,                                \pageref{moment assumption}
  \\[.1cm]
  \PCQ~~primal constraint qualification,            \pageref{primal constraint qualification}
  \\[.1cm]
  pm~~probability measure,              \pageref{probmeas}
  \\[.1cm]
  primal/dual values,                   \pageref{primal value for a},
                                        \pageref{dual value for}
  \\[.1cm]
  $\mathsfsl{ri}$~~relative interior,   \pageref{relative interior of}
  \\[.1cm]
  $\sg{r}$~~sign of $r$, equals $+$ if $r=0$, \hspace{.3cm}  \pageref{sign of r}
  \\[.1cm]
  $(\TT,\calT)$~~underlying measurable space,                \pageref{underlying measurable space}
  \\[.1cm]
  $\cofi$, $\asli$~~special subsets of $\TT$,                \pageref{special subsets of TT}
\end{minipage}

\bigskip
\normalsize
\makecontacts

\end{document}